\documentclass{amsart}
\usepackage{amsmath,amssymb,amsthm,mathrsfs}

\newcommand{\sA}{{\mathcal A}}
\newcommand{\sB}{{\mathcal B}}
\newcommand{\sC}{{\mathcal C}}
\newcommand{\sD}{{\mathcal D}}

\newcommand{\sF}{{\mathcal F}}

\newcommand{\sL}{{\mathcal L}}
\newcommand{\sM}{{\mathcal M}}
\newcommand{\sN}{{\mathcal N}}
\newcommand{\sP}{{\mathcal P}}
\newcommand{\sR}{{\mathcal R}}
\newcommand{\sT}{{\mathcal T}}
\newcommand{\sU}{{\mathcal U}}
\newcommand{\sV}{{\mathcal V}}
\newcommand{\sW}{{\mathcal W}}
\newcommand{\sX}{{\mathcal X}}
\newcommand{\sY}{{\mathcal Y}}
\newcommand{\sZ}{{\mathcal Z}}

\newcommand{\nn}{\nonumber}

\def\a{{\alpha}}
\def\b{\beta}
\def\c{\gamma}
\def\d{\delta}
\def\de{\Delta}

\def\l{\lambda}

\def\t{\tau}

\def\va{\varphi}

\def\z{\zeta}
\def\ts{\times}

\def\iy{\infty}

\def\ds{\,\textup{d}s\, }

\def\dt{\,\textup{d}t\, }

\def\wt{\widetilde}

\def\BC{{\mathbb C}}

\def\BR{{\mathbb R}}
\def\BT{{\mathbb T}}

\renewcommand{\theequation}{\arabic{section}.\arabic{equation}}

\newcommand{\bpr}{{\noindent\textbf{Proof.}\ \ }}
\newcommand{\epr}{{$\mbox{\ }$ \hfill $\Box$}}

\theoremstyle{plain}
\newtheorem{thm}{Theorem}[section]
\newtheorem{prop}[thm]{Proposition}
\newtheorem{lem}[thm]{Lemma}
\newtheorem{cor}[thm]{Corollary}

\theoremstyle{definition}
\newtheorem{DEF}[thm]{Definition}
\newtheorem{rem}[thm]{Remark}
\newtheorem{exa}[thm]{Example}

\newcommand{\ands}{\quad\mbox{and}\quad}

\newcommand{\supp}{\textup{supp}}

\newcommand{\bL}{\textup{\textbf{L}}}
\newcommand{\bT}{\textup{\textbf{T}}}
\newcommand{\bH}{\textup{\textbf{H}}}
\newcommand{\bI}{\textup{\textbf{I}}}
\newcommand{\bV}{\textup{\textbf{V}}}
\newcommand{\bP}{\textup{\textbf{P}}}
\newcommand{\bR}{\textup{\textbf{R}}}

\newcommand{\sob}{\textup{SB}}
\newcommand{\diam}{\diamond}

\newcommand{\bOm}{\boldsymbol\Omega}

\begin{document}

\title[The twofold Ellis-Gohberg inverse problem]{The twofold Ellis-Gohberg inverse problem in an abstract setting  and applications}

\author[S. ter Horst]{S. ter Horst}
\address{S. ter Horst, Department of Mathematics, Unit for BMI, North West
University,
Potchefstroom, 2531 South Africa}
\email{Sanne.TerHorst@nwu.ac.za}

\author[M.A. Kaashoek]{M.A. Kaashoek}
\address{M.A. Kaashoek, Department of Mathematics,
Vrije Universiteit  Amsterdam,
De Boelelaan 1081a, 1081 HV Amsterdam, The Netherlands}
\email{m.a.kaashoek@vu.nl}

\author[F. van Schagen]{F. van Schagen}
\address{F. van Schagen, Department of Mathematics,
Vrije Universiteit Amsterdam,
De Boelelaan 1081a, 1081 HV Amsterdam, The Netherlands}
\email{f.van.schagen@vu.nl}

\thanks{This work is based on the research supported in part by the
National Research Foundation of South Africa (Grant Numbers 90670 and 93406).}

\begin{abstract}
In this paper we consider a twofold Ellis-Gohberg type inverse problem in an abstract $\ast$-algebraic setting. Under natural assumptions, necessary and sufficient conditions for the existence of a solution are obtained, and it is shown that in case a solution exists, it is unique. The main result relies strongly on an inversion formula for  a $2\ts 2$ block operator matrix  whose off diagonal  entries are Hankel operators while the diagonal  entries are identity operators.  Various special cases are presented, including the cases of  matrix-valued $L^1$-functions on the real line and matrix-valued Wiener functions  on the unit circle of the complex plane. For the latter case, it is shown how the results obtained in an earlier publication by the authors can be recovered.
\end{abstract}



\maketitle


\section{Introduction}\label{sec:intro}

In the present paper we consider a twofold inverse problem related to orthogonal matrix function equations  considered by R.J. Ellis and I. Gohberg for the scalar-valued  case and  mainly in discrete time; see \cite{EG92} and the book \cite{EG03}. The problem is referred to as the twofold EG inverse problem for short. Solutions of the onefold version of the problem, both in  discrete   and continuous time setting, have been obtained in \cite{KvSch14, KvSch16}. For the discrete time setting a solution  of the twofold  problem is given in  \cite{tHKvS1}. One of our aims  is to solve the twofold  problem for the   case of $L^1$-matrix functions on the real line which has not be done yet. More generally,  we will  solve  an abstract $\ast$-algebraic version  of the  twofold EG inverse problem  that  contains various special cases, including the case of $L^1$-matrix functions on the real line. Our abstract setting will include  an abstract  inversion theorem which plays an important role in various concrete cases as well.

The abstract version of   the twofold EG inverse problem we shall be dealing with is presented in Section \ref{sec:setting}.  Here, for  convenience  of the reader,  we  consider  the twofold  EG inverse problem  for $L^1$-matrix functions on the real line, and present the two main theorems for this case,  Theorem \ref{thm:mainL1} and Theorem \ref{thm:mainL2} below. This requires some notation and terminology.

Throughout  $ \BC^{r \times s} $ denotes  the linear space of  all $ r \times s $ matrices with complex entries and $ L^1(\BR)^{r \ts s}$ denotes the space of  all $r \ts s$ matrices of which the entries are Lebesgue integrable functions on  the real line $\BR$. Furthermore
\begin{align*}
& L^1(\BR_+)^{r \ts s}=\{f\in L^1(\BR)^{r \ts s}\mid \supp (f)  \subset \BR_+=[0, \iy)\},\\
& L^1(\BR_-)^{r \ts s}=\{f\in L^1(\BR)^{r \ts s}\mid  \supp (f) \subset  \BR_-=(-\iy,0 ] \}.
\end{align*}
Here $\supp (f)$ indicates the support of the function $f$.  Now assume  we are given
\begin{align}
&a\in L^1(\BR_+)^{p \ts p}, \quad c\in L^1(\BR_-)^{q \ts p},\label{eq:data1}\\
&b\in L^1(\BR_+)^{p \ts q},\quad d\in L^1(\BR_-)^{q \ts q}. \label{eq:data2}
\end{align}
Given these data  the  \emph{twofold EG inverse problem} referred to in the title  is the problem to find $g\in L^1(\BR_+)^{p \ts q}$ satisfying
\begin{align}
&a+g \star c\in L^1(\BR_-)^{p \ts p}, \quad g^*+g^* \star a+c\in L^1(\BR_+)^{q \ts p},\label{eq:probl1}\\
&d+g^* \star b\in L^1(\BR_+)^{q \ts q}, \quad g+g \star d+b\in  L^1(\BR_-)^{p \ts q}. \label{eq:probl2}
\end{align}
Here $g^*(t)=g(-t)^*$ for each $t\in \BR$, and as usual  $f \star h$ denotes the convolution product of  $L^1(\BR)$ matrix functions  $f$ and $h$.

The onefold version of the problem, when only $a$ and $c$ in \eqref{eq:data1} have been given and the problem is to find $g$ such that \eqref{eq:probl1} is satisfied, has been dealt with in \cite{KvSch16}.

To see the EG inverse problem from an operator point of view, let  $g\in L^1(\BR_+)^{p \ts q}$, and  let $G$ and $G_*$  be the Hankel operators  defined by
\begin{align}
&G: L^1(\BR_-)^{q}\to  L^1(\BR_+)^{p}, \quad
(Gf)(t)=\int_{-\iy}^0  g(t-s) f(s)\, \ds, \quad t\geq 0; \label{defG}\\
&G_{*}: L^1(\BR_+)^{p}\to  L^1(\BR_-)^{q}, \ \
(G_{*}h)(t)=\int_{0}^\iy g^*(t-s) h(s)\, \ds, \ \  t\leq 0.  \label{defG*}
\end{align}
Here $L^1(\BR_\pm)^{r}=L^1(\BR_\pm)^{r\ts 1}$. Using these Hankel operators, the conditions  in \eqref{eq:probl1} and \eqref{eq:probl2} are equivalent to
\begin{align}
&a + G c = 0, \qquad G_* a + c = -g^*, \label{eq:probl1a}\\
&d + G_* b = 0, \qquad G d + b = -g. \label{eq:probl2a}
\end{align}
To understand the above identities let us mention  that   we follow the  convention that  an operator   acting  on columns can be extended in a canonical way to an operator  acting on matrices.  We do this without changing  the notation. For instance,  in the first identity  in \eqref{eq:probl1a}   the operator  $G$ acts on each of the $p$ columns of  the $q\ts p$ matrix function $c$, and $Gc$ is the resulting $p\ts p$ matrix function.   Thus  the first condition in \eqref{eq:probl1a} is equivalent to  the first condition in \eqref{eq:probl1}. Similarly,  the  second condition in \eqref{eq:probl1} is equivalent to  the second condition in  \eqref{eq:probl1a} and so on.

Hence the four conditions in \eqref{eq:probl1a} and \eqref{eq:probl2a}  can be summarized by
\begin{equation}\label{eq:probl-v2}
\begin{bmatrix}I&G\\ G_*&I \end{bmatrix}\begin{bmatrix} a\\ c\end{bmatrix}=\begin{bmatrix}0\\-g^* \end{bmatrix}\ands
\begin{bmatrix}I&G\\ G_*&I \end{bmatrix}\begin{bmatrix} b\\ d\end{bmatrix}=\begin{bmatrix}-g\\0 \end{bmatrix}.
\end{equation}
In other words, in this context the inverse problem is to reconstruct, if possible,  a (block) Hankel operator and its associate, from the given data $\{a,b,c,d\}$.

To describe the main theorems in the present context, we need some further preliminaries  about Laurent, Hankel and  Wiener Hopf  operators. Let $\rho$ be a function on $\BR$ given by
\begin{equation}
\rho(t)= r_0 +r(t), t\in \BR,  \hspace{.2cm}
\mbox{where $r\in L^1(\BR)^{k\ts m}$ and $r_0\in \BC^{k\ts m}$}\label{defrho1}.
\end{equation}
With  $\rho$ in \eqref{defrho1} we associate the Laurent operator $L_\rho: L^1(\BR)^m \to L^1(\BR)^k$ which is defined by
\begin{equation} \label{defLrho}
 (L_\rho f)(t)=r_0 f(t)+\int_{-\iy}^\iy r(t-s)f(s) \, ds, \quad  (t\in \BR).
\end{equation}
Furthermore, we write  $L_\rho$ as a $2\ts 2$ operator matrix relative to the direct sum decompositions $L^1(\BR)^\ell= L^1(\BR_-)^\ell\dot{+}  L^1(\BR_+)^\ell$, $\ell= m, k$, as follows:
\[
L_\rho=\begin{bmatrix} T_{-, \rho}& H_{-, \rho}\\[.2cm]  H_{+, \rho} &T_{+, \rho} \end{bmatrix}: \begin{bmatrix}  L^1(\BR_-)^m\\[.2cm]  L^1(\BR_+)^m \end{bmatrix}\to
\begin{bmatrix}  L^1(\BR_-)^k\\[.2cm]  L^1(\BR_+)^k \end{bmatrix}.
\]
Thus $T_{-, \rho}$ and $T_{+, \rho}$ are the Wiener Hopf  operators given by
\begin{align}
\left(T_{-, \rho}f\right)(t)&= r_0f(t)+\int_{-\iy}^0 r(t-s)f(s) \, ds, \quad t\leq 0, \quad  f\in L^1(\BR_-)^m, \label{defTrhoplus} \\
\left(T_{+, \rho}f\right)(t)&= r_0f(t)+\int_{0}^\iy r(t-s)f(s) \, ds, \quad t\geq 0, \quad f\in L^1(\BR_+)^m, \label{defTrhominus}
\end{align}
and $H_{-, \rho}$ and $H_{+, \rho}$ are the Hankel operators given by
\begin{align}
\left(H_{-, \rho}f\right)(t)&=\int_{0}^\iy r(t-s)f(s) \, ds, \quad t\leq 0,\quad  f\in L^1(\BR_+)^m, \label{defHrhoplus}   \\
\left(H_{+, \rho}f\right)(t)&=\int_{-\iy}^0 r(t-s)f(s) \, ds, \quad t\geq 0 \quad  f\in L^1(\BR_-)^m. \label{defHrhominus}
\end{align}
In particular, the Hankel operators $G$ and $G_*$ appearing in \eqref{eq:probl-v2} are equal to $G=H_{+,g}$ and  $G_*=H_{-,g^*}$, respectively.

In what follows, instead of  the  data set $\{a,b,c,d\}$ we will often use the equivalent data set $\{\a, \b, \c,\d\}$,  where
\begin{equation}\label{def:abcd-greek}
\a=e_p +a, \quad \b=b, \quad \c=c, \quad \d=e_q+d.
\end{equation}
Here $e_p$ and $e_q$ are the functions on $\BR$ identically equal to the unit matrix $ I_p$  and $I_q$, respectively.  Using the  data in \eqref{def:abcd-greek} and the definitions of  Toeplitz and Hankel operators in \eqref{defTrhoplus} -- \eqref{defHrhominus}, we define the following operators:
\begin{align}
M_{11} &= T_{+, \a} T_{+, \a^*}- T_{+, \b} T_{+, \b^*}: L^1(\BR_+)^{p} \to L^1(\BR_+)^{p},\label{M11a}\\
M_{21} &= H_{-, \c}T_{+, \a^*} - H_{-, \d} T_{+, \b^*}: L^1(\BR_+)^{p} \to L^1(\BR_-)^{q},\label{M21a}\\
M_{12} &= H_{+, \b}T_{-, \d^*} - H_{+, \a} T_{-, \c^*}: L^1(\BR_-)^{q} \to L^1(\BR_+)^{p},\label{M12a}\\
M_{22}& = T_{-, \d} T_{-, \d^*}- T_{-, \c} T_{-, \c^*}: L^1(\BR_-)^{q} \to L^1(\BR_-)^{q}.\label{M22a}
\end{align}
Notice that these four operators are uniquely determined by the data.

We are now ready to state, in the present context,  our two main theorems. In the abstract setting  these theorems appear in Sections \ref{sec:inversionthm}  and  \ref{sec:solinverseprobl}, respectively. The first is an inversion theorem and the second presents the solution of the EG inverse problem.

\begin{thm}\label{thm:mainL1}
Let $g\in L^1(\BR+)^{p\ts q}$, and let   $W$ be the operator given by
\begin{equation}
W:=\begin{bmatrix} I & H_{+,g}\\[.1cm] H_{-,g^*} & I \end{bmatrix}:
\begin{bmatrix}L^1(\BR_+)^{p}\\[.1cm]   L^1(\BR_-)^{q}\end{bmatrix}\to \begin{bmatrix} L^1(\BR_+)^{p}\\[.1cm]   L^1(\BR_-)^{q}\end{bmatrix}.\label{defZ1}
\end{equation}
Then $W$ is invertible if and only if   $g$ is a solution to a twofold  EG inverse problem for some data set $\{a, b,c,d \}$ as in \eqref{eq:data1} and \eqref{eq:data2}, that is, if and only if  the following two equations are solvable:\begin{equation}\label{eq:abcd-v2}
W\begin{bmatrix} a\\ c\end{bmatrix}=\begin{bmatrix}0\\-g^* \end{bmatrix}\ands
W\begin{bmatrix} b\\ d\end{bmatrix}=\begin{bmatrix}-g\\0 \end{bmatrix}.
\end{equation}
In that case the  inverse of $W$ is given by
\begin{equation}\label{formula-Winv}
W^{-1}=\begin{bmatrix} M_{11} &  M_{12} \\ M_{21} &  M_{22} \end{bmatrix},
\end{equation}
where $M_{ij}$, $1\leq i,j\leq 2$, are the operators defined by \eqref{M11a} -- \eqref{M22a} with $\a$, $\b$, $\c$, $\d$ being given by \eqref{def:abcd-greek} where  $a, b,c,d $ are given by \eqref{eq:abcd-v2}. Furthermore, the operators $ M_{11} $ and $ M_{22} $ are invertible and
\begin{align}
& H_{+,g} = - M_{11}^{-1} M_{12} = - M_{12} M_{22}^{-1} , \label{formulaHg} \\
& H_{-,g^\ast } = - M_{21} M_{11}^{-1} = - M_{22}^{-1} M_{21},  \label{formulaHgast}\\
& g  = - M_{11}^{-1} b  ,\quad g^\ast = - M_{22}^{-1} c. \label{formgg}
\end{align}
\end{thm}

For the second theorem  we need a generalization of the convolution  product  $\star$, which we shall  denote  by the  symbol $\diamond$.  In fact, given the data set $\{a, b,c,d \}$ and  the equivalent data set $\{\a, \b, \c,\d\}$ given by \eqref{def:abcd-greek}, we define the following $\diamond$-products:
\begin{align*}
& \a^* \diamond \a := e_p + a^* + a + a^* \star a, \quad \c^* \diamond \c: = c^* \star c,\\
& \d^* \diamond \d := e_q+ d^* + d + d^* \star  d, \quad \b^* \diamond \b: = b^* \star b,\\
& \a^* \diamond \b: = b + a^* \star b, \quad  \c^* \diamond \d: = c^* + c^* \star d.
\end{align*}

\begin{thm}\label{thm:mainL2}
Let $\{a, b,c,d \}$ be the functions given  by \eqref{eq:data1}
and \eqref{eq:data2}, let $\a$, $\b$, $\c$, $\d$ be the functions given by \eqref{def:abcd-greek}, and let $e_p$ and $e_q$ be the functions on $\BR$ identically equal to the unit matrix $ I_p$  and $I_q$, respectively. Then the twofold  EG inverse problem associated with the data  set  $\{a, b,c,d \}$ has a solution if and only the following two conditions are satisfied:
\begin{itemize}
\item[\textup{(L1)}] $\a^\ast \diamond \a - \c^\ast \diamond \c = e_p $, \ $\d^\ast \diamond \d - \b^\ast \diamond \b = e_q$, \ $\a^\ast \diamond \b = \c^\ast \diamond \d  $; \vspace{.15cm}
\item[\textup{(L2)}]  the operators $ M_{11} $ and $ M_{22} $  defined by  \eqref{M11a} and \eqref{M22a} are one-to-one.
\end{itemize}
In that case $M_{11}$ and $M_{22}$  are invertible, and  the (unique) solution $g$  and its adjoint $g^*$  are given by
\begin{equation} \label{eq:1defgg*}
g = - M_{11}^{-1} b \ands g^* = - M_{22}^{-1} c.
\end{equation}
Here $b$ and $c$ are the matrix functions appearing in   \eqref{eq:data2} and  \eqref{eq:data1}, respectively.
\end{thm}

Assuming that condition (L1) above is satisfied, the invertibility of  the operator $M_{11}$  is equivalent to the injectivity of the operator  $M_{11}$, and the invertibility  $M_{22}$  is equivalent to the injectivity of the operator $M_{22}$. \
To prove these equivalences we use the fact (cf.,  formulas (4.18) and (4.19) in \cite[Section 4.3]{KvSch12})   that $M_{11}$ and $M_{22}$ are also given by
\begin{align}
M_{11}&= I+H_{+, \b}H_{-, \b^*}- H_{+, \a}H_{-, \a^*}, \label{M11b}\\
M_{22}&=I+H_{-, \c}H_{+, \c^*} - H_{-, \d}H_{+, \d^*}.  \label{M22b}
\end{align}
Since the Hankel operators appearing in these formulas are all compact operators,
$M_{11}$ and $M_{22}$ are Fredholm operators, and thus invertible if and only if they are one-to-one.

We shall see in Lemma \ref{lem:altMij}, again assuming that condition (L1) above is satisfied, that the operators $M_{21}$ and $M_{12}$ are also given by
\begin{align}
M_{21}&=T_{-,\d}H_{-, \b^*}-T_{-,\c}H_{-, \a^*},  \label{M21b}\\
M_{12}&=T_{+,\a}H_{+, \c^*}-T_{+,\b}H_{+, \d^*}. \label{M12b}
\end{align}

Since the functions $a$, $b$, $c$, $d$ are $L^1(\BR)$ matrix functions, the operators $M_{ij}$, $1\leq i,j \leq 2$, are   also well-defined as bounded linear operators on the corresponding $L^2$ spaces. It follows that Theorems \ref{thm:mainL1} and   \ref{thm:mainL2} remain true if  the $L^1$ spaces in \eqref{defZ1} are replaced by corresponding $L^2$ spaces. In this $L^2$-setting  Theorems \ref{thm:mainL1} and   \ref{thm:mainL2} are the continuous analogs of Theorems 3.1 and 4.1 in \cite{tHKvS1}. Furthermore, in this $L^2$-setting the adjoints of the operators $M_{ij}$, $1\leq i,j \leq 2$, as operators on $L^2$-spaces, are well-defined as well. In fact, assuming   condition (L1)  is satisfied and  using \eqref{M11a} -- \eqref{M22a} and the identities \eqref{M21b}, \eqref{M12b}, we see that in the $L^2$ setting we have
\begin{equation}
M_{11}^*=M_{11}, \quad M_{21}^*=M_{12},  \quad M_{12}^*=M_{21},  \quad M_{22}^*=M_{22}.
\end{equation}

Theorem  \ref{thm:mainL1} belongs to the wide class of inversion theorems for structured operators. In particular, the theorem can be viewed as  an analogue of the Gohberg-Heinig  inversion theorem for convolution operators on a finite interval \cite{GH75}.
In its present form  Theorem \ref{thm:mainL1} can be seen as an  addition to Theorem 12.2.4  in  \cite{EG03}, where, using a somewhat different notation,  the invertibility of $W$ is proved. The formula for the inverse of $W$ could be obtained from \cite [Theorem~0.1]{GrK05}, where the formula for $M_{11}$ appears in a somewhat different notation. Note that \cite [Theorem 0.1]{GrK05} also solves the asymmetric version of the inversion problem. Formulas \eqref{formulaHg}--\eqref{formgg} seem to be new.

As mentioned before, in  the present paper we put the twofold EG problem in an abstract  $\ast$-algebraic setting.  This allows us to consider  and solve non-stationary  twofold EG problems (see Subsection \ref{subsec:Mat} for an example).  Furthermore,  Theorems  \ref{thm:mainL1} and \ref{thm:mainL2} are obtained as  corollaries of  the two abstract   theorems, Theorem \ref{thm:inversion1}  and  Theorem \ref{thm:mainthm}, derived  in this paper. Also, as we shall prove in Section~\ref{sec:wienerC}, Theorems 3.1 and 4.1 in \cite{tHKvS1} appear as corollaries of  our main theorems.

The paper consists of ten sections (including the present introduction) and an appendix.  In Section \ref{sec:setting} we introduce the  abstract $\ast$-algebraic  setting and state the main problem. Section \ref{examples} presents a numerical example and a number of illustrative special cases, including various Wiener algebra  examples. Sections \ref{sec:preliminaries} and \ref{sec:prelim} have a preliminary character. Here we introduce Toeplitz-like and Hankel-like operators, which play an important role  in the abstract setting, and we derive  a number of identities and lemmas that are used in the proofs of the main results.  In  Section \ref{sec:inversionthm} the abstract inversion theorem  (Theorem \ref{thm:inversion1}) is proved, and in Section  \ref{sec:solinverseprobl}  the solution  to the  abstract  twofold  EG inverse problem (Theorem \ref{thm:mainthm}) is presented  and proved. Theorems \ref{thm:mainL1} and   \ref{thm:mainL2} are proved in Section \ref{sec:wienerL} using the results of Section \ref{sec:inversionthm}  and  Section  \ref{sec:solinverseprobl}. In Section  \ref{sec:can-case} we further specify Theorem \ref{thm:mainthm} for the case when there are additional invertibility  conditions on the underlying data.  The proof in this section is direct and does not use Theorem \ref{thm:mainthm}.  As mentioned in the previous paragraph,   Theorems 3.1 and 4.1 in \cite{tHKvS1} are derived  in Section~\ref{sec:wienerC}  as corollaries of  our main theorems in Sections \ref{sec:inversionthm} and  \ref{sec:solinverseprobl}.

Finally, in Appendix \ref{appendix} we review a number of results that play an important  role in Section \ref{sec:wienerL}, where  we have to relate Hankel-type and Toeplitz-type operators used in Section \ref{sec:inversionthm} and Section  \ref{sec:solinverseprobl} to classical Hankel and Wiener-Hopf integral operators. Appendix \ref{appendix} consists of three subsections.  In Subsection~\ref{ss:prelim} we recall the definition of a Hankel operator  on  $L^2(\BR_+)$ and review some basic facts. In Subsection~\ref{ss:Hankel-L1} we present  a theorem (partially new) characterizing  classical Hankel integral operators  mapping  $L^1(\BR_+)^p$ into $L^1(\BR_+)^q$. Two auxiliary   results are presented  in the final subsection.

\setcounter{equation}{0}
\section{General setting and main problem}\label{sec:setting}
We first describe the general  $\ast$-algebraic  setting  that  we will be working with. To do this we  use the notation introduced on pages 109 and 110  of  \cite{KvSch13};  see also the first two pages of \cite[Section II.1]{GKW89}.  Throughout $\sA$, $\sB$, $\sC$ and $\sD$ are complex linear vector spaces such that the following set of $2 \times 2$ block matrices form an algebra:
\begin{equation}
\label{defR1}
\sM= \sM_{\sA, \sB, \sC, \sD} = \left\{f=\begin{bmatrix}
 a     &   b \\
   c   &  d
\end{bmatrix}\mid a\in \sA, \    b\in \sB, \   c\in \sC, \   d\in \sD\right\}.
\end{equation}
Furthermore,  we assume  $\sA$ and $\sD$ are $*$-algebras (see \cite[Chapter IV]{Rickart60} for the definition of this notion) with units $e_{\sA}$  and $e_{\sD}$, respectively, and endowed with involutions~${}^*$. The diagonal
\[
e_\sM = \begin{bmatrix} e_\sA & 0 \\ 0 & e_\sD \end{bmatrix}
\]
is the unit element of  $ \sM$. Moreover, $\sC$ is a linear space isomorphic to $\sB$ via
a conjugate linear transformation~${}^*$ whose inverse is also denoted by ${}^*$.  We require $ \sM $ to be a $\ast$-algebra with respect to  the usual matrix multiplication and with the involution   given by
\begin{equation*}
\begin{bmatrix} a & b \\ c & d \end{bmatrix}^\ast =
\begin{bmatrix} a^\ast & c^\ast \\ b^\ast & d^\ast \end{bmatrix}.
\end{equation*}
The algebras $\sA$ and $\sD$ are assumed to admit direct sum decompositions:
\begin{equation}
\label{decomp1}
\sA=\sA_-^0\dot{+}\sA_d \dot{+}\sA_+^0 , \quad
\sD=\sD_-^0 \dot{+}\sD_d \dot{+}\sD_+^0.
\end{equation}
In these two direct sum  decompositions  the summands are assumed to be subalgebras
of $\sA$ and $\sD$, respectively.
Furthermore, we require
\begin{equation}\label{requir1}
\begin{array}{ccc}
e_\sA\in \sA_d ,  & (\sA_-^0)^*=\sA_+^0 , & (\sA_d)^*=\sA_d,\\[.2cm]
e_\sD\in \sD_d , & (\sD_-^0)^*=\sD_+^0 , &(\sD_d)^*=\sD_d.
\end{array}
\end{equation}
Set
\[
\sA_-=\sA_-^0\dot{+}\sA_d,\quad
\sA_+=\sA_d \dot{+}\sA_+^0,\quad
\sD_-=\sD_-^0 \dot{+}\sD_d,\quad
\sD_+=\sD_d \dot{+}\sD_+^0.
\]
We also assume that $\sB$ and $\sC$ admit direct sum decompositions:
\begin{equation}\label{decomp2}
\sB=\sB_-\dot{+}\sB_+,\quad \sC=\sC_-\dot{+}\sC_+, \quad\mbox{such that}\quad  \sC_-=\sB_+^*, \quad  \sC_+=\sB_-^*.
\end{equation}
These direct sum decompositions yield a  direct sum decomposition of $\sM$, namely
$\sM = \sM_-^0 \dot{+} \sM_d \dot{+} \sM_+^0$, where
\begin{equation}\label{decompM}
\sM_-^0 = \begin{bmatrix} \sA_-^0 & \sB_- \\ \sC_- & \sD_-^0 \end{bmatrix}, \quad
 \sM_d = \begin{bmatrix} \sA_d & 0  \\ 0 & \sD_d \end{bmatrix}, \quad
\sM_+^0 = \begin{bmatrix} \sA_+^0 & \sB_+ \\ \sC_+ & \sD_+^0 \end{bmatrix}.
\end{equation}
Note that
\[
(\sM_-^0)^*=\sM_+^0, \quad  (\sM_+^0)^*=\sM_-^0, \quad  \sM_d ^*= \sM_d.
\]
Finally, we assume that the products of elements from the summands in $\sM= \sM_-^0 \dot{+}  \sM_d \dot{+} \sM_{+}^0$ satisfy the rules of the  following
\begin{center}
\textbf{Multiplication table}

\smallskip
\begin{tabular}{l||l|l|l}
$\ \ts$ & $\sM_-^0$ & $\sM_d $ & $\sM_+^0$ \\
\hline
\hline
 $\sM_-^0$& $\sM_-^0$  & $\sM_-^0$ & $\sM$ \\
\hline
 $\sM_d$  & $\sM_-^0$ &   $ \sM_d $&$\sM_+^0$ \\
\hline
 $\sM_+^0$  & $ \sM $  &  $\sM_+^0$&$\sM_+^0$
\end{tabular}
\end{center}
We say  that   the algebra $\sM=\sM_{\sA, \sB, \sC, \sD}$ defined by (2.1)  is \emph{admissible} if all the conditions listed in the above paragraph are satisfied.

\paragraph{Main problem.} We are now ready to state the main problem  that we shall be dealing with.  Let $\a \in \sA_+ $, $\b \in \sB_+ $, $\c \in \sC_-$ and $\d \in \sD_- $ be given.
We call  $g \in \sB_+ $ a \emph{solution to the  twofold  EG inverse problem associated with  $ \a $, $\b$, $\c$, and $\d$} whenever
\begin{align}
\a + g \c  -e_\sA \in  \sA_-^0   & \ands g^\ast \a  + \c \in \sC_+ , \label{inclu12} \\
g \d +\b \in \sB_- & \ands  \d +g^\ast \b   - e_\sD \in \sD_+^0. \label{inclu34}
\end{align}
Our main aim is to determine  necessary and sufficient  conditions for this inverse problem to be solvable  and to derive explicit formulas for its solution.
We shall show that the solution, if it exists, is unique. The following result is a special case of \cite[Theorem 1.2]{KvSch13}.

\begin{prop}\label{prop:cond123}
If the twofold  EG inverse problem associated with  $ \a $, $\b$, $\c$ and $\d$ has a solution, then
\begin{itemize}
\item[\textup{(C1)}] \qquad\qquad  $\a^\ast \a - \c^\ast \c = P_{\sA_d} \a $,
\item[\textup{(C2)}] \qquad\qquad  $\d^\ast \d - \b^\ast \b = P_{\sD_d} \d $,
\item[\textup{(C3)}] \qquad\qquad  $\a^\ast \b = \c^\ast \d  $.
\end{itemize}
Here $ P_{\sA_d} $ and $ P_{\sD_d}$  denote  the projections of  $ \sA$ and $\sD $  onto $\sA_d $ and $ \sD_d$, respectively, along $\sA^0=\sA_-^0\dot{+} \sA_+^0$ and  $\sD^0=\sD_-^0\dot{+} \sD_+^0$, respectively.
\end{prop}

Notice that (C1) and (C2) imply that
\begin{equation}\label{selfadj-ad0}
{a_0 :=} P_{\sA_d} \a= (P_{\sA_d} \a)^\ast =a_0^* \ands  d_0 := P_{\sD_d} \d = (P_{\sD_d} \d )^\ast = d_0^*.
\end{equation}
Furthermore, together the  three conditions (C1)--(C3) are equivalent to
\begin{equation}\label{C1-3matrix}
\begin{bmatrix} \a^\ast & \c^\ast  \\ \b^\ast & \d^\ast \end{bmatrix}
\begin{bmatrix} e_\sA & 0 \\ 0 & -e_\sD \end{bmatrix}
\begin{bmatrix} \a  & \b \\ \c & \d \end{bmatrix}=
\begin{bmatrix} a_0 & 0 \\ 0 & - d_0 \end{bmatrix}.
\end{equation}

\begin{rem}\label{rem:inv0}
Since $a_0$ and $d_0$ belong to $\sA_d$ and $\sD_d$, respectively, invertibility of $a_0$ in $\sA$ and of $d_0$ in $\sD$ imply that  $ a_0^{-1} \in \sA_d $ and $ d_0^{-1} \in \sD_d $. In other words, $a_0 $ and $d_0 $ are invertible in  $\sA_d $ and $\sD_d $,  respectively.
\end{rem}

\begin{rem}\label{rem:cond456}
In the sequel it will often be assumed that $a_0$ and $d_0$ are invertible.  In that case    the following three conditions are well defined.
\begin{itemize}
\item[\textup{(C4)}] \qquad\qquad  $\a a_0^{-1} \a^\ast - \b d_0^{-1} \b^\ast  = e_\sA $,
\item[\textup{(C5)}] \qquad\qquad  $\d d_0^{-1} \d^\ast - \c a_0^{-1} \c^\ast = e_\sD $,
\item[\textup{(C6)}] \qquad\qquad  $ \a a_0^{-1} \c^\ast = \b d_0^{-1} \d^\ast $.
\end{itemize}
In solving the twofold  EG inverse problem referred to above we shall always assume that
$a_0$ and $d_0$ are invertible and that the six conditions (C1)--(C6) are fulfilled.
\end{rem}

The next lemma shows that in many cases (C4)--(C6) are satisfied whenever
conditions (C1)--(C3) are satisfied.

\begin{lem}\label{lem:C1-C6}
Let $\a \in \sA_+ $, $\b \in \sB_+ $, $\c \in \sC_-$, $\d \in \sD_- $,  and let
\[
Q = \begin{bmatrix} \a & \b \\ \c & \d \end{bmatrix}.
\]
Assume that $a_0 $ and $ \a $ are invertible in $ \sA $,
and that $d_0 $ and $ \d $ are invertible in $ \sD $.
If, in addition,  $ \a $, $ \b $, $ \c $ and $ \d $ satisfy conditions \textup{(C1)--(C3)},
then  $ Q $  is invertible, and conditions \textup{(C4)--(C6)} are satisfied.
\end{lem}

\bpr
Since $\d$ is invertible, a classical Schur complement argument
(see, e.g., formula (2.3) in \cite[Chapter 2]{BGKR08}) shows that
\[
Q  = \begin{bmatrix} \a  & \b  \\ \c &  \d  \end{bmatrix} =
\begin{bmatrix} e_\sA & \b\d^{-1}  \\ 0&  e_\sD  \end{bmatrix}
\begin{bmatrix} \de & 0  \\ 0&  \d  \end{bmatrix}
\begin{bmatrix} e_\sA & 0 \\ \c\d^{-1}  & e_\sD  \end{bmatrix}, \quad \mbox{with } \de=\a -\b \d^{-1}\c.
\]
Using  the invertibility of $\a$ and $\d$, we can rewrite (C3) as $\b\d^{-1}=\a^{-*}\c^*$. The latter identity together with (C1) yields:
\[
\de=\a -\b \d^{-1}\c=\a - \a^{-*}\c^*\c=\a^{-*}\left(\a^*\a-\c^*\c\right)=\a^{-*}a_0.
\]
It follows that the Schur complement $\de$ is invertible. But then $ Q $  is invertible too,
and the identity \eqref{C1-3matrix} shows that the inverse  $Q^{-1}$ of $Q$ is given by
\[
Q^{-1} =
\begin{bmatrix} a_0^{-1} & 0 \\ 0 & - d_0^{-1} \end{bmatrix}
\begin{bmatrix} \a^\ast & \c^\ast  \\ \b^\ast & \d^\ast \end{bmatrix}
\begin{bmatrix} e_\sA & 0 \\ 0 & -e_\sD \end{bmatrix}.
\]
Since  $Q Q^{-1}$  is a $2 \ts 2$  block identity matrix, we conclude  that
\[
\begin{bmatrix} \a  & \b  \\ \c &  \d  \end{bmatrix}
\begin{bmatrix} a_0^{-1} & 0 \\ 0 & - d_0^{-1} \end{bmatrix}
\begin{bmatrix} \a^\ast & \c^\ast  \\ \b^\ast & \d^\ast \end{bmatrix}
\begin{bmatrix} e_\sA & 0 \\ 0 & -e_\sD \end{bmatrix}=
\begin{bmatrix} e_\sA & 0 \\ 0 & e_\sD \end{bmatrix}
\]
This yields
\begin{equation}\label{C1-3matrix*}
\begin{bmatrix} \a  & \b  \\ \c &  \d  \end{bmatrix}
\begin{bmatrix} a_0^{-1}  & 0 \\ 0 & - d_0^{-1} \end{bmatrix}
\begin{bmatrix} \a^\ast & \c^\ast  \\ \b^\ast & \d^\ast \end{bmatrix}
= \begin{bmatrix} e_\sA & 0 \\ 0 & - e_\sD \end{bmatrix},
\end{equation}
and hence (C4)--(C6) are satisfied.
\epr

\setcounter{equation}{0}
\section{A numerical example and  some illustrative  special cases}\label{examples}

In this section we present a few inverse problems which are  special cases of the abstract problem presented in the previous section.

\subsection{A numerical example}\label{ex:simple1}
As a first illustration  we  consider a simple example of a  problem for $3\ts 3$ matrices. Given
\begin{align}
& \a =  \textstyle { \frac{1}{8} }
\begin{bmatrix} -2 & 2 & 0 \\ 0 & -3 & -4 \\ 0 & 0 & 6 \end{bmatrix},  \qquad
\b = -  \textstyle { \frac{1}{8} }
\begin{bmatrix} -2 & -6 & 0 \\ 0 & 1 & -4 \\ 0 & 0 & -2 \end{bmatrix}, \label{defab} \\ \noalign{\medskip}
& \c = -  \textstyle { \frac{1}{8} }
\begin{bmatrix} -2 & 0 & 0 \\ -4 & 1 & 0 \\ 0 & -6 & -2 \end{bmatrix}, \qquad
\d =  \textstyle  \frac{1}{8}
\begin{bmatrix} 6 & 0 & 0 \\ -4 & -3 & 0 \\ 0 & 2 & -2 \end{bmatrix}
,\label{defcd}
\end{align}
we seek a $3 \ts 3 $ upper triangular matrix $ g $  such that
\begin{align}
& \a + g \c = \begin{bmatrix} 1 & 0 & 0 \\ \star & 1 & 0 \\ \star & \star & 1 \end{bmatrix}, \quad
g^* \a + \c = \begin{bmatrix} 0 & \star & \star \\ 0 & 0 & \star \\ 0 & 0 & 0 \end{bmatrix}, \label{inclu12b} \\
& \b + g \d = \begin{bmatrix} 0 & 0 & 0 \\ \star & 0 & 0 \\ \star & \star & 0 \end{bmatrix}, \quad
g^* \b + \d = \begin{bmatrix} 1 & \star & \star \\ 0 &1 & \star \\ 0 & 0 & 1 \end{bmatrix}.
\label{inclu34b}
\end{align}
Here the symbols $ \star $ denote unspecified entries.
By direct checking it is easy to see that the matrix $ g_\circ$ given by
\begin{equation*}
g_\circ =  \begin{bmatrix} 1 & 2 & 0 \\ 0 & 1 & 2 \\ 0 & 0 & 1  \end{bmatrix}
\end{equation*}
is upper triangular and satisfies \eqref{inclu12b} and \eqref{inclu34b}. From the general
results about  existence of solutions  and methods to determine solutions, which will be
presented  in this paper, it   follows that $g_\circ$ is the only  solution.
For  this example it also straightforward to check that conditions (C1)--(C3) presented\
in the previous section are satisfied.

\subsection{A class of finite dimensional matrix examples}\label{subsec:Mat}
We will put the problem considered in the preceding example into the general setting
considered in the previous section. Let $ p \geq 1 $ be an integer (in the example above we took $ p = 3 $),  and let
\begin{equation}\label{ExampleMatrices2}
\sA = \sB = \sC=\sD=\BC^{ p \times p }.
\end{equation}
The involution $ {}^\ast $ is given by the usual transposed conjugate  of a matrix.
Let $ \sA_+^0 = \sB_+^0 = \sC_+^0 = \sD_+^0 $ be the subspace of $ \BC^{ p \times p } $
of the strictly upper triangular matrices and $ \sA_-^0 = \sB_-^0 = \sC_-^0 = \sD_-^0 $
the subspace of the strictly lower triangular matrices. Furthermore, let
$ \sA_d = \sB_d = \sC_d = \sD_d $ be the subspace consisting of the $p\ts p$ diagonal matrices,  that is,  matrices with all   entries off the  main diagonal being equal to zero.
We set
\begin{align}
& \sA_- = \sA_-^0 \dot{+} \sA_d, \quad
\sA_+ = \sA_d \dot{+} \sA_+^0, \quad
\sD_- = \sD_-^0 \dot{+} \sD_d, \quad
\sD_+ = \sD_d \dot{+} \sD_+^0, \label{def:ADplusmin}\\
&\hspace{1cm} \sB_- = \sB_-^0, \quad
\sB_+ = \sB_d \dot{+} \sB_+^0, \quad
\sC_- = \sC_-^0 \dot{+} \sC_d, \quad
\sC_+ = \sC_+^0 .  \label{def:BCplusmin}
\end{align}
The  problem we consider in this setting  is the following.
Let $\a$, $\b$,  $\c$, $\d$  be given $p \ts p$ matrices,
and assume that $\a \in \sA_+ $, $ \b \in \sB_+ $, $ \c \in \sC_- $ and $ \d \in \sD_-  $.
Then a  $p \ts p$ matrix $ g \in \sB_+ $ is said to be a solution to
the EG inverse problem for the given data $\a$, $\b$,  $\c$, $\d$
whenever the four inclusions in \eqref{inclu12} and \eqref{inclu34} are satisfied.
In the numerical example considered above  this amounts to the conditions
\eqref{inclu12b} and \eqref{inclu34b} being fulfilled.

If a solution exists, then the conditions (C1)--(C3) are satisfied, what in this setting means that
\begin{equation}\label{C1-3mMatrix}
\begin{bmatrix} \a^\ast & \c^\ast \\ \b^\ast & \c^\ast \end{bmatrix}
\begin{bmatrix} I_p & 0 \\ 0 & -I_p \end{bmatrix}
\begin{bmatrix} \a & \b \\ \c & \d \end{bmatrix}  =
\begin{bmatrix} \a_d & 0 \\ 0 & -\d_d \end{bmatrix} .
\end{equation}
Here $I_p$  is the $p\ts p$ identity matrix, $ \a_d $ is the diagonal matrix whose
main  diagonal coincides with the one of $ \a $, and $ \d_d $ is the diagonal matrix whose
main  diagonal coincides with the one of $ \d $.
If $ \a_d $ and $ \d_d $ are invertible, then $ \a $ and $ \d $ are invertible as
$ p \times p $ (lower or upper) triangular matrices.
In this case, as we shall see in Theorem \ref{thm:cancase}, the twofold EG inverse problem
is solvable, the solution is unique, and the solution is given by
$ g = - P_{\sB_+} \left( ( \a^\ast )^{-1} \c^\ast  \right)$.

The above special case is an example of a non-stationary EG inverse problem. We intend to deal with other non-stationary problems in a later publication, using elements of \cite{GKW89}; see also \cite[Section 5]{KL12}.

\subsection{Wiener algebra examples}\label{Wieneralgebra}
Let $\sN$ be a unital  $*$-algebra with unit $e_{\sN}$ and  involution ${}^*$.
We assume that  $\sN$ admits  a direct sum decomposition:
\[
\sN = \sN_{-,0} \dot{+} \sN_d \dot{+} \sN_{+,0}.
\]
In this direct sum  decomposition the summands are subalgebras of $ \sN $, and we require
\begin{align*}
&e_\sN \in \sN_d, \quad (\sN_d)^* = \sN_d, \quad  (\sN_{-,0} )^\ast = \sN_{+,0},\\
\noalign{\smallskip}
&  \quad
 \sN_d \sN_{\pm,0} \subset \sN_{\pm,0}, \quad \sN_{\pm,0} \sN_d \subset \sN_{\pm,0}.
\end{align*}
Given $\sN$ we construct two admissible algebras  $\sM_{\sA, \sB, \sC, \sD}$
 using the following two translation tables:
\begin{center}
{\bf Table 1} \\
\medskip
\begin{tabular}{|l|l|l|l|l|l|l|}
\hline \noalign{\smallskip}
 $\sA$ & $ \sA_+^0 $ & $ \sA_d $ & $ \sA_-^0 $ & $\sB$ & $\sB_+ $ & $\sB_- $ \\
\hline \noalign{\smallskip}
 $ \sN^{p \times p} $ \qquad & $ \sN^{p \times p}_{+,0} $ &
 $ \sN^{p \times p}_d  $  & $ \sN^{p \times p}_{-,0}$ &  $ \sN^{p \times q}$ &
$ \sN^{p \times q}_+ $ & $ \sN^{p \times q}_{-,0} $ \\
\hline \hline \noalign{\smallskip}
 $\sC$ & $\sC_- $ & $\sC_+ $ & $\sD$  & $ \sD_+^0 $ & $ \sD_d $ & $ \sD_-^0 $   \\
\hline \noalign{\smallskip}
$ \sN^{q \times p} $ & $ \sN^{q \times p}_- $ & $ \sN^{q \times p}_{+,0} $ &
$ \sN^{q \times q} $ \quad & $ \sN^{q \times q}_{+,0}  $ & $ \sN^{q \times q}_d  $ &
$ \sN^{q \times q}_{-,0}  $  \\
\hline
\end{tabular}
\end{center}

\bigskip
\begin{center}
{\bf Table 2}  \\
\medskip
\begin{tabular}{|l|l|l|l|l|l|l|}
\hline
\noalign{\smallskip}
$\sA$ & $ \sA_+^0 $ & $ \sA_d $ & $ \sA_-^0 $ & $\sB$ & $\sB_+ $ & $\sB_- $ \\
\hline \noalign{\smallskip}
 $ \sN^{p \times p} $ \qquad & $ \sN^{p \times p}_{+,0} $ &
$ \sN^{p \times p}_d $ & $ \sN^{p \times p}_{-,0} $ & $ \sN^{p \times q}_0 $ &
$ \sN^{p \times q}_{+,0} $ & $ \sN^{p \times q}_{-,0} $ \\
\hline\hline \noalign{\smallskip}
$\sC$ & $\sC_- $ & $\sC_+ $ & $\sD$  & $ \sD_+^0 $ & $ \sD_d $ & $ \sD_-^0 $   \\
\hline \noalign{\smallskip}
 $ \sN^{q \times p}_0 $ & $ \sN^{q \times p}_{-,0}$ &  $ \sN^{q \times p}_{+,0} $ &
$ \sN^{q \times q} $ \quad  & $ \sN^{q \times q}_{+,0} $ & $ \sN^{q \times q}_d  $ &
$ \sN^{q \times q}_{-,0} $  \\
\hline
\end{tabular}
\end{center}
\medskip
\noindent
In  subsequent special cases  we make these examples more concrete.

\medskip
\subsubsection{The Wiener algebra on the real line.}\label{ssec:WR}
Recall that the Wiener algebra on the real line $ \sW(\BR) $ consists of the functions
$ \varphi $ of the form
\begin{equation}\label{def:WR}
\varphi(\lambda) = f_0 + \int_{-\infty}^\infty e^{i \lambda t} f(t) \dt,\quad \l\in \BR,
\end{equation}
with $ f_0 \in\BC $ and $ f \in L^1(\BR) $.
The subspaces $ \sW(\BR)_{\pm,0} $ consist of the functions $\varphi$ in $ \sW(\BR) $ for which in the  representation \eqref{def:WR} the constant $ f_0 = 0 $ and $ f \in L^1(\BR_\pm) $.
A function $ \varphi $ belongs to the subspace $ \sW(\BR)_0 $  if and only if
$ f =0 $ in the representation in \eqref{def:WR}.
With $ \sN = \sW(\BR) $, it is straightforward  to check that the spaces $\sA, \sB, \sC$, and $\sD$
and their subspaces defined by Table~2 have all the properties listed in the first
two paragraphs of  Section \ref{sec:setting}, that is, $\sM=\sM_{\sA, \sB, \sC, \sD}$ is
admissible. Indeed, let
\begin{equation*}
\a \in e_p + \sW(\BR)_{+,0}^{p \times p} , \quad
\b \in \sW(\BR)_{+,0}^{p \times q} , \quad
\c \in \sW(\BR)_{-,0}^{q \times p} , \quad
\d \in e_q + \sW(\BR)_{-,0}^{q \times q},
\end{equation*}
where  $e_p$  and $e_q$ are the functions identically equal to the unit matrix $ I_p $ and $I_q $, respectively.
Then the twofold EG inverse problem is to find $g \in \sW(\BR)_{+,0}^{p\ts q}$
such that the following four inclusions are satisfied:
\begin{align*}
& \a + g \c  - e_p \in  \sW(\BR)_{-,0}^{ p \times p }  \ands
g^\ast \a  + \c    \in \sW(\BR)_{+ ,0}^{ q \times p };  \\
&  g\d+\b  \in \sW(\BR)_{-,0}^{p\times q}  \ands
\d+ g^*\b   -e_q \in \sW(\BR)_{+,0}^{q\ts q }.
\end{align*}
Notice that these inclusions are just the same as the inclusions in
\eqref{inclu12} and  \eqref{inclu34}.
In this way this twofold EG inverse problem is put in the abstract setting of the
twofold EG inverse problem defined in Section \ref{sec:setting}.

\begin{rem}
The version of the twofold EG inverse problem considered in this subsubsection is isomorphic to
the twofold EG inverse problem considered in the introduction.
This follows from the definition of the Wiener algebra $ \sW(\BR)$ in \eqref{def:WR}.
The solution of the twofold EG inverse problem as described in this subsubsection
 follows from Theorems \ref{thm:mainL1}  and \ref{thm:mainL2}.
The latter two theorems will be proved in Section \ref{sec:wienerL}.
\end{rem}

Note that in this special case the algebra $\sM=\sM_{\sA, \sB, \sC, \sD}$ appearing in
\eqref{defR1} can be considered as a subalgebra of
$ \sW(\BR)^{ (p+q) \times (p+q)} $.  Indeed,
\begin{equation*}
\sM =  \sW(\BR)^{(p+q) \times (p+q)}_{-,0} \dot{+} \sM_d \dot{+}
 \sW(\BR)^{(p+q) \times (p+q)}_{+,0}
\end{equation*}
with
\[
\sM_d  = \left\{ \begin{bmatrix} a_0 & 0 \\ 0 & d_0 \end{bmatrix} \mid a_0 \in \BC^{p \times p} ,
 d_0 \in \BC^{q \times q} \right\} .
\]

\medskip
\paragraph{The case $ \sN = \sR\sW(\BR) $.}
Let $ \sR\sW(\BR) $ be the subalgebra of $ \sW(\BR)$ consisting of all rational
functions in $ \sW(\BR) $.
With $ \sN = \sR\sW(\BR) $ it is straightforward  to check that the resulting
$\sA, \sB, \sC, \sD$ defined in Table 2 have all the properties listed in
the first two paragraphs of Section \ref{sec:setting}, that is,
$\sM=\sM_{\sA, \sB, \sC, \sD}$ is admissible. Let
\begin{align*}
&\a \in e_p + \sR\sW(\BR)^{p \times p}_{+,0}, \quad  \b \in \sR\sW(\BR)^{p \times q}_{+,0},\quad \c \in \sR\sW(\BR)^{q \times p}_{-,0},\\
&\hspace{2.5cm} \d \in e_q + \sR\sW(\BR)^{q \times q}_{-,0}.
\end{align*}
The twofold EG inverse problem is to find $g \in \sR\sW(\BR)_{+,0}^{p\ts q}$
such that the following four inclusions are satisfied:
\begin{align*}
& \a + g \c  -e_p \in  \sR\sW(\BR)_{-,0}^{ p \times p },   \ands
g^\ast \a  + \c    \in \sR\sW(\BR)_{+ ,0}^{ q \times p };  \\
&   g\d+\b  \in \sR\sW(\BR)_{-,0}^{p\ts q} \ands
\d+ g^*\b   -e_q\in \sR\sW(\BR)_{+,0}^{q\ts q} .
\end{align*}
In a forthcoming paper we plan to deal with the twofold EG inverse problem for rational functions in $ \sW(\BR) $, using minimal realizations of the rational functions involved and related state space techniques. The latter will lead to new explicit formulas for the solution.

\medskip
\paragraph{The case $ \sN = \sF\sW(\BR)$.} Let $\sF\sW(\BR)$ denote the subalgebra of $\sW(\BR)$ of functions in $\sW(\BR)$ whose inverse Fourier transforms are elements in $L^1(\BR)$ with finite support. Hence if $\rho\in\sW(\BR)$ is given by
\[
\rho(\l)= r_0 +\int_{-\iy}^\iy e^{i \l t} r(t)\,dt\quad (\l\in \BR),
\]
with $r\in L^1(\BR)$ and $r_0\in\BC$, then $\rho\in \sF\sW(\BR)$ in case there are real numbers $\tau_1<\tau_2$ so that $r(t)=0$ for all $t\not\in [\tau_1,\tau_2]$.
In this case, one easily verifies that  $\sM=\sM_{\sA, \sB, \sC, \sD}$ with $\sA, \sB, \sC, \sD$ as in Table 2 is admissible.
The twofold EG inverse problem specified for these choices can be stated as follows.
Let
\begin{align*}
&\a \in e_p + \sF\sW(\BR)^{p \times p}_{+,0}, \quad  \b \in \sF\sW(\BR)^{p \times q}_{+,0}, \quad  \c \in \sF\sW(\BR)^{q \times p}_{-,0}, \\
&\hspace{2.5cm}\d \in e_q + \sF\sW(\BR)^{q \times q}_{-,0}.
\end{align*}
The twofold EG inverse problem now is to find $g \in \sR\sW(\BR)_{+,0}^{p\ts q}$
such that the following four inclusions are satisfied:
\begin{align*}
&\a + g \c  -e_p \in  \sF\sW(\BR)_{-,0}^{ p \times p }  \ands
g^\ast \a  + \c    \in \sF\sW(\BR)_{+ ,0}^{ q \times p };  \\
&  g\d+\b  \in \sF\sW(\BR)_{-,0}^{p\ts q} \ands
 \d+ g^*\b   -e_q\in \sF\sW(\BR)_{+,0}^{q\ts q}.
\end{align*}
We plan to return to this case in a forthcoming paper.

\subsubsection{The Wiener algebra on the unit circle.}\label{ssec:WT}
Let $ \sN = \sW(\BT) $, where
 $ \sW(\BT) $ is the Wiener algebra of functions on the unit circle  $ \BT $, that
is,  the algebra of  all functions on $\BT $ with absolutely converging Fourier series.
Define $\sA, \sB, \sC, \sD$ as in Table~1. In this case $\sM=\sM_{\sA, \sB, \sC, \sD}$ is admissible too. Note that the Fourier transform defines an isomorphism  between $ \sW(\BT) $ and the algebra $ \ell^1 $ of absolutely converging complex sequences. The version of the twofold EG inverse problem for $ \ell^1 $ has been  solved in  \cite{tHKvS1}.
In Section \ref{sec:wienerC} we give a new proof of the main theorems in \cite{tHKvS1} by
putting the inversion theorem, \cite[Theorem~3.1]{tHKvS1}, and the solution of the
twofold EG inverse problem, \cite[Theorem~4.1]{tHKvS1}, into the general setting
of Section \ref{sec:setting} and using the results  of  Sections  \ref{sec:preliminaries}--\ref{sec:solinverseprobl}.

\medskip
\paragraph{The case $ \sN = \sR\sW(\BT) $.}
Let $ \sR\sW(\BT) $ be the subalgebra of $ \sW(\BT)$ consisting of all rational
functions in $\sW(\BT)$.
With $ \sN = \sR\sW(\BT)$ and $\sA, \sB, \sC, \sD$ as in Table~1, the resulting algebra
$\sM=\sM_{\sA, \sB, \sC, \sD}$ is admissible. Let
\begin{align*}
&\a \in  \sR\sW(\BT)^{p \times p}_{+}, \quad \b \in \sR\sW(\BT)^{p \times q}_{+}, \quad
\c \in \sR\sW(\BT)^{q \times p}_{-}, \\
&\hspace{2cm}\d \in \sR\sW(\BT)^{q \times q}_{-}.
\end{align*}
The twofold EG inverse problem is to find $g \in \sR\sW(\BT)_{+}^{p\ts q}$ such that the following
four inclusions are satisfied:
\begin{align*}
& \a + g \c  -e_p \in  \sR\sW(\BT)_{-,0}^{ p \times p }  \ands
g^\ast \a  + \c    \in \sR\sW(\BT)_{+ ,0}^{ q \times p };  \\
&  g\d+\b  \in \sR\sW(\BT)_{-,0}^{p\ts q} \ands
 \d+ g^*\b   -e_q\in \sR\sW(\BT)_{+,0}^{q\ts q}.
\end{align*}
The onefold EG inverse problem for rational matrix functions on $ \BT $ is treated in
\cite[Section 6]{KvSch16}. Minimal realizations of the functions involved play an important
role in the approach in \cite{KvSch16}. We intend to work on the twofold EG inverse problem for rational matrix functions on the unit circle in a later publication, again  using  minimal state space realizations of the functions involved.

\medskip
\paragraph{The case $ \sN = \sT\sP$.} Let $\sT\sP$ be the set consisting of the trigonometric polynomials in $z$ viewed as a subalgebra of $\sW(\BT)$, and write $\sT\sP_+$ and $\sT\sP_-$ for the subalgebras of polynomials in $z$ and in $z^{-1}$, respectively. With  $\sT\sP_{+,0}$ and $\sT\sP_{-,0}$ we denote the corresponding spaces with the constant functions left out. Again, $\sA, \sB, \sC, \sD$ are defined as in Table 1, and the algebra $\sM=\sM_{\sA, \sB, \sC, \sD}$ is admissible. Let
\begin{equation*}
\a \in  \sT\sP^{p \times p}_{+} , \quad
\b \in \sT\sP^{p \times q}_{+} , \quad
\c \in \sT\sP^{q \times p}_{-} , \quad
\d \in \sT\sP^{q \times q}_{-}.
\end{equation*}
In this context the twofold EG inverse problem is to find $g \in \sT\sP^{p \times q}_{+} $ such that the following
 four inclusions are satisfied:
\begin{align*}
&\a + g \c  -e_p \in  \sT\sP_{-,0}^{ p \times p }  \ands
g^\ast \a  + \c    \in \sT\sP_{+ ,0}^{ q \times p };  \\
&  g\d+\b  \in \sT\sP_{-,0}^{p\ts q}  \ands
\d+ g^*\b   -e_q\in \sT\sP_{+,0}^{q\ts q} .
\end{align*}
For this case, a solution to the twofold EG inverse problem has been obtained in
\cite[Section 9]{tHKvS1}.

\setcounter{equation}{0}
\section{Preliminaries about Toeplitz-like and Hankel-like operators}\label{sec:preliminaries}
In this section we define Toeplitz-like and Hankel-like operators and derive some of their properties.
First some notation. In what follows the direct sum  of  two linear spaces $\sN$ and $\sL$ will be denoted by   $\sN\dot{+} \sL$. Thus ( see \cite[pages 37, 38]{DS64})  the space $\sN\dot{+} \sL$ consists of all
$ (n,\ell) $ with $n \in \sN $ and $\ell \in \sL $ and its the linear structure is given by
\[
 (n_1,\ell_1)+ (n_2,\ell_2)=(n_1+n_2,\ell_1+\ell_2), \quad \l (n,\ell)=(\l n,\l\ell)\quad (\l\in \BC).
\]
In a canonical way $\sN$ and $\sL$ can be identified with the linear spaces
\[
\{(n, \ell)\mid n \in \sN,  \ell=0\in \sL\} \ands \{(n, \ell)\mid n=0 \in \sN,  \ell\in \sL\},
\]
respectively. We will use these identifications   without further  explanation.

Throughout this section   $\sA$,  $\sB $, $\sC $ and $\sD$ are  as in Section~\ref{sec:setting}, and we assume that $\sM_{\sA, \sB, \sC, \sD}$ is admissible. Put $ \sX = \sA \dot{+} \sB $ and $ \sY = \sC \dot{+} \sD $. Thus
$\sX $ is the direct sum of $\sA$ and $\sB$, and  $\sY $ is the direct sum of $\sC$ and $\sD$.
Furthermore, let
\begin{align}
& \sX_+ = \sA_+ \dot{+} \sB_+ , \quad \sX_- =  \sA_-^0 \dot{+} \sB_- ,\label{defX} \\
& \sY_+ = \sC_+ \dot{+} \sD_+^0 , \quad \sY_- = \sC_- \dot{+} \sD_- . \label{defY}
\end{align}
With these direct sums we associate four projections, denoted by
\[
\bP_{\sX_+}, \quad \bP_{\sX_-}, \quad \bP_{\sY_+}, \quad \bP_{\sY_-}.
\]
By definition,  $\bP_{\sX_+}$ is  the projection of $\sX$ onto $\sX_+$ along  $\sX_-$, and $\bP_{\sX_-}$ is the  projection of $\sX$ onto $\sX_-$ along  $\sX_+$ . The two other projections $\bP_{\sY_+}$ and $\bP_{\sY_-}$ are defined in a similar way, replacing $\sX$ by $\sY$.

We proceed with  defining multiplication (or Laurent-like) operators  and related Toeplitz-like and Hankel-like operators distinguishing  four cases. In each case the Toeplitz- and Hankel-like operators are compressions of the multiplication operators. Our terminology differs from the one used in \cite{RS85} and \cite{RS18}. Intertwining relations with shift-like operators appear later in the end of Section \ref{sec:solinverseprobl}, in Section \ref{sec:wienerL}, and in the Appendix.

\paragraph{1.  The case when $ \rho \in \sA $.}  Assume  $ \rho \in \sA $.
Then $ \rho \sA \subset \sA $ and $ \rho \sB \subset \sB $ and therefore  we have
for $ x = (\a, \b) \in \sX $ that
$ \rho x = (\rho\a, \rho\b ) \in \sX $, i.e.,  $ \rho \sX \subset \sX $.
We define the multiplication operator $ \bL_\rho : \sX \to \sX $ by
putting ${\bL}_\rho x = \rho x $ for $ x \in \sX$.
With respect to the decomposition $ \sX = \sX_- \dot{+} \sX_+ $ we write $ \bL_\rho $
as a $ 2 \times 2 $ operator matrix as follows
\begin{equation}\label{defAXX}
\bL_\rho = \begin{bmatrix} \bT_{-,\rho} & \bH_{-,\rho} \\  \bH_{+,\rho} & \bT_{+,\rho} \end{bmatrix} :
\begin{bmatrix} \sX_- \\ \sX_+ \end{bmatrix} \to \begin{bmatrix}  \sX_- \\ \sX_+ \end{bmatrix}.
\end{equation}
Thus for $ x_- \in \sX_- $ and $ x_+ \in \sX_+ $ we have
\begin{align*}
&  \bT_{-,\rho} x_- = \bP_{\sX_-} (\rho x_-) , \quad  \bT_{+,\rho} x_+ = \bP_{\sX_+} (\rho x_+),\\
&  \bH_{+,\rho} x_- = \bP_{\sX_+} (\rho x_-) , \quad  \bH_{-,\rho} x_+ = \bP_{\sX_-} (\rho x_+).
\end{align*}
We have $ \bL_\rho [\sA] \subset \sA $ and $ \bL_\rho [\sB] \subset \sB $.
Similarly, one has the inclusions
\begin{align}
&\bT_{\pm,\rho}[ \sA_\pm ] \subset \sA_\pm , \quad
\bT_{\pm,\rho}[ \sB_\pm ] \subset \sB_\pm,  \label{invAB1} \\
&\bH_{\pm,\rho}[ \sA_\mp ] \subset \sA_\pm \quad
\bH_{\pm,\rho}[ \sB_\mp ] \subset \sB_\pm .  \label{invAB2}
\end{align}
Furthermore, as expected from the classical theory of Hankel operators,  we have
\begin{equation}\label{Hmp_of_Apm}
\rho \in \sA_+ \Rightarrow \bH_{-,\rho} = 0\ands
\phi \in \sA_- \Rightarrow \bH_{+,\phi} = 0.
\end{equation}

\paragraph{2.  The case when $ \rho \in \sB $.}
For $ \rho \in \sB $ we have $ \rho \sC \subset \sA $ and $ \rho \sD \subset \sB $, and therefore $ \rho \sY \subset \sX $. We define the multiplication operator $ \bL_\rho : \sY \to \sX $ by putting $ \bL_\rho y = \rho y $ for $ y \in \sY $. With respect to the decompositions $ \sY = \sY_- \dot{+} \sY_+ $ and $ \sX = \sX_- \dot{+} \sX_+ $  we write $ \bL_\rho $ as a $ 2 \times 2 $ operator matrix as follows
\begin{equation}\label{defBYX}
{\bL_\rho} = \begin{bmatrix} \bT_{-,\rho} & \bH_{-,\rho} \\  \bH_{+,\rho} & \bT_{+,\rho} \end{bmatrix} :
\begin{bmatrix} \sY_- \\ \sY_+ \end{bmatrix} \to \begin{bmatrix}  \sX_- \\ \sX_+ \end{bmatrix}.
\end{equation}
Thus for $ y_- \in \sY_- $ and $ y_+ \in \sY_+ $ we have
\begin{align*}
&  \bT_{-,\rho} y_- = \bP_{\sX_-} (\rho y_-) , \quad  \bT_{+,\rho} y_+ = \bP_{\sX_+} (\rho y_+),\\
&  \bH_{+,\rho} y_- = \bP_{\sX_+} (\rho y_- ), \quad  \bH_{-,\rho} y_+ = \bP_{\sX_-} (\rho y_+ ).
\end{align*}
We have $ \bL_\rho [\sC] \subset \sA $ and $ \bL_\rho [\sD] \subset \sB $.
Similarly, one has
\begin{align}
&\bT_{\pm,\rho}[ \sC_\pm ] \subset \sA_\pm, \quad    \bT_{\pm,\rho}[ \sD_\pm ] \subset \sB_\pm, \label{invAB3}\\
&\bH_{\pm,\rho}[ \sC_\mp ] \subset \sA_\pm , \quad
\bH_{\pm,\rho}[ \sD_\mp ] \subset \sB_\pm . \label{invAB4}
\end{align}
Furthermore, we have
\begin{equation}\label{Hmp_of_Bpm}
\rho \in \sB_+ \Rightarrow \bH_{-,\rho} = 0 \ands
\phi \in \sB_- \Rightarrow \bH_{+,\phi} = 0.
\end{equation}

\paragraph{3.  The case when $ \rho \in \sC $.}
Let $ \rho \in \sC $. Then $ \rho \sA \subset \sC $ and $ \rho \sB \subset \sD $, and
therefore $ \rho \sX \subset \sY $.
We define the multiplication operator $ \bL_\rho : \sX \to \sY $ by
putting $ \bL_\rho x = \rho x $ for $ x \in \sX $.
With respect to the decomposition $ \sX = \sX_- \dot{+} \sX_+ $ and
$ \sY = \sY_- \dot{+} \sY_+ $ we write $ \bL_\rho $ as a $ 2 \times 2 $ operator matrix as follows
\begin{equation}\label{defCXY}
{\bL_\rho} = \begin{bmatrix} \bT_{-,\rho} & \bH_{-,\rho} \\  \bH_{+,\rho} & \bT_{+,\rho} \end{bmatrix} :
\begin{bmatrix} \sX_- \\ \sX_+ \end{bmatrix} \to \begin{bmatrix}  \sY_- \\ \sY_+ \end{bmatrix}.
\end{equation}
Thus for $ x_- \in \sX_- $ and $ x_+ \in \sX_+ $ we have
\begin{align*}
&  \bT_{-,\rho} x_- = \bP_{\sY_-} (\rho x_-) , \quad  \bT_{+,\rho} x_+ = \bP_{\sY_+}(\rho x_+), \\
&  \bH_{+,\rho} x_- = \bP_{\sY_+} (\rho x_-) , \quad  \bH_{-,\rho} x_+ = \bP_{\sY_-} (\rho x_+).
\end{align*}
We have $ \bL_\rho [\sA] \subset \sC $ and $ \bL_\rho [\sB] \subset \sD $.
Similarly, one has the inclusions
\begin{align}
&\bT_{\pm,\rho}[ \sA_\pm ] \subset \sC_\pm , \quad
\bT_{\pm,\rho}[ \sB_\pm ] \subset \sD_\pm ,\label{invCD1}\\
&\bH_{\pm,\rho}[ \sA_\mp ] \subset \sC_\pm,  \quad
\bH_{\pm,\rho}[ \sB_\mp ] \subset \sD_\pm .\label{invCD2}
\end{align}
Furthermore, we have
\begin{equation}\label{Hmp_of_Cpm}
\rho \in \sC_+ \Rightarrow  \bH_{-,\rho} = 0 \ands
\phi \in \sC_- \Rightarrow  \bH_{+,\phi} = 0.
\end{equation}

\paragraph{4.  The case when $ \rho \in \sD $.}
For $ \rho \in \sD $ we have $ \rho \sC \subset \sC $ and $ \rho \sD \subset \sD $, and
therefore $ \rho \sY \subset \sY $.
We define the multiplication operator $ \bL_\rho : \sY \to \sY $ by
putting $ \bL_\rho y = \rho y $ for $ y \in \sY $.
With respect to the decomposition $ \sY = \sY_- \dot{+} \sY_+ $  we write $ \bL_\rho $ as a $ 2 \times 2 $ operator matrix as follows
\begin{equation}\label{defDYY}
{\bL_\rho} = \begin{bmatrix} \bT_{-,\rho} & \bH_{-,\rho} \\  \bH_{+,\rho} & \bT_{+,\rho} \end{bmatrix} :
\begin{bmatrix} \sY_- \\ \sY_+ \end{bmatrix} \to \begin{bmatrix}  \sY_- \\ \sY_+ \end{bmatrix}.
\end{equation}
Thus for $ y_- \in \sY_- $ and $ y_+ \in \sY_+ $ we have
\begin{align*}
&  \bT_{-,\rho} y_- = \bP_{\sY_-} (\rho y_-) , \quad  \bT_{+,\rho} y_+ = \bP_{\sY_+} (\rho y_+),\\
&  \bH_{+,\rho} y_- = \bP_{\sY_+} (\rho y_- ), \quad  \bH_{-,\rho} y_+ = \bP_{\sY_-} (\rho y_+ ).
\end{align*}
We have $ \bL_\rho [\sC] \subset \sC $ and $ \bL_\rho [\sD] \subset \sD $.
Similarly, we have the inclusions
\begin{align}
&\bT_{\pm,\rho} [ \sC_\pm ] \subset \sC_\pm, \quad
\bT_{\pm,\rho} [ \sD_\pm ] \subset \sD_\pm, \label{invCD3}\\
&\bH_{\pm,\rho}[ \sC_\mp ] \subset \sC_\pm,   \quad
\bH_{\pm,\rho}[ \sD_\mp ] \subset \sD_\pm. \label{invCD4}
\end{align}
Furthermore, we have
\begin{equation}\label{Hmp_of_Dpm}
\rho \in \sD_+ \Rightarrow \bH_{-,\rho} = 0,  \qquad
\phi \in \sD_- \Rightarrow \bH_{+,\phi} = 0.
\end{equation}

\paragraph{{5. Multiplicative identities.}}
Let $ \sU $, $ \sV $ and $\sZ$ each be one of the spaces $ \sX $ or $\sY $ defined above. The corresponding decomposition of the spaces we denote as $ \sU = \sU_- \dot{+} \sU_+ $ and similarly for $ \sV $ and $\sZ $.
Let $ \phi $ be such that for $ u \in \sU $ we have $ \phi u \in \sV $ and
$ \rho $ be such that for $ v \in \sV $ we have $ \rho v \in \sZ $.
Then we have that $ \bL_{\rho \phi} = \bL_\rho \bL_\phi $, which gives
\begin{equation}\label{Laurent}
\begin{bmatrix} \bT_{-, \rho \phi } & \bH_{-, \rho \phi } \\
\bH_{+, \rho \phi } & \bT_{+, \rho \phi }   \end{bmatrix} =
\begin{bmatrix} \bT_{-, \rho } & \bH_{-, \rho} \\ \bH_{+, \rho } & \bT_{+, \rho  } \end{bmatrix}
\begin{bmatrix} \bT_{-, \phi } & \bH_{-, \phi } \\ \bH_{+, \phi } & \bT_{+,\phi }  \end{bmatrix}
: \begin{bmatrix} \sU_- \\ \sU_+ \end{bmatrix} \to
\begin{bmatrix} \sZ_- \\ \sZ_+ \end{bmatrix}.
\end{equation}
In particular, we have the following identities:
\begin{align}
& \bT_{+, \rho \phi } = \bT_{+, \rho } \bT_{+, \phi } + \bH_{+, \rho} \bH_{-, \phi }
: \sU_+ \to \sZ_+ ,\label{Tplusprod} \\
& \bH_{+, \rho \phi } = \bT_{+, \rho } \bH_{+, \phi } + \bH_{+, \rho} \bT_{-, \phi }
: \sU_- \to \sZ_+ , \label{Hplusprod} \\
& \bH_{-, \rho \phi } = \bH_{-, \rho } \bT_{+, \phi } + \bT_{-, \rho} \bH_{-, \phi }
: \sU_+ \to \sZ_- ,\label{Hminprod} \\
& \bT_{-, \rho \phi } = \bH_{-, \rho } \bH_{+, \phi } + \bT_{-, \rho} \bT_{-, \phi }
: \sU_- \to \sZ_- .\label{Tminprod}
\end{align}

\setcounter{equation}{0}
\section{Further notations and auxiliary results}\label{sec:prelim}

In this section we bring together a number of identities and lemmas that will be used in the proofs of the main results. Throughout this section $\a \in \sA_+ $, $\b \in \sB_+ $, $\c \in \sC_-$, and $\d \in \sD_- $. Furthermore, $g $ is an arbitrary element in $\sB_+$.  We split this section into two parts.

\smallskip
\noindent\textsc{Part 1.} With $g$  we associate the operator $\boldsymbol\Omega $ given by
\begin{equation}\label{defZg}
\boldsymbol\Omega = \begin{bmatrix} \bI_{\sX_+} & \bH_{+,g} \\ \bH_{-,g^*} & \bI_{\sY_-} \end{bmatrix}:
\begin{bmatrix} \sX_+ \\ \sY_- \end{bmatrix} \to \begin{bmatrix} \sX_+ \\ \sY_- \end{bmatrix}.
\end{equation}
Here $\sX_{\pm}$ and $\sY_{\pm}$ are as in \eqref{defX} and \eqref{defY}, respectively.
Using the  properties of  Hankel-like operators  given in the previous section we see that
\begin{align}
\eqref{inclu12} \ &\Longleftrightarrow\
\boldsymbol\Omega
\begin{bmatrix} \a \\ \c \end{bmatrix} = \begin{bmatrix} e_\sA \\ 0 \end{bmatrix},
\label{inclu12a} \\
\eqref{inclu34} \ &\Longleftrightarrow\
\boldsymbol\Omega
\begin{bmatrix} \b \\ \d \end{bmatrix} = \begin{bmatrix} 0 \\ e_\sD \end{bmatrix}.
\label{inclu34a}
\end{align}
Summarizing this yields  the following corollary.

\begin{cor}\label{cor:g-solinv} The element   $g \in \sB_+ $ is a solution to the  twofold  EG inverse problem associated with  $ \a $, $\b$, $\c$, and $\d$   if and only if
\begin{equation} \label{Omega-basic-eqs1}
\boldsymbol\Omega\begin{bmatrix} \a \\ \c \end{bmatrix} =
\begin{bmatrix} e_\sA \\ 0 \end{bmatrix} \ands
\boldsymbol\Omega\begin{bmatrix} \b \\ \d \end{bmatrix} = \begin{bmatrix} 0 \\ e_\sD \end{bmatrix},
\end{equation}
\end{cor}

We also have the following implications:
\begin{align}
\a + g \c \in \sA_- \ & \Longrightarrow \ \bH_{+,g \c} = -\bH_{+,\a}, \label{inclu1b}  \\
g^* \a + \c \in \sC_+ \ & \Longleftrightarrow \ \bH_{-, g^* \a} = -\bH_{-,\c } \label{inclu2b}  \\
\b + g \d \in \sB_- \ & \Longleftrightarrow \ \bH_{+,g\d} = -\bH_{+,\b},\label{inclu3b} \\
g^\ast \b + \d \in \sD_+ \ & \Longrightarrow \ \bH_{-,g^* \b} = - \bH_{-,\d}.
\label{inclu4b}
\end{align}
After taking adjoints in the left hand inclusions above we obtain
\begin{align}
\a + g \c \in \sA_- \ & \Longrightarrow \ \bH_{-,\c^* g^*} = -\bH_{-,\a^*}, \label{inclu1c}\\
g^\ast \a + \c \in \sC_+ \ & \Longleftrightarrow \ \bH_{+,\a^* g} = -\bH_{+,\c^*} \label{inclu2c}  \\
\b + g \d \in \sB_- \ & \Longleftrightarrow \ \bH_{-,\d^* g^*} = -\bH_{-,\b^*},\label{inclu3c} \\
g^\ast \b + \d \in \sD_+ \ & \Longrightarrow \ \bH_{+,\b^* g} = - \bH_{+,\d^*}.
\label{inclu4c}
\end{align}
Notice that the first inclusion in \eqref{inclu12} implies that $ \a + g \c \in \sA_- $ and
the second inclusion in \eqref{inclu34} implies $ g^\ast \b + \d \in \sD_+$. The implications from left to right are obvious.
To prove the implications from right to left in \eqref{inclu2b}, \eqref{inclu3b},
\eqref{inclu2c} and \eqref{inclu3c} one reasons as follows.
For example, for \eqref{inclu2b} one uses that $ e_\sA \in \sX_+ $, such that
\[
0 =\bH_{- ,g^\ast \a + \c } e_\sA = \bP_{ \sX_- } (g^\ast \a + \c)e_\sA =
\bP_{ \sC_- } (g^\ast \a + \c).
\]
Hence $g^\ast \a + \c\in \sC_+$, as claimed. Since $ e_\sA \not\in \sX_- $ and $ e_\sD \not\in \sY_+ $, the reverse implications in
\eqref{inclu1b}, \eqref{inclu4b}, \eqref{inclu1c} and \eqref{inclu4c}
cannot be derived in this way.

Note that $\a \in \sA_+ $, $\b \in \sB_+ $, $\c \in \sC_-$, $\d \in \sD_- $ implies  that
\begin{align}
&\bH_{-,\a} = 0, \quad \bH_{-,\b} = 0, \quad \bH_{+,\c} = 0, \quad \bH_{+,\d} = 0, \label{Hzero1}\\
&\bH_{+,\a^*} = 0, \quad \bH_{+,\b^*} = 0, \quad \bH_{-,\c^*} = 0, \quad \bH_{-,\d^*} = 0. \label{Hzero2}
\end{align}
Using the identities \eqref{Hzero1} and \eqref{Hzero2} together  with the product formulas at the end of Section~\ref{sec:preliminaries} we obtain the following eight identities:
\begin{align}
\bH_{+,\a^*g} & = \bT_{+,\a^*} \bH_{+,g}, \quad \bH_{+,\b^*g} = \bT_{+,\b^*} \bH_{+,g}, \label{Hplus*g} \\
\bH_{-,\c^*g^*} & = \bT_{-,\c^*} \bH_{-,g^*}, \quad \bH_{-,\d^*g^*} = \bT_{-,\d^*} \bH_{-,g^*},
\label{Hmin*g*} \\
\bH_{+,g\c} & = \bH_{+,g} \bT_{-,\c} , \quad \bH_{+,g\d} = \bH_{+,g} \bT_{-,\d} \label{Hplusg*}  \\
\bH_{-,g^\ast \a} & = \bH_{-,g^\ast} \bT_{+,\a} , \quad \bH_{-,g^\ast \b} = \bH_{-,g^\ast} \bT_{+,\b}
\label{Hming**}
\end{align}

The next  lemma is an immediate consequence of the definitions.

\begin{lem}\label{lem:Hunit} For $g\in \sB_+$  and $h\in \sC_-$ we have
\[
\bH_{+,g}e_\sD=g \ands  \bH_{-,h}e_\sA=h.
\]
\end{lem}
We conclude this part with the following lemma.

\begin{lem}\label{basicids5}
Assume that conditions \textup{(C1)--(C3)} are satisfied. Then
\begin{equation}\label{matrixTHHT}
\begin{bmatrix}  \bT_{+,\a^\ast} \\ \bT_{+,\b^\ast}\end{bmatrix}
\begin{bmatrix} \bH_{+,\b} & \bH_{+,\a} \end{bmatrix} =
\begin{bmatrix}  \bH_{+,\c^\ast} \\ \bH_{+,\d^\ast}\end{bmatrix}
\begin{bmatrix} \bT_{-,\d} & \bT_{-,\c} \end{bmatrix}.
\end{equation}
and
\begin{equation}\label{matrixTHHT2}
\begin{bmatrix}  \bT_{-,\d^\ast} \\ \bT_{-,\c^\ast}\end{bmatrix}
\begin{bmatrix} \bH_{-,\c} & \bH_{-,\d} \end{bmatrix} =
\begin{bmatrix}  \bH_{-,\b^\ast} \\ \bH_{-,\a^\ast}\end{bmatrix}
\begin{bmatrix} \bT_{+,\a} & \bT_{+,\b} \end{bmatrix}.
\end{equation}
\end{lem}

\bpr
In the course of the proof  we repeatedly use the product rules \eqref{Tplusprod}--\eqref{Tminprod}.
Using the first identity in \eqref{Hzero2}, condition (C3) and the fourth identity in
\eqref{Hzero1} we see that
\begin{align*}
 \bT_{+, \a^\ast} \bH_{+,\b}& = \bH_{+,\a^\ast \b } - \bH_{+,\a^\ast } \bT_{-,\b}\\
 & = \bH_{+,\a^\ast \b }= \bH_{+,\c^\ast \d }\\
& = \bH_{+,\c^\ast } \bT_{-,\d} + \bT_{+,\c^\ast } \bH_{+,\d} =
\bH_{+,\c^\ast } \bT_{-,\d}.
\end{align*}
It follows that
\begin{equation}\label{THHT1}
\bT_{+, \a^\ast} \bH_{+,\b} = \bH_{+,\c^\ast} \bT_{-,\d}.
\end{equation}
Next, using the first identity in \eqref{Hzero2} and the third in \eqref{Hzero1} we obtain
\begin{align*}
\bT_{+, \a^\ast} \bH_{+,\a}& = \bH_{+,\a^\ast \a} - \bH_{+, \a^\ast}\bT_{-,\a} = \bH_{+,\a^\ast \a},\\
\bH_{+,\c^\ast} \bT_{-,\c}& = \bH_{+,\c^\ast \c} - \bT_{+,\c^\ast} \bH_{+,\c} = \bH_{+,\c^\ast \c}.
\end{align*}
On the other hand, using  condition (C1)  and  the second identity in \eqref{Hmp_of_Apm}
with $\phi= a_0 = \bP_{\sA_d} \a $, we see that
$ \bH_{+,\a^\ast \a} - \bH_{+,\c^\ast \c}  = \bH_{+,a_0} = 0  $.  We proved
\begin{equation}\label{THHT2}
\bT_{+, \a^\ast} \bH_{+,\a} = \bH_{+,\c^\ast} \bT_{-,\c}.
\end{equation}
The next two equalities are proved in a similar way as the previous two:
\begin{equation}\label{THHT3}
\bT_{+, \b^\ast} \bH_{+,\b} = \bH_{+,\d^\ast} \bT_{-,\d} ,
\end{equation}
\begin{equation}\label{THHT4}
\bT_{+, \b^\ast} \bH_{+,\a} = \bH_{+,\d^\ast} \bT_{-,\c} .
\end{equation}
Observe that \eqref{THHT1}, \eqref{THHT2}, \eqref{THHT3}  and
\eqref{THHT4} can be rewritten as  \eqref{matrixTHHT}.

The equality \eqref{matrixTHHT2} is proved  similarly,
using (C2)  instead of (C1).
\epr

\bigskip
\noindent\textsc{Part 2.}  In the second part of this section we assume that $ a_0 = \bP_{\sA_d} \a $ and $ d_0 = \bP_{\sA_d} \d $ are invertible in $ \sA_d $ and $ \sD_d $, respectively.
Using the notations introduced in the previous section we associate with the elements  $\a, \b, \c, \d$ the following  operators:
\begin{align}
& \bR_{11} = \bT_{+,\a} a_0^{-1} \bT_{+,\a^\ast} - \bT_{+,\b} d_0^{-1} \bT_{+,\b^\ast} : \sX_+ \to \sX_+ ,
\label{defR11} \\
& \bR_{21} = \bH_{-,\c} a_0^{-1} \bT_{+,\a^\ast} - \bH_{-,\d} d_0^{-1} \bT_{+,\b^\ast} : \sX_+ \to \sY_- , \label{defR21}\\
& \bR_{12} = \bH_{+,\b} d_0^{-1} \bT_{-,\d^\ast} - \bH_{+,\a }a_0^{-1} \bT_{-,\c^\ast} : \sY_- \to \sX_+ , \label{defR12}  \\
& \bR_{22} = \bT_{-,\d} d_0^{-1} \bT_{-,\d^\ast} - \bT_{-,\c} a_0^{-1} \bT_{-,\c^\ast} : \sY_- \to \sY_- .
\label{defR22}
\end{align}

\begin{lem}\label{lem:unitsR}
Assume that conditions \textup{(C1)} and \textup{(C2)} are satisfied and that $a_0 $ and $ d_0 $ are invertible in $ \sA_d $ and $ \sD_d $, respectively.
Then the following identities hold true:
\begin{equation} \label{unitsR}
\bR_{11}e_\sA=\a, \quad  \bR_{12}e_\sD=\b, \quad \bR_{21}e_\sA=\c, \quad  \bR_{22}e_\sD=\d.
\end{equation}
\end{lem}

\bpr
Note that $\b^*\in \sC_-$.
Thus
\[\bT_{+, \b^*}e_\sA=\bP_{\sY_+}(\b^*e_\sA)=\bP_{\sY_+}\b^*=0.
\]
{Since $\a^*\in\sA_-^*$, we have}
\[
{\bT_{+,\a^\ast}e_\sA
=\bP_{\sX_+}(\a^*e_\sA)=\bP_{\sX_+} \a^* =a_0^*.}
\]
Using $a_0=a_0^*$ (by the first part of \eqref{selfadj-ad0}), it follows that
\begin{align*}
\bR_{11}e_\sA &= \bT_{+,\a} a_0^{-1} \bT_{+,\a^\ast}e_\sA
{= \bT_{+,\a} a_0^{-1}a_0^*=\bT_{+,\a} e_\sA=P_{\sX_+}(\a e_\sA)=\a.}
\end{align*}
Notice that we used condition (C1).
This proves the first identity \eqref{unitsR}.

Next, using $\c\in \sC_-$, $\bT_{+,\a^\ast}e_\sA=a_0^*$, and
$\bT_{+,\b^*}e_{\sA}=\bP_{\sY_+} \b^\ast = 0 $ we obtain
\[
\bR_{21}e_\sA = \bH_{-,\c} a_0^{-1} \bT_{+,\a^\ast}e_\sA=\bH_{-, \c} e_\sA=\bP_{\sY_-}(\c e_\sA)=\bP_{\sY_-}\c=\c,
\]
which proves  the third  identity in \eqref{unitsR}.
The two other identities in \eqref{unitsR}, involving $\bR_{12}$  and $\bR_{22}$,
are obtained in a similar way, using (C2), \eqref{selfadj-ad0} and
\[
\bT_{-, \c^*} e_\sD = 0, \quad  \bT_{-, \d^*} e_\sD = d_0^*, \quad
\bH_{+,\b} e_\sD = \b , \quad \bT_{-, \d} e_\sD = \d.
\]
This proves the lemma.
\epr

\medskip
The next lemma presents alternative formulas for the operators $\bR_{ij}$, $1\leq i,j \leq 2$, given by \eqref{defR11}--\eqref{defR22}, assuming conditions  (C4)--(C6) are satisfied.

\begin{lem}\label{lemaltRij}
Assume that $a_0 $ and $ d_0 $ are invertible in $ \sA_d $ and $ \sD_d $, respectively,
and that conditions \textup{(C4), (C5)}, and \textup{(C6)} are satisfied.
Then
\begin{align}
& \bR_{11} = \bI_{\sX_+} - \bH_{+,\a} a_0^{-1} \bH_{-,\a^\ast} +
\bH_{+,\b} d_0^{-1} \bH_{-,\b^\ast}  : \sX_+ \to \sX_+,
\label{formulaR11} \\
& \bR_{21} =  \bT_{-,\d} d_0^{-1} \bH_{-,\b^\ast } - \bT_{-,\c} a_0^{-1} \bH_{-,\a^\ast}
: \sX_+ \to \sY_- , \label{formulaR21}\\
& \bR_{12} =  \bT_{+,\a} a_0^{-1} \bH_{+,\c^\ast} - \bT_{+,\b} d_0^{-1}  \bH_{+,\d^\ast}
 : \sY_- \to \sX_+ , \label{formulaR12} \\
& \bR_{22} = \bI_{\sY_-} - \bH_{-,\d} d_0^{-1} \bH_{+,\d^\ast} +
\bH_{-,\c} a_0^{-1} \bH_{+,\c^\ast} : \sY_- \to \sY_-.
\label{formulaR22}
\end{align}
\end{lem}

\bpr
First notice that    $ a_0^{-1} \in \sA_d $ and $ d_0^{-1} \in \sD_d $ yield the
following identities
\begin{align}
& \bT_{+,\a a_0^{-1}} = \bT_{+,\a} a_0^{-1},\quad \bH_{+,\a a_0^{-1}} = \bH_{+,\a} a_0^{-1}, \label{pla0}\\
& \bT_{+,\b d_0^{-1}} = \bT_{+,\b} d_0^{-1},\quad \bH_{+,\b d_0^{-1}} = \bH_{+,\b} d_0^{-1}, \label{plb0}\\
& \bT_{-,\d d_0^{-1}} = \bT_{-,\d} d_0^{-1},\quad \bH_{-,\d d_0^{-1}} = \bH_{-,\d} d_0^{-1}, \label{mid0}\\
& \bT_{-,\c a_0^{-1}} = \bT_{-,\c} a_0^{-1},\quad \bH_{-,\c a_0^{-1}} = \bH_{-,\c} a_0^{-1}. \label{mic0}
\end{align}
Next, note that condition (C4) implies that
$ \bT_{+,\a a_0^{-1} \a^\ast - \b d_0^{-1} \b^\ast - e_\sA } = 0 $.
It follows that
\[
 \bT_{+,\a a_0^{-1} \a^\ast } - \bT_{+,\b d_0^{-1} \b^\ast } - \bI_{\sX_+} =0 .
\]
Applying the product rule  \eqref{Tplusprod} and the identities in \eqref{pla0} and \eqref{plb0} we see that
\[
\bT_{+,\a } a_0^{-1}\bT_{+, \a^\ast }-\bT_{+,\b}  d_0^{-1} \bT_{+, \b^\ast }= \bI_{\sX_+} -  \bH_{+,\a}  a_0^{-1}\bH_{-, \a^\ast }+\bH_{+,\b } d_0^{-1}\bH_{-, \b^\ast}.
\]
It follows that the operator $\bR_{11}$ defined by \eqref{defR11} is also given by
\eqref{formulaR11}.
In a similar way one shows that condition (C5) yields the identity \eqref{formulaR22}.

Since (C6) states $ \a a_0^{-1} \c^\ast = \b d_0^{-1} \d^\ast $, we have
the equality
$ \bH_{+,\a a_0^{-1} \c^\ast  } = \bH_{+, \b d_0^{-1} \d^\ast }$.
Applying the product rule \eqref{Hplusprod} and the identities in \eqref{pla0}
and \eqref{plb0} it follows that
\[
\bH_{+,\a } a_0^{-1}\bT_{-, \c^\ast } + \bT_{+,\a } a_0^{-1} \bH_{+, \c^\ast }=
\bH_{+, \b } d_0^{-1}\bT_{-,\d^\ast } + \bT_{+, \b}  d_0^{-1}\bH_{+,\d^\ast }.
\]
This yields
\[
\bR_{12} = \bH_{+,\b } d_0^{-1} \bT_{-,\d^\ast} - \bH_{+,\a }a_0^{-1} \bT_{-,\c^\ast} =
\bT_{+,\a} a_0^{-1} \bH_{+,\c^\ast} - \bT_{+,\b} d_0^{-1}  \bH_{+,\d^\ast},
\]
which proves \eqref{formulaR12}.

Finally, to prove the identity \eqref{formulaR21}, note that, by taking adjoints,
condition  (C6)  yields that $ \d d_0^{-1} \b^\ast = \c a_0^{-1} \a^\ast $.
But then using the identities in \eqref{mid0} and \eqref{mic0},
arguments similar to the ones used in the previous paragraph,
yield  the identity \eqref{formulaR21}.
\epr

\bigskip
The following lemma contains some useful formulas that we will prove by direct verification.

\begin{lem}\label{lem:inverseR}
Assume that $a_0 $ and $ d_0 $ are invertible in $ \sA_d $ and $ \sD_d $, respectively,
and that the conditions \textup{(C1)--(C6)} are satisfied.
Let $ \bR_{ij} $, $ i,j=1,2$,  be given by \eqref{defR11}--\eqref{defR22}. Then
\begin{equation}\label{equality1}
\begin{bmatrix} \bR_{11} & \bR_{12} \\ \bR_{21} &  \bR_{22}  \end{bmatrix}
\begin{bmatrix} \bI_{\sX_+}  & 0 \\ 0 &  -\bI_{\sY_-}  \end{bmatrix}
\begin{bmatrix} \bR_{11} & \bR_{12} \\ \bR_{21} &  \bR_{22}  \end{bmatrix}  =
\begin{bmatrix} \bR_{11} & 0 \\ 0 &  -\bR_{22}  \end{bmatrix}.
\end{equation}
This implies that
\begin{equation}\label{invform1}
\bR = \begin{bmatrix} \bR_{11} & \bR_{12} \\ \bR_{21} &  \bR_{22}  \end{bmatrix} :
\begin{bmatrix} \sX_+ \\ \sY_- \end{bmatrix}  \to \begin{bmatrix} \sX_+ \\ \sY_- \end{bmatrix}
\end{equation}
is invertible if and only if $ \bR_{11} $ and $ \bR_{22} $ are invertible.
Furthermore, in that case
\begin{equation}\label{invform2}
\bR^{-1}  =
\begin{bmatrix} \bI_{\sX_+} & -\bR_{12} \bR_{22}^{-1} \\ -\bR_{21} \bR_{11}^{-1} & \bI_{\sY_-} \end{bmatrix}
=
\begin{bmatrix} \bI_{\sX_+} & -\bR_{11}^{-1} \bR_{12}  \\ -\bR_{22}^{-1} \bR_{21}  & \bI_{\sY_-} \end{bmatrix}.
\end{equation}
\end{lem}

\bpr
To check \eqref{equality1} we will prove the four identities
\begin{align}
& \bR_{11} \bR_{12} = \bR_{12} \bR_{22}, \quad \bR_{22} \bR_{21} = \bR_{21} \bR_{11}, \label{eq1} \\
& \bR_{11} \bR_{11} - \bR_{12} \bR_{21} = \bR_{11} , \quad \bR_{22} \bR_{22} - \bR_{21} \bR_{12} = \bR_{22} .
\label{eq2}
\end{align}

From \eqref{defR11}--\eqref{defR22} and \eqref{formulaR11}--\eqref{formulaR22} it follows  that
\begin{align*}
& \bR_{11} \bR_{12} = \begin{bmatrix} \bT_{+,\a} & \bT_{+,\b} \end{bmatrix}
\begin{bmatrix} a_0^{-1} & 0 \\ 0 & -d_0^{-1} \end{bmatrix}
\begin{bmatrix}  \bT_{+,\a^\ast} \\ \bT_{+,\b^\ast}\end{bmatrix} \times\\
&\hspace{4cm} \times
\begin{bmatrix} \bH_{+,\b} & \bH_{+,\a} \end{bmatrix}
\begin{bmatrix} d_0^{-1} & 0 \\ 0 & -a_0^{-1} \end{bmatrix}
\begin{bmatrix}  \bT_{-,\d^\ast} \\ \bT_{-,\c^\ast}\end{bmatrix},
\end{align*}
and
\begin{align*}
& \bR_{12} \bR_{22} = \begin{bmatrix} \bT_{+,\a} & \bT_{+,\b} \end{bmatrix}
\begin{bmatrix} a_0^{-1} & 0 \\ 0 & -d_0^{-1} \end{bmatrix}
\begin{bmatrix}  \bH_{+,\c^\ast} \\ \bH_{+,\d\ast}\end{bmatrix}\times \\
&\hspace{4cm} \times\begin{bmatrix} \bT_{-,\d} & \bT_{-,\c} \end{bmatrix}
\begin{bmatrix} d_0^{-1} & 0 \\ 0 & -a_0^{-1} \end{bmatrix}
\begin{bmatrix}  \bT_{-,\d^\ast} \\ \bT_{-,\c^\ast}\end{bmatrix}.
\end{align*}
But then \eqref{matrixTHHT} shows that $\bR_{11} \bR_{12}=\bR_{12} \bR_{22}$. In a
similar way, using \eqref{matrixTHHT2} one proves that $ \bR_{22} \bR_{21} = \bR_{21} \bR_{11} $.

Next observe that
\begin{align*}
& \bR_{11} ( \bR_{11} - \bI_{\sX_+} ) =  - \begin{bmatrix} \bT_{+,\a} & \bT_{+,\b} \end{bmatrix}
\begin{bmatrix} a_0^{-1} & 0 \\ 0 & -d_0^{-1} \end{bmatrix}
\begin{bmatrix}  \bT_{+,\a^\ast} \\ \bT_{+,\b^\ast}\end{bmatrix}\times \\
&\hspace{4cm}\times  \begin{bmatrix} \bH_{+,\a} & \bH_{+,\b} \end{bmatrix}
\begin{bmatrix} a_0^{-1} & 0 \\ 0 & -d_0^{-1} \end{bmatrix}
\begin{bmatrix}  \bH_{-,\a^\ast} \\ \bH_{-,\b^\ast}\end{bmatrix},
\end{align*}
and
\begin{align*}
 \bR_{12} \bR_{21} &=  \begin{bmatrix} \bT_{+,\a} & \bT_{+,\b} \end{bmatrix}
\begin{bmatrix} a_0^{-1} & 0 \\ 0 & -d_0^{-1} \end{bmatrix}
\begin{bmatrix}  \bH_{+,\c^\ast} \\ \bH_{+,\d\ast}\end{bmatrix}\times\\
&\hspace{4cm}\times\begin{bmatrix} \bT_{-,\d} & \bT_{-,\c} \end{bmatrix}
\begin{bmatrix} d_0^{-1} & 0 \\ 0 & -a_0^{-1} \end{bmatrix}
\begin{bmatrix}  \bH_{-,\b^\ast} \\ \bH_{-,\a^\ast}\end{bmatrix}  \\
&= - \begin{bmatrix} \bT_{+,\a} & \bT_{+,\b} \end{bmatrix}
\begin{bmatrix} a_0^{-1} & 0 \\ 0 & -d_0^{-1} \end{bmatrix}
\begin{bmatrix}  \bH_{+,\c^\ast} \\ \bH_{+,\d\ast}\end{bmatrix}\times\\
&\hspace{4cm}\times\begin{bmatrix} \bT_{-,\d} & \bT_{-,\c} \end{bmatrix}
\begin{bmatrix} a_0^{-1} & 0 \\ 0 & -d_0^{-1} \end{bmatrix}
\begin{bmatrix}  \bH_{-,a^\ast} \\ \bH_{-,\b^\ast}\end{bmatrix}
\end{align*}
But then \eqref{matrixTHHT} implies that  $\bR_{11} \bR_{11} - \bR_{12} \bR_{21} = \bR_{11}$.
Similarly, using \eqref{matrixTHHT2} one proves that $ \bR_{22} \bR_{22} - \bR_{21} \bR_{12} = \bR_{22} $.

The final statements \eqref{invform1} and \eqref{invform2} are immediate from
\eqref{eq1} and \eqref{eq2}.
\epr

\setcounter{equation}{0}
\section{An abstract inversion theorem}\label{sec:inversionthm}
Let  $\sM=\sM_{\sA, \sB, \sC, \sD}$ be  an  admissible  algebra. Fix $g\in \sB_+$, and let $\boldsymbol\Omega$ be the operator given by
\begin{equation}\label{defZg1}
\boldsymbol\Omega: = \begin{bmatrix} I_{\sX_+} & \bH_{+,g} \\[.1cm] \bH_{-,g^*} & \bI_{\sY_-} \end{bmatrix}:
\begin{bmatrix} \sX_+ \\[.1cm] \sY_- \end{bmatrix} \to \begin{bmatrix} \sX_+\\[.1cm] \sY_- \end{bmatrix}.
\end{equation}
We shall prove the following inversion theorem.

\begin{thm}\label{thm:inversion1}
Let  $\sM=\sM_{\sA, \sB, \sC, \sD}$ be  an  admissible  algebra and  let $g\in \sB_+$. Then the operator $\boldsymbol\Omega$ defined  by \eqref{defZg1} is invertible if there exist $\a \in \sA_+ $, $\b \in \sB_+ $, $\c \in \sC_-$, $\d \in \sD_- $ such that
\begin{equation} \label{Omega-basic-eqs}
\boldsymbol\Omega\begin{bmatrix} \a \\ \c \end{bmatrix} =
\begin{bmatrix} e_\sA \\ 0 \end{bmatrix} \ands
\boldsymbol\Omega\begin{bmatrix} \b \\ \d \end{bmatrix} = \begin{bmatrix} 0 \\ e_\sD \end{bmatrix},
\end{equation}
and  the following two conditions are satisfied:
\begin{itemize}
\item[\textup{(a)}]   $ a_0 := P_{\sA_d} \a $ and $ d_0 := P_{\sD_d} \d $ are invertible in $ \sA_d $ and $ \sD_d $, respectively;
\item[\textup{(b)}]  conditions \textup{(C4)--(C6)} are satisfied
\end{itemize}
In that case  the inverse of $\boldsymbol\Omega$ is given by
\begin{equation}\label{IGGIinverse}
\boldsymbol\Omega^{-1}
=\begin{bmatrix} \bR_{11} &  \bR_{12} \\ \bR_{21} &  \bR_{22} \end{bmatrix},
\end{equation}
where $\bR_{ij}$, $1\leq i,j \leq 2$, are the operators  defined by  \eqref{defR11}--\eqref{defR22}. Furthermore, the operators $\bR_{11}$ and $\bR_{22}$ are invertible and
\begin{align}
& \bH_{+,g} = - \bR_{11}^{-1} \bR_{12} = - \bR_{12} \bR_{22}^{-1} , \label{formula-bHg} \\
& \bH_{-,g^\ast } = - \bR_{21} \bR_{11}^{-1} = - \bR_{22}^{-1} \bR_{21},  \label{formula-bHgast}\\
& g  = - \bR_{11}^{-1} \b  ,\quad g^\ast = - \bR_{22}^{-1} \c. \label{formggR}
\end{align}
\end{thm}

\begin{rem}\label{rem:not-iff}
In contrast to Theorem \ref{thm:mainL1}, the above theorem is not an ``if and only if"  statement.  In this general  setting we only have the following partial converse:  if  the operator $\boldsymbol\Omega$ given by \eqref{defZg1} is invertible, then there exist $\a \in \sA_+ $, $\b \in \sB_+ $, $\c \in \sC_-$, $\d \in \sD_- $  such that the equations \eqref{Omega-basic-eqs} are satisfied.  It can happen that   the operator $\boldsymbol\Omega$ is invertible and item (a)  is not satisfied; see Example \ref{exa:not-inv} given at the end of the present section.
\end{rem}

Note that the operators $\bR_{ij}$, $1\leq i,j\leq 2$, appearing in \eqref{IGGIinverse} do  not depend on the particular choice of $g$, but on $\a$, $\b$, $\c$, $\d$ only. It follows that Theorem \ref{thm:inversion1} yields the following corollary.

\begin{cor}\label{cor:uniq}
Let $\a \in \sA_+ $, $\b \in \sB_+ $, $\c \in \sC_-$, $\d \in \sD_- $, and assume that
\begin{itemize}
\item[\textup{(a)}] $ a_0 = P_{\sA_d} \a $ and $ d_0 = P_{\sD_d} \d $ are invertible
in $ \sA_d $ and $ \sD_d $, respectively;
\item[\textup{(b)}] conditions \textup{ (C4)--(C6) } are satisfied.
\end{itemize}
Under these conditions,  if the twofold  EG inverse problem associated with
$\a \in \sA_+ $, $\b \in \sB_+ $, $\c \in \sC_-$,  $\d \in \sD_- $ has a  solution,
then  the solution is unique.
\end{cor}

\bpr
Assume that   the twofold  EG inverse problem associated with $\a \in \sA_+ $, $\b \in \sB_+ $, $\c \in \sC_-$,  $\d \in \sD_- $ has a solution, $g$ say. Then (see Corollary \ref{cor:g-solinv})   the two identities in \eqref{Omega-basic-eqs} are satisfied. Furthermore, by assumption, items  (a) and (b) in Theorem \ref{thm:inversion1} are satisfied too. We conclude that  \eqref{IGGIinverse} holds, and hence  $\bOm$  is uniquely determined by   the operators $\bR_{ij}$, $1\leq i,j\leq 2$. But these  $\bR_{ij}$, $1\leq i,j\leq 2$,  do  not depend on  $g$, but  on $\a$, $\b$, $\c$, $\d$ only. It follows   that the same is true for $ \bH_{+,g}$.  But   $ \bH_{+,g} e_\sD =P_{\sX_+} g e_\sD = P_{\sX_+} g = g $. Thus $g$ is uniquely determined by the data.
\epr

\medskip
The following lemma will be be useful in the proof of  Theorem \ref{thm:inversion1}.

\begin{lem}\label{lem:basicids6}
Let $g\in \sB_+$ satisfy the inclusions \eqref{inclu12} and \eqref{inclu34}. Then the following identities hold:
\begin{align}
&\bT_{+,\a^*} \bH_{+,g} = - \bH_{+, \c^*}, \quad  \bT_{+,\b^*} \bH_{+,g} = - \bH_{-, \d^*},\label{ident6a}\\
&\bT_{-,\c^*} \bH_{-,g^*} = - \bH_{-, \a^*}, \quad \bT_{-,\d^*} \bH_{-,g^*} = - \bH_{-, \b^*},\label{ident6b}\\
&\bH_{+, g} \bT_{-, \d} = - \bH_{+, \b}, \quad  \bH_{+, g} \bT_{-, \c}= - \bH_{+, \a},\label{ident6c}\\
&\bH_{-, g^*} \bT_{+, \a} = - \bH_{-, \c}, \quad \bH_{-, g^*} \bT_{+, \b}= - \bH_{-,\d}.\label{ident6d}
\end{align}
\end{lem}
\bpr The above   identities follow by using the implications
in \eqref{inclu1b}--\eqref{inclu4b} and \eqref{inclu1c}--\eqref{inclu4c}
together  with the identities in \eqref{Hplus*g} -- \eqref{Hming**}.
Let us illustrate this by proving the first identity in \eqref{ident6a}.

From the first identity in \eqref{Hplus*g}  we know  that $\bT_{+,\a^*}\bH_{+,g}=\bH_{+, \a^*g}$.
Since $g\in \sB_+$ satisfies the first inclusion in \eqref{inclu12}, the equivalence
in \eqref{inclu2c} tells us that $\bH_{+, \a^*g}=-\bH_{+, \c^*}$.
Hence $\bT_{+,\a^*}\bH_{+,g}=-\bH_{+, \c^*}$, and the first identity \eqref{ident6a} is proved.
\epr

\medskip\noindent
{\bf Proof of Theorem \ref{thm:inversion1}.}
Recall that the two  identities in \eqref{Omega-basic-eqs} together are equivalent to  $ g \in \sB_+ $ being a solution to the  twofold  EG inverse problem associated with $\a \in \sA_+ $, $\b \in \sB_+ $, $\c \in \sC_-$, $\d \in \sD_- $, and hence the two  identities  in  \eqref{Omega-basic-eqs}  imply that the conditions (C1)--(C3) are satisfied. Given item (b) in Theorem \ref{thm:inversion1} we conclude that all conditions (C1)--(C6) are satisfied.

The remainder of  the proof is divided into three parts.

\noindent\textsc{Part 1.}
First we will prove that
\begin{equation}\label{RLinvers}
\begin{bmatrix} \bR_{11} &  \bR_{12} \\ \bR_{21} &  \bR_{22} \end{bmatrix}
\begin{bmatrix} I_{\sX_+} &  \bH_{+,g} \\ \bH_{-,g^\ast} &  I_{\sY_-} \end{bmatrix} =
\begin{bmatrix} I_{\sX_+} & 0 \\ 0 &  I_{\sY_-} \end{bmatrix}.
\end{equation}
We start with the identity $ \bR_{11} + \bR_{12} \bH_{-,g^\ast} = I_{\sX_+}$. Using the two identities in \eqref{ident6b} we have
\begin{align*}
\bR_{12} \bH_{-, g^\ast } & = \bH_{+,\b } d_0^{-1} \bT_{-,\d^\ast } \bH_{-,g^\ast }
- \bH_{+,\a }a_0^{-1} \bT_{-,\c^\ast} \bH_{-,g^\ast } \\
& = - \bH_{+,\b } d_0^{-1} \bH_{-,\b^\ast} + \bH_{+,\a }a_0^{-1} \bH_{-,\a^\ast} \\
& = - \bR_{11} + I_{\sX_+},
\end{align*}
which proves $ \bR_{11} + \bR_{12} \bH_{-,g^\ast} = I_{\sX_+} $.

Similarly, using the two identities  in \eqref{ident6a}, we obtain
\begin{align*}
\bR_{11}  \bH_{+,g} &= \bT_{+,\a} a_0^{-1} \bT_{+,\a^\ast}  \bH_{+,g} -
\bT_{+,\b} d_0^{-1} \bT_{+,\b^\ast} \bH_{+,g}\\
&= - \bT_{+,\a} a_0^{-1} \bH_{ +, \c^\ast } + \bT_{+,\b} d_0^{-1} \bH_{ +, \d^\ast }\\
&=- \bR_{12}.
\end{align*}
Thus $\bR_{11}  \bH_{+,g}+\bR_{12}=0$.

The equalities $ \bR_{21} \bH_{+,g} + \bR_{22} = I_{\sY} $ and  $ \bR_{21} + \bR_{22} \bH_{-,g^\ast}=0 $ are proved in a similar way.

\medskip
\noindent\textsc{Part 2.}
In this part  we prove that
\begin{equation}\label{RRinvers}
\begin{bmatrix} I_{\sX_+} &  \bH_{+,g} \\ \bH_{-,g^\ast} &  I_{\sY_-} \end{bmatrix}
\begin{bmatrix} \bR_{11} &  \bR_{12} \\ \bR_{21} &  \bR_{22} \end{bmatrix}=
\begin{bmatrix} I_{\sX_+} & 0 \\ 0 &  I_{\sY_-} \end{bmatrix}.
\end{equation}
To see  this we  first show that $ \bR_{11} + \bH_{+,g} \bR_{21} = I_{\sX_+}$. We use \eqref{formulaR21} and the two identities  in \eqref{ident6c}. This yields
\begin{align*}
\bH_{+,g} \bR_{21}& = \bH_{+,g} \bT_{-,\d} d_0^{-1} \bH_{-,\b^\ast } - \bH_{+,g} \bT_{-,\c} a_0^{-1} \bH_{-,\a^\ast}\\
&=- \bH_{+, \b} d_0^{-1} \bH_{-,\b^\ast } + \bH_{ +, \a} a_0^{-1} \bH_{-,\a^\ast}\\
&=    I_{\sX_+ } - \bR_{11}.
\end{align*}
where the last equality follows from \eqref{formulaR11}. We proved $ \bR_{11} + \bH_{+,g} \bR_{21} = I_{\sX_+}$.

Next we will prove that $\bR_{12} + \bH_{+,g} \bR_{22} =0 $. Using  \eqref{defR22} and the identities in \eqref{ident6c} we obtain
\begin{align*}
\bH_{+,g} \bR_{22} &=
\bH_{+,g} \bT_{-,\d } d_0^{-1} \bT_{-,\d^\ast} - \bH_{+,g} \bT_{-,\c} a_0^{-1} \bT_{-,\c^\ast}\\
&= - \bH_{+, \b} d_0^{-1} \bT_{-,\d^\ast} + \bH_{+,\a } a_0^{-1} \bT_{-,\c^\ast}\\
&= - \bR_{12}.
\end{align*}
We proved that $ \bR_{12} + \bH_{+,g}\bR_{22} =0 $.

The equalities $ \bH_{-,g^\ast} \bR_{11} + \bR_{21} =0 $ and
$ \bH_{-,g^\ast} \bR_{12} + \bR_{22} = I_{\sY_-} $
are proved in a similar way.

\medskip
\noindent\textsc{Part 3.}
To finish the proof we note that \eqref{RLinvers} and \eqref{RRinvers} imply that the operator $\boldsymbol\Omega$ is invertible and that its inverse   is given by
\[
\boldsymbol\Omega^{-1}=\begin{bmatrix} \bR_{11} & \bR_{12} \\ \bR_{21} &  \bR_{22} \end{bmatrix},
\]
which completes the proof. \epr

\medskip
As mentioned in the introduction,  Theorem \ref{thm:inversion1} has been many predecessors.
See also Sections \ref{sec:wienerL} and  \ref{sec:wienerC}.

\begin{exa}\label{exa:not-inv}
We conclude this section with an example of the type announced in Remark \ref{rem:not-iff}, i.e.,   the operator  $ \boldsymbol{\Omega} $
is invertible and item (a) in Theorem \ref{thm:inversion1} is not satisfied.
We use a special case of the example in Subsection \ref{subsec:Mat}.
Let $ p = 2 $ and
\[
g = \begin{bmatrix} 1 & 1 \\ 0 & 1 \end{bmatrix}.
\]
Then $ \boldsymbol{\Omega} $ is invertible.
To see this we choose bases for the upper triangular and the lower triangular matrices
and determine the matrix of $ \boldsymbol{\Omega} $ with respect to these bases.
The basis we choose for the upper triangular matrices is
\[
\left\{
\begin{bmatrix} 1 & 0 \\ 0 & 0 \end{bmatrix} , \quad
\begin{bmatrix} 0 & 1 \\ 0 & 0 \end{bmatrix} , \quad
\begin{bmatrix} 0 & 0 \\ 0 & 1 \end{bmatrix} \right\}
\]
and the basis we choose for the lower triangular matrices is
\[
\left\{
\begin{bmatrix} 1 & 0 \\ 0 & 0 \end{bmatrix}, \quad
\begin{bmatrix} 0 & 0 \\ 1 & 0 \end{bmatrix}, \quad
\begin{bmatrix} 0 & 0 \\ 0 & 1 \end{bmatrix} \right\}.
\]
Then it follows that the matrix for $ \boldsymbol\Omega $ with respect to these bases is
\[
\begin{bmatrix}
1 & 0 & 0 & 1 & 1 & 0 \\
0 & 1 & 0 & 0 & 0 & 1 \\
0 & 0 & 1 & 0 & 0 & 1 \\
1 & 0 & 0 & 1 & 0 & 0 \\
1 & 0 & 0 & 0 & 1 & 0 \\
0 & 1 & 1 & 0 & 0 & 1
\end{bmatrix},
\]
which is an invertible matrix. The solution of the two equations \eqref{Omega-basic-eqs} is
\[
\a = \begin{bmatrix} -1 & -1 \\ 0 & 0 \end{bmatrix} , \quad
\c = \begin{bmatrix} 1 & 0 \\ 1 & 1 \end{bmatrix} , \quad
\b = \begin{bmatrix} 1 & 1 \\ 0 & 1 \end{bmatrix} , \quad
\d = \begin{bmatrix} 0 & 0 \\ -1 & -1 \end{bmatrix} .
\]
We see that $ \a $ and $ \d $ are not invertible and then
the diagonals $ \a_d $ and $ \d_d $ are also  not invertible.
It is also easy to check that $ \a, \b, \c $ and $ \d $ satisfy the inclusions
\eqref{inclu12} and \eqref{inclu34}.
\end{exa}

\setcounter{equation}{0}
\section{Solution to the abstract  twofold  EG inverse problem }\label{sec:solinverseprobl}
The next theorem is the main result of this section.

\begin{thm} \label{thm:mainthm}
Let $\a \in \sA_+ $, $\b \in \sB_+ $, $\c \in \sC_-$, $\d \in \sD_- $, and assume that
\begin{itemize}
\item[\textup{(a)}] $ a_0 = P_{\sA_d} \a $ and $ d_0 = P_{\sD_d} \d $ are invertible in $ \sA_d $ and $ \sD_d $, respectively;
\item[\textup{(b)}] conditions \textup{ (C1)--(C6) } are satisfied.
\end{itemize}
Furthermore, let $ \bR_{11} $,  $ \bR_{12} $, $ \bR_{21} $, $ \bR_{22} $ be the operators defined by
\eqref{defR11}--\eqref{defR22}. Then the twofold  EG inverse problem associated with the data set $\{ \a, , \b, \c, \d\}$  has a solution if and only if
\begin{itemize}
\item[\textup{(i)}]
$ \bR_{11} : \sX_+ \to \sX_+ $ and $ \bR_{22} : \sY_- \to \sY_- $ are invertible;
\item[\textup{(ii)}]
$ \left( \bR_{11}^{-1} \b \right)^\ast = \bR_{22}^{-1} \c $;
\item[\textup{(iii)}]
$ \bR_{11}^{-1} \bR_{12} = \bH_{+,\rho} $ for some $\rho \in \sB $ and
$ \bR_{22}^{-1} \bR_{21} = \bH_{-,\eta} $ for some $ \eta \in \sC $.
\end{itemize}
In that case the solution $ g $ of the twofold   EG inverse problem associated with
$ \a$, $ \b $, $ \c $ and $ \d $ is unique and is given by
\begin{equation}\label{formula-gg*}
  g = - \bR_{11}^{-1} \b =-(\bR_{22}^{-1} \c)^*.
\end{equation}
\end{thm}

\bpr
The proof is divided into two parts.  Note that the uniqueness statement is already covered by Corollary \ref{cor:uniq}.  In the first part of the proof we prove the necessity of the conditions (i), (ii), (iii).

\smallskip
\noindent\textsc{Part 1.} Assume $g\in \sB_+$ is a solution to the  twofold  EG inverse problem associated with the data set $\{ \a,  \b, \c, \d\}$. Note that conditions  (a) and (b) in Theorem \ref{thm:mainthm} imply conditions (a) and (b) in Theorem \ref{thm:inversion1}. Furthermore, from Corollary \ref{cor:g-solinv} we know that the identities in  \eqref{Omega-basic-eqs} are satisfied.  Thus Theorem \ref{thm:inversion1} tells us that  operator  $\bOm$ defined by \eqref{defZg1} is invertible and its inverse is given by  \eqref{IGGIinverse}. In particular,  the operator $\bR$ defined by
\[
\bR= \begin{bmatrix} \bR_{11} &  \bR_{12} \\ \bR_{21} &  \bR_{22} \end{bmatrix}
\]
is invertible. But then    the second part of Lemma \ref{lem:inverseR} tells us that the operators $\bR_{11}$ and  $\bR_{22}$ are invertible, i.e., condition (i) is fulfilled.  Furthermore, again using  the second part of Lemma \ref{lem:inverseR},
we have
\[
\bOm=\bR^{-1}=\begin{bmatrix} \bI_{\sX_+} & -\bR_{12} \bR_{22}^{-1} \\ -\bR_{21} \bR_{11}^{-1} & \bI_{\sY_-} \end{bmatrix}=
\begin{bmatrix} \bI_{\sX_+} & -\bR_{11}^{-1} \bR_{12}  \\ -\bR_{22}^{-1} \bR_{21}  & \bI_{\sY_-} \end{bmatrix}.
\]
In particular, we have
\begin{equation*}
\bH_{+,g}=-  \bR_{11}^{-1}  \bR_{12} =
- \bR_{12}  \bR_{22}^{-1}, \quad \bH_{-,g^*}=-  \bR_{21} \bR_{11}^{-1}   =
- \bR_{22}^{-1}\bR_{21}.
\end{equation*}
The preceding two identities show that    item (iii) holds with $\rho=-g$  and $\eta=-g^*$. Finally, since $\rho=-g$  and $\eta=-g^*$,  the identities in  \eqref{formula-gg*} imply that item (ii) is satisfied.

\smallskip
\noindent\textsc{Part 2.} In this part we assume that conditions (i), (ii), (iii)  are satisfied and we show that the twofold  EG inverse problem associated with the data set $\{ \a,  \b, \c, \d\}$ has a solution.

Put $g=-P_{\sB_+}\rho$ and $h=-P_{\sC_-}\eta$. We shall show that $h=g^*$ and for this choice of $g$ the inclusions \eqref{inclu12} and \eqref{inclu34} are fulfilled. Note that $P_{\sB_+}(g+\rho)=0$, so that $g+\rho\in\sB_-$. From the second part of  \eqref{Hmp_of_Bpm} we then obtain that  $\bH_{+, g}=-\bH_{+, \rho}$, and, by a similar argument, from the first part of \eqref{Hmp_of_Cpm} it follows that $\bH_{-, h}=-\bH_{-, \eta}$. Using these identities  and  those given by Lemma \ref{lem:Hunit} together with the second and third identity in \eqref{unitsR} we see that condition  (iii) yields
\begin{align*}
g&= \bH_{+, g}e_\sD=-\bH_{+, \rho}e_\sD=-\bR_{11}^{-1} \bR_{12}e_\sD=-\bR_{11}^{-1}\b,\\
h&= \bH_{-, h}e_\sA=-\bH_{-, \eta}e_\sA=-\bR_{22}^{-1} \bR_{21}e_\sA=-\bR_{22}^{-1} \c.\end{align*}
But then  (ii)  implies that $h=g^*$. Furthermore, (iii) tells us that
\begin{equation}
\label{HgHg*7}
\bR_{11}^{-1} \bR_{12}=-\bH_{+,g} \ands \bR_{22}^{-1} \bR_{21}=-\bH_{-,g^*}.
\end{equation}

According to Lemma \ref{lem:inverseR}  condition (i) implies that the operator $\bR$ given by \eqref{invform1} is invertible, and its inverse is given by \eqref{invform2}. This together with  the identities in \eqref{HgHg*7} implies that
\[
{\bR^{-1}}=\begin{bmatrix} \bR_{11} & \bR_{12} \\ \bR_{21} &  \bR_{22}  \end{bmatrix}^{-1}= \begin{bmatrix} I_{\sX_+} & -\bR_{11}^{-1} \bR_{12}  \\ -\bR_{22}^{-1} \bR_{21}  & I_{\sY_-} \end{bmatrix}= \begin{bmatrix} I_{\sX_+} & \bH_{+g}  \\ \bH_{-,g^*}  & I_{\sY_-} \end{bmatrix}.
 \]
Next note that the identities in  \eqref{unitsR} can be rephrased as
\[
\begin{bmatrix} \bR_{11} & \bR_{12} \\ \bR_{21} &  \bR_{22}  \end{bmatrix}
\begin{bmatrix} e_\sA & 0 \\ 0 &   e_\sD \end{bmatrix}=
\begin{bmatrix} \a &\b  \\ \c &   \d \end{bmatrix}.
\]
But then
 \[
 \begin{bmatrix} I_{\sX_+} & \bH_{+g}  \\ \bH_{-,g^*}  & I_{\sY_-} \end{bmatrix}
 \begin{bmatrix} \a &\b  \\ \c &   \d \end{bmatrix}=\begin{bmatrix} e_\sA & 0 \\ 0 &   e_\sD \end{bmatrix},
\]
and  the equivalences in \eqref{inclu12a} and   \eqref{inclu34a} tell us that with our choice of $g$ the inclusions \eqref{inclu12} and \eqref{inclu34} are fulfilled. Hence  $g$ is a solution  to the twofold  EG inverse problem associated with the data set $\{ \a,  \b, \c, \d\}$.
Since $g=- \bR_{11}^{-1} \b$, the proof is complete.
\epr

\medskip

A variation on condition (iii) in Theorem \ref{thm:mainthm} does not appear in the solution to the twofold EG inverse problem in $L^1(\BR)$ as formulate in the introduction, e.g., in Theorem \ref{thm:mainL2}, and neither in the solution to the discrete twofold EG inverse problem in \cite{tHKvS1}. This is because in the abstract setting presented in this paper we do not have a characterization of Hankel-type operators via an intertwining condition as in the discrete case as well as in the continuous case (where an extra condition is needed, as shown in the Appendix). Lemma \ref{isHankel} below provides, at the abstract level, a result that will be useful in proving that condition (iii) is implied by the assumptions made for the special cases we consider.

\medskip
Assume we have operators ${\bV_{\sZ,\pm}} : \sZ_\pm \to \sZ_\pm $
and ${\bV_{\ast,\sZ,\pm}}: \sZ_\pm \to \sZ_\pm $, with $ \sZ $ either $ \sX $ or $\sY $,
that are such that $ \bV_{\ast,\sZ,\pm} \bV_{\sZ,\pm} = I_{\sZ_\pm} $ and
\begin{align*}
  \mbox{for any $ \phi \in \sA$:}\quad
& \bV_{\ast,\sX,\pm} \bH_{\pm,\phi} = \bH_{\pm,\phi}   \bV_{\sX,\mp} ,\quad
\bV_{\ast,\sX,\pm} \bT_{\pm,\phi} \bV_{\sX,\pm} = \bT_{\pm,\phi} ;\\
  \mbox{for any $ \phi \in \sB $:}\quad
& \bV_{\ast,\sX,\pm} \bH_{\pm,\phi} = \bH_{\pm,\phi}  \bV_{\sY,\mp}, \quad
\bV_{\ast,\sX,\pm} \bT_{\pm,\phi} \bV_{\sY,\pm} = \bT_{\pm,\phi} ; \\
  \mbox{for any $ \psi \in \sD $:}\quad
& \bV_{\ast,\sY,\pm} \bH_{\pm,\psi} = \bH_{\pm,\psi}  \bV_{\sY,\mp}, \quad
\bV_{\ast,\sY,\pm} \bT_{\pm,\phi} \bV_{\sY,\pm} = \bT_{\pm,\phi} ; \\
  \mbox{for any $\psi \in \sC $:}\quad
& \bV_{\ast,\sY,\pm} \bH_{\pm,\psi} = \bH_{\pm,\psi}  \bV_{\sX,\mp}, \quad
\bV_{\ast,\sY,\pm} \bT_{\pm,\phi} \bV_{\sX,\pm} = \bT_{\pm,\phi} ,
\end{align*}
and
\begin{align*}
\mbox{for any $ \phi \in \sA_\pm $:}\quad
&  \bT_{\pm,\phi} \bV_{\sX,\pm} = \bV_{\sX,\pm} \bT_{\pm,\phi};\\
\mbox{for any $ \phi \in \sB_\pm $:}\quad
&  \bT_{\pm,\phi} \bV_{\sY,\pm} = \bV_{\sX,\pm} \bT_{\pm,\phi};\\
\mbox{for any $ \psi \in \sC_\pm $:}\quad
& \bT_{\pm,\psi} \bV_{\sY,\pm} = \bV_{\sY,\pm} \bT_{\pm,\psi};\\
\mbox{for any $ \phi \in \sB_\pm $:}\quad
& \bT_{\pm,\psi} \bV_{\sX,\pm} = \bV_{\sY,\pm} \bT_{\pm,\psi}.
\end{align*}

\begin{lem}\label{isHankel}
With $ \bR_{ij} $ defined as above one has the equalities
\[
\bR_{11} \bV_{\ast,\sX,+} \bR_{12} = \bR_{12} \bV_{\sY,-} \bR_{22}
\ands  \bR_{22} \bV_{\ast,\sY,-} \bR_{21} = \bR_{21} \bV_{\sX,+} \bR_{11} .
\]
Moreover, if $ \bR_{11} $ and $ \bR_{22} $ are invertible, then
\begin{equation}\label{eq:intertw12-21}
\bV_{\ast,\sX,+} \bR_{12}  \bR_{22}^{-1}= \bR_{11}^{-1} \bR_{12} \bV_{\sY,-}.
\end{equation}
\end{lem}

\bpr
First we will prove that $ \bR_{11} \bV_{\ast,\sX,+} \bR_{12} = \bR_{12} \bV_{\sY,-} \bR_{22} $.
We start with deriving the equality
\begin{equation}\label{matrixTHSHT}
\begin{bmatrix}  \bT_{+,\a^\ast} \\ \bT_{+,\b^\ast}\end{bmatrix} \bV_{\ast,\sX,+}
\begin{bmatrix} \bH_{+,\b} & \bH_{+,\a} \end{bmatrix} =
\begin{bmatrix}  \bH_{+,\c^\ast} \\ \bH_{+,\d\ast}\end{bmatrix} \bV_{\sY,-}
\begin{bmatrix} \bT_{-,\d} & \bT_{-,\c} \end{bmatrix}.
\end{equation}
To obtain \eqref{matrixTHSHT}, first notice that
\begin{align*}
 \bV_{\ast,\sX,+} \begin{bmatrix} \bH_{+,\b} & \bH_{+,\a} \end{bmatrix}
 &= \begin{bmatrix} \bH_{+,\b}  & \bH_{+,\a}  \end{bmatrix}
 \begin{bmatrix} \bV_{\sY,-} & 0  \\ 0 &  \bV_{\sX,-} \end{bmatrix};\\
 \bV_{\sY,-}
\begin{bmatrix} \bT_{-,\d} & \bT_{-,\c} \end{bmatrix}
&= \begin{bmatrix} \bT_{-,\d} & \bT_{-,\c} \end{bmatrix}
\begin{bmatrix} \bV_{\sY,-} & 0  \\ 0 &  \bV_{\sX,-} \end{bmatrix}.
\end{align*}
Then use \eqref{matrixTHHT} to get that
\begin{align*}
\begin{bmatrix}  \bT_{+,\a^\ast} \\ \bT_{+,\b^\ast} \end{bmatrix} \bV_{\ast,\sX,+}
\begin{bmatrix} \bH_{+,\b} & \bH_{+,\a} \end{bmatrix}
& =
{\begin{bmatrix}  \bT_{+,\a^\ast} \\ \bT_{+,\b^\ast}\end{bmatrix}
\begin{bmatrix} \bH_{+,\b} & \bH_{+,\a} \end{bmatrix}
\begin{bmatrix} \bV_{\sY,-} & 0  \\ 0 &  \bV_{\sX,-} \end{bmatrix}} \\
& =
\begin{bmatrix}  \bH_{+,\c^\ast} \\ \bH_{+,\d\ast}\end{bmatrix}
\begin{bmatrix} \bT_{-,\d} & \bT_{-,\c} \end{bmatrix}
\begin{bmatrix} \bV_{\sY,-} & 0  \\ 0 &  \bV_{\sX,-} \end{bmatrix} \\
& =
\begin{bmatrix}  \bH_{+,\c^\ast} \\ \bH_{+,\d\ast}\end{bmatrix} \bV_{\sY,-}
\begin{bmatrix} \bT_{-,\d} & \bT_{-,\c} \end{bmatrix}.
\end{align*}
We proved \eqref{matrixTHSHT}. By multiplying \eqref{matrixTHSHT} on the left  and the right by
\[
\begin{bmatrix}  \bT_{+,\a} & \bT_{+,\b} \end{bmatrix}  \quad \mbox{and \ }
\begin{bmatrix} \bT_{-,\d^\ast } \\ \bT_{-,\c^\ast } \end{bmatrix},
\]
respectively, one gets $ \bR_{11} \bV_{\ast,\sX,+} \bR_{12} = \bR_{12} \bV_{\sY,-} \bR_{22} $.
Furthermore, the  equality $ \bR_{22} \bV_{\ast,\sY,-} \bR_{21} = \bR_{21} \bV_{\sX,+} \bR_{11} $ can be proved in a similar way.

Given the invertibility of $\bR_{11}$ and $\bR_{22}$ the preceding two identities yield the identity \eqref{eq:intertw12-21}  trivially.
\epr

\setcounter{equation}{0}
\section{Proof of Theorems \ref{thm:mainL1} and \ref{thm:mainL2} }\label{sec:wienerL}

In this section we will prove Theorems \ref{thm:mainL1} and   \ref{thm:mainL2}. Recall that in this  case the data are given by \eqref{eq:data1} and  \eqref{eq:data2}, and the twofold EG inverse problem is to find $g\in L^1(\BR_+)^{p\ts q}$  such that \eqref{eq:probl1} and  \eqref{eq:probl2} are satisfied.

As a  first step,  the  above problem will be put into the  general setting  introduced in Section \ref{sec:setting} using  a particular choice for $ \sA$, $ \sB$, $\sC$, $ \sD$, namely as follows:
\begin{align}
&\sA=\{f \mid f =\eta e_p +f _0,\mbox{ where }  \eta \in \BC^{p \times  p},\  f_0\in L^1(\BR)^{p \times p}\},\label{def:A8}\\
&\sB=L^1(\BR)^{p \times q}, \qquad  \sC=L^1(\BR)^{q \times p},\label{def:BC8}\\
&\sD=\{h\mid h=\zeta e_q +h_0,\mbox{ where }  \zeta \in \BC^{q\ts q},\  h_0\in L^1(\BR)^{q \times q}\}.\label{def:D8}
\end{align}
Furthermore, $ \sA$, $ \sB$, $\sC$, $ \sD $  admit  decompositions as in \eqref{decomp1} and \eqref{decomp2} using
\begin{align*}
& \sA_+^0 = L^1(\BR_+)^{p \times p},  \quad \sA_-^0 = L^1(\BR_-)^{p \times p},  \quad
\sA_d = \{\eta e_p \mid \eta \in \BC^{p \times p} \}, \\
& \sB_+ = L^1(\BR_+)^{p \times q},  \quad \sB_- = L^1(\BR_-)^{p \times q},  \\
& \sC_+ = L^1(\BR_+)^{q \times p},  \quad \sC_- = L^1(\BR_-)^{q \times p},  \\
& \sD_+^0 = L^1(\BR_+)^{q \times q},  \quad \sD_-^0 = L^1(\BR_-)^{q \times q},  \quad
\sD_d = \{ \zeta e_q \mid \zeta \in \BC^{q \times q} \}.
\end{align*}
Here $e_m$, for $m=p,q$,
is the constant $ m \times m $ matrix function on $\BR$ whose value is the $m\ts m$ identity matrix $I_m$.  Thus given $  \eta \in \BC^{m \times  m}$, the  symbol $ \eta e_m $  denotes the  constant matrix function on $\BR$ identically equal to $\eta$.

We proceed by defining the algebraic structure. The addition is the usual addition of functions and is denoted by $+$. For the product we use the symbol $\diam$ which in certain cases is just the usual  convolution product $\star$.  If $f=\eta_f e_p+ f_0 \in \sA$ and  $\tilde{f} = \eta_{\tilde{f}}e_p+ \tilde{f}_0\in \sA$,  then the $\diamond$ product  is defined by
\[
 f \diamond \tilde{f} := \eta_f \eta_{\tilde{f}}e_p   + \big(\eta_f  \tilde{f}_0   +  f_0\eta_{\tilde{f}} e_p+  f_0 \star  \tilde{f}_0\big.
 \]
Thus for  $ f \in L^1(\BR)^{n \times m} $ and $ h \in L^1(\BR)^{m \times k} $  the product   $f\diamond h$ is the  convolution product $ f \star h $.
The product of elements  $f= \eta_f e_p  + f_0 \in \sA $ and $ h_0 \in \sB $ is defined as $ f \diam  h_0 = \eta_f h_0 + f_0 \star h_0 $. Other products are defined likewise.
One only needs the matrix dimension to allow the multiplication.
The units  $ e_\sA $ and $ e_\sD $   in $ \sA $ and  $\sD $ are given by  $ e_\sA =e_p $ and $ e_\sD=e_q  $, respectively.  Finally, the adjoint  $ f^\ast $ for $ f \in L^1(\BR)^{r \times s}$  is defined by   $ f^\ast(\l) = f(-\l)^{\ast}$, $\l\in \BR$, so that $f^\ast \in L^1(\BR)^{s \times r} $.  For $ f = \eta e_s+f_0$   with  $ \eta\in \BC^{s\ts s}$ and $f_0\in L^1(\BR)^{s \times s} $ we define $f^\ast$ by  $ f^\ast ={\eta}^\ast e_s+f_0^\ast $, where $\eta^\ast$ is the adjoint  of the matrix  $\eta$.  It easily follows that all conditions of the first paragraph of Section \ref{sec:setting} are satisfied. We conclude that $\sM_{\sA,\sB,\sC,\sD}$ is admissible.

\begin{rem}\label{rem:coincide} Observe that for a data set $\{a, b,c,d \}$ as in
\eqref{eq:data1} with $\a$, $\b$, $\c$, $\d$ the  functions given by \eqref{def:abcd-greek},
the inclusions for $a$, $b$, $c$ and $d$ in \eqref{eq:probl1} and \eqref{eq:probl2} are
equivalent to the inclusions \eqref{inclu12} and \eqref{inclu34} for $\a$, $\b$, $\c$ and $\d$.
Thus the solutions $g\in L^1(\BR_+)^{p\times q}$ for the twofold EG inverse problem
formulated in the introduction  coincide with the solutions
of the abstract twofold EG inverse problem of Section~\ref{sec:setting} using
the specification given in the present section.
Furthermore, in this case
\begin{align}
&\textup{(C1)} \quad \Longleftrightarrow \quad \a^\ast \diamond \a - \c^\ast \diamond \c = e_p;\label{eq:L1a}\\
&\textup{(C2)} \quad \Longleftrightarrow \quad  d^\ast \diamond \d - \b^\ast \diamond \b = e_q; \label{eq:L1b}\\
&\textup{(C3)} \quad \Longleftrightarrow \quad \a^\ast \diamond \b = \c^\ast \diamond \d. \label{eq:L1c}
\end{align}
Thus  (C1)--(C3) are satisfied if and only if the following three  identities  hold true:
\begin{equation}\label{eq:L1abc}
\a^\ast \diamond \a - \c^\ast \diamond \c = e_p, \quad d^\ast \diamond \d - \b^\ast \diamond \b = e_q, \quad\a^\ast \diamond \b = \c^\ast \diamond \d.
\end{equation}
\end{rem}

\subsection{Proof of Theorem \ref{thm:mainL1}}
Note that Theorem  \ref{thm:mainL1} is an ``if and only if '' theorem. We first proof the ``only if'' part.  Let $g \in L^1(\BR_+)^{p\ts q}$, and assume that  the operator $W$ given by  \eqref{defZ1} is invertible. Note that
\[
\begin{bmatrix}0\\[.1cm]-g^* \end{bmatrix}\in \begin{bmatrix} L^1(\BR_+)^{p\ts p}\\[.1cm]L^1(\BR_+)^{q\ts p}\end{bmatrix}\ands
\begin{bmatrix}-g\\[.1cm]0 \end{bmatrix}\in \begin{bmatrix} L^1(\BR_+)^{p\ts q}\\[.1cm]
L^1(\BR_+)^{q\ts q}\end{bmatrix}
\]
Since $W$ is invertible, we see that there exist
\begin{equation*}
a\in L^1(\BR_+)^{p \ts p}, \quad c\in L^1(\BR_-)^{q \ts p},\quad  b\in L^1(\BR_+)^{p \ts q},  \quad d\in L^1(\BR_-)^{q \ts q}
\end{equation*}
such that
\begin{equation*}
W\begin{bmatrix} a\\ c\end{bmatrix}=\begin{bmatrix}0\\-g^* \end{bmatrix}\ands
W\begin{bmatrix} b\\ d\end{bmatrix}=\begin{bmatrix}-g\\0 \end{bmatrix}.
\end{equation*}
But this implies that $g$ is a solution to the twofold    EG inverse problem defined by the   data set $\{a, b,c,d \}$. Thus the ``only if''  part of Theorem  \ref{thm:mainL1} is proved.

\medskip
Next we  prove the ``if''  part of Theorem~\ref{thm:mainL1}.  We assume that $g \in L^1(\BR_+)^{p\ts q}$  is a solution to the  twofold    EG inverse problem defined by  the data set $\{a, b,c,d \}$ given by \eqref{eq:data1} and \eqref{eq:data2}. Furthermore, $\a$, $\b$, $\c$, and $\d$ are given by \eqref{def:abcd-greek}, and $\sM=\sM_{\sA,\sB,\sC,\sD}$  is the admissible algebra defined in the beginning of this section. Our aim is to obtain the ``if''  part of Theorem~\ref{thm:mainL1} as a corollary of Theorem \ref{thm:inversion1}.  For that purpose    various results of  Section \ref{sec:setting} and Sections \ref{sec:preliminaries}--\ref{sec:inversionthm} have to be specified further  for the case when $ \sA$, $ \sB$, $\sC$, $ \sD$ are given by \eqref{def:A8}--\eqref{def:D8} in the beginning of this section.  This will be done in four  steps.

\smallskip
\noindent
\textsc{Step 1. Results  from Section \ref{sec:setting}.} Since $g\in L^1(\BR_+)^{p\ts q}=\sB_+$  is a solution to the  twofold EG inverse problem   associated with  the data $\{\a, \b, \c, \d\}$, we know from  Proposition \ref{prop:cond123} that conditions (C1)-(C3) are satisfied.   Furthermore,
\begin{equation}\label{eq:a0-d0}
a_0=P_{\sA_d}\a=e_p \ands d_0=P_{\sD_d}\d=e_q.
\end{equation}
But then the fact that  $\a$, $\b$, $\c$, and $\d$ are matrix functions implies that conditions (C4)-(C6) are also satisfied. Indeed,  using the identities in \eqref{eq:a0-d0}, we see from  \eqref{C1-3matrix}  that
\begin{equation}\label{eq:C1-3matrix2a}
\begin{bmatrix} \a^\ast(\l) & \c^\ast(\l)   \\ \b^\ast(\l)  & \d^\ast(\l) \end{bmatrix}
\begin{bmatrix}I_p & 0 \\ 0 & -I_q \end{bmatrix}
\begin{bmatrix} \a(\l)  & \b(\l) \\ \c(\l) & \d(\l) \end{bmatrix}=
\begin{bmatrix} I_p & 0 \\ 0 & - I_q \end{bmatrix}, \quad \l\in \BR.
\end{equation}
In particular,  the first matrix  in the left hand side  of  \eqref{eq:C1-3matrix2a}  is surjective and third matrix in the left hand side  of  \eqref{eq:C1-3matrix2a} is injective. But all   matrices in  \eqref{eq:C1-3matrix2a} are  finite square matrices. It  follows   that all these matrices are invertible. Hence
\[
\begin{bmatrix} \a  & \b \\ \c & \d \end{bmatrix}^{-1}
\begin{bmatrix} e_p & 0 \\ 0 & -e_q \end{bmatrix}
\begin{bmatrix} \a^\ast & \c^\ast  \\ \b^\ast & \d^\ast \end{bmatrix}^{-1}=\begin{bmatrix} e_p & 0 \\ 0 & -e_q \end{bmatrix},
\]
which yields
\[
\begin{bmatrix} \a  & \b \\ \c & \d \end{bmatrix}
\begin{bmatrix} e_p & 0 \\ 0 & -e_q \end{bmatrix}
\begin{bmatrix} \a^\ast & \c^\ast  \\ \b^\ast & \d^\ast \end{bmatrix}=
\begin{bmatrix} e_p & 0 \\ 0 & -e_q \end{bmatrix}.
\]
The latter implies  that conditions (1.4)--(1.6) are satisfied. In particular, we have proved that
\begin{itemize}
\item[\textup{(i)}]   $ a_0 = P_{\sA_d} \a $ and $ d_0 = P_{\sD_d} \d $ are invertible in $ \sA_d $ and $ \sD_d $, respectively;
\item[\textup{(ii)}]  conditions  (C1)--(C6)  are satisfied.
\end{itemize}

\smallskip
\noindent\textsc{Step 2. Results  from Section \ref{sec:preliminaries}.}   In the present context the spaces $\sX$ and $\sY$, $\sX_+$ and $\sY_+$, and $\sX_-$ and $\sY_-$  defined in the first paragraph of  Section~\ref{sec:preliminaries} are given  by
\begin{align*}
\sX=\sA\dot{+}\sB&=\left(\BC^{p\ts p}e_p +L^1(\BR)^{p\ts p}\right) \dot{+}  L^1(\BR)^{p\ts q}, \\
\sX_+= \sA_+\dot{+}\sB_+&=\left(\BC^{p\ts p}e_p +L^1(\BR_+)^{p\ts p}\right) \dot{+}  L_+^1(\BR)^{p\ts q},\\
\sX_-=\sA_-^0\dot{+}\sB_-&=L^1(\BR_-)^{p\ts p}  \dot{+}  L^1(\BR_-)^{p\ts q},
\end{align*}
and
\begin{align*}
\sY=\sC \dot{+} \sD&= L^1(\BR)^{q\ts p} \dot{+} \left(L^1(\BR)^{q\ts q} +
\BC^{q\ts q}e_q\right),\\
\sY_+=\sC_+ \dot{+} \sD_+^0&=L^1(\BR_+)^{q\ts p} \dot{+} L^1(\BR_+)^{q\ts q},\\
\sY_-= \sC_- \dot{+}\sD_-&=  L^1(\BR_-)^{q\ts p} \dot{+}
\left(L^1(\BR_-)^{q\ts q} +\BC^{q\ts q}e_q\right).
\end{align*}
In the sequel we write $x\in  \sX$ as $ x = ( f, g )$, where $ f =\eta_f e_p + f_0 \in \sA $ and $ g  \in \sB $.  In a similar way  vectors $x_+ \in \sX_+$ and $x_- \in \sX_-$  will be written as
\begin{align*}
x_+&=( f_+, g_+ ), \hspace{.1cm}\mbox{where $f_+= \eta_{f_+} e_p + f_{+,0} \in \sA_+$ and $g_+\in \sB_+$},\\
x_-&=( f_-, g_- ), \hspace{.1cm}\mbox{where $f_-\in \sA_-^0$ and $g_-\in \sB_-$}.
\end{align*}
Analogous notations will be used  for vectors $y\in \sY$, $y_+\in \sY_+$, and $y_-\in \sY_-$.  Indeed, $y\in \sY$ will be written as $(h,k)$, where $h\in \sC$, and $k= \zeta e_q+k_0\in \sD$, and
\begin{align*}
y_+&=(h_+,k_+), \hspace{.1cm}\mbox{where $h_+\in \sC_+$ and  $k_+\in\sD_+^0$},\\
y_-&=(h_-,k_-), \hspace{.1cm}\mbox{where $h_-\in \sC_-$  and  $k_-=\zeta_{k_-}e_q +k_{-,0}\in  \sD_-$}.
\end{align*}
Furthermore, in what follows  $0_{p \times q} $ and $0_{q \times p} $ denote  the linear spaces consisting only  of  the zero $p\ts q$  and  zero $q\ts p$ matrix, respectively.

Using the above notation  we define the following operators:
\begin{align*}
&J_{\sX_+} : \sX_+ \to \begin{bmatrix} \BC^{p \times p} \dot{+} 0_{p \times q}
\\[.15cm] L^1(\BR_+)^{p \times p} \dot{+} L^1(\BR_+)^{p \times q}  \end{bmatrix},
\quad J_{\sX_+} x_+ = \begin{bmatrix} ( \eta_{f_+ } , 0) \\ ( f_{+,0} , g_+) \end{bmatrix},\\[.2cm]
& J_{\sX_-} : \sX_- \to
L^1(\BR_-)^{p \times p} \dot{+} L^1(\BR_-)^{p \times q},
\quad J_{\sX_-} x_- = ( f_- , g_-  ) ,
\end{align*}
and
\begin{align*}
&J_{\sY_+} : \sY_+ \to
 L^1(\BR_+)^{q \times p} \dot{+} L^1(\BR_+)^{q \times q} ,
\quad J_{\sY_+} y_+ =  (h_+ , k_+ ),\\[.2cm]
&J_{\sY_-} : \sY_- \to \begin{bmatrix} 0_{q \times p} \dot{+} \BC^{q \times q}
\\[.15cm]  L^1(\BR_-)^{q \times p} \dot{+} L^1(\BR_-)^{q \times q}  \end{bmatrix},
\quad J_{\sY_-} y_- =  \begin{bmatrix} ( 0 , \zeta_{k_-} )\\ (h_- , k_{-,0}) \end{bmatrix}.
\end{align*}
Note that all four operators defined above are invertible operators.

Next, in our present setting where $\a, \b, \c,\d$ are given by given by \eqref{def:abcd-greek},  we relate the Toeplitz-like and Hankel-like operators introduced in  Section \ref{sec:preliminaries} to ordinary Wiener-Hopf and  Hankel integral operators.

\smallskip\noindent
Let $ \a = a_+ + e_p + a_-  \in \sA_+ + \sA_d + \sA_- = \sA $. Then
\begin{align*}
& J_{\sX_+} \bT_{+,\a} = \begin{bmatrix} I_p & 0 \\ a_+ &  T_{+,\a} \end{bmatrix} J_{\sX_+} ,
\quad
J_{\sX_-} \bT_{-,\a} = T_{-,\a} J_{\sX_-}, \\
& J_{\sX_+} \bH_{+,\a} = \begin{bmatrix} 0 \\  H_{+,\a} \end{bmatrix} J_{\sX_-},
\quad
J_{\sX_-} \bH_{-,\a}  = \begin{bmatrix} a_- & H_{-,\a} \end{bmatrix} J_{\sX_+}.
\end{align*}
For $ \b = b_+ + b_- \in \sB_+ + \sB_- = \sB $ we get
\begin{align*}
& J_{\sX_+} \bT_{+,\b} = \begin{bmatrix} 0 \\  T_{+,\b} \end{bmatrix} J_{\sY_+},
\quad
J_{\sX_-} \bT_{-,\b} = \begin{bmatrix} b_- & T_{-,\b} \end{bmatrix} J_{\sY_-},
\\ &
J_{\sX_+} \bH_{+,\b}  = \begin{bmatrix} 0 & 0 \\ b_+ & H_{+,\b} \end{bmatrix} J_{\sY_-},
\quad
J_{\sX_-} \bH_{-,\b} = H_{-,\b} J_{\sY_+}.
\end{align*}
Let $ \c = c_+ + c_- \in \sC_+ + \sC_- = \sC $. We have the  equalities
\begin{align*}
& J_{\sY_-} \bT_{-,\c} = \begin{bmatrix} 0 \\ T_{-,\c} \end{bmatrix} J_{\sX_-},
\quad
J_{\sY_+} \bT_{+,\c} = \begin{bmatrix} c_+ & T_{+,\c} \end{bmatrix} J_{\sX_+},
\\ &
J_{\sY_-} \bH_{-,\c} = \begin{bmatrix} 0 & 0 \\ c_- & H_{-,\c} \end{bmatrix} J_{\sX_+},
\quad
J_{\sY_+} \bH_{+,\c} = H_{+,\c} J_{\sX_-}.
\end{align*}
Let $ \d = d_+ + e_q + d_-  \in \sD_+ + \sD_d + \sD_- = \sD $.
Then
\begin{align*}
& J_{\sY_-} \bT_{-,\d} = \begin{bmatrix} I_p & 0 \\ d_- &  T_{-,\d} \end{bmatrix} J_{\sY_-},
\quad
J_{\sY_+} \bT_{+,\d} = T_{+,\d} J_{\sX_-},
\\ &
J_{\sY_-} \bH_{-,\d} = \begin{bmatrix} 0 \\  H_{-,\d} \end{bmatrix} J_{\sY_+},
\quad
J_{\sY_+} \bH_{+,\d} = \begin{bmatrix} d_+ & H_{+,\d} \end{bmatrix} J_{\sY_-}.
\end{align*}
The following lemma is an immediate consequence of the above relations.

 \smallskip
\begin{lem}\label{WandOmega}
For $ g \in \sB_+ $ one has
\begin{equation}\label{OandW}
\begin{bmatrix} J_{\sX_+} \!\!&\!\! 0  \\ 0 \!\!&\!\! J_{\sY_-} \end{bmatrix}
\begin{bmatrix} \bI_{\sX_+} \!&\! \bH_{+,g} \\ \bH_{-,g^\ast} \!&\! \bI_{\sY_-}  \end{bmatrix}
=
\begin{bmatrix} I_p & 0 & 0 & 0 \\ 0 & I & g & H_{+,g} \\ 0 & 0 & I_q & 0 \\
g^\ast & H_{-,g^\ast} & 0 & I \end{bmatrix}
\begin{bmatrix} J_{\sX_+} \!\!&\!\! 0  \\ 0 \!\!&\!\! J_{\sY_-} \end{bmatrix}.
\end{equation}
Furthermore, if  $\boldsymbol\Omega $ is the operator  defined   by \eqref{defZg1} using the present data,  and if  $W$ is the operator  defined by \eqref{defZ1}, then   \eqref{OandW} shows that $ \boldsymbol\Omega $ is invertible if and only if  $W$ is invertible.
\end{lem}

\medskip
\noindent\textsc{Step 3. Results from Section \ref{sec:prelim}.}
As before we assume that $ a  $, $ b $, $ c $ and $ d $ are given by
\eqref{eq:data1} and \eqref{eq:data2} and $ \a $, $ \ b $, $ \c $ and $ \d $ by
\eqref{def:abcd-greek}, and that $g\in L^1(\BR_+)^{p\ts q}=\sB_+$  is a solution to the  twofold EG inverse problem   associated with  the data $\{\a, \b, \c, \d\}$.  Thus we know from \textsc{Step 1} that
\begin{itemize}
\item[\textup{(i)}]   $ a_0 = P_{\sA_d} \a $ and $ d_0 = P_{\sD_d} \d $ are invertible in $ \sA_d $ and $ \sD_d $, respectively;
\item[\textup{(ii)}]  conditions  (C1)--(C6)  are satisfied.
\end{itemize}
In particular, all conditions underlying the lemmas proved in Section \ref{sec:prelim} are fulfilled.

The following lemma is an immediate consequence of Lemma~\ref{basicids5}.

\begin{lem}\label{lem:basicids9}
Since  conditions \textup{(C1)-(C3) } are satisfied, we have
\begin{equation}\label{matrixTHHT19}
\begin{bmatrix}  T_{+,\a^\ast} \\ T_{+,\b^\ast}\end{bmatrix}
\begin{bmatrix} H_{+,\b} & H_{+,\a} \end{bmatrix} =
\begin{bmatrix}  H_{+,\c^\ast} \\ H_{+,\d^\ast}\end{bmatrix}
\begin{bmatrix} T_{-,\d} & T_{-,\c} \end{bmatrix}.
\end{equation}
and
\begin{equation}\label{matrixTHHT29}
\begin{bmatrix}  T_{-,\d^\ast} \\ T_{-,\c^\ast}\end{bmatrix}
\begin{bmatrix} H_{-,\c} & H_{-,\d} \end{bmatrix} =
\begin{bmatrix}  H_{-,\b^\ast} \\ H_{-,\a^\ast}\end{bmatrix}
\begin{bmatrix} T_{+,\a} & T_{+,\b} \end{bmatrix}.
\end{equation}
\end{lem}

\bpr
The above equalities \eqref{matrixTHHT19} and \eqref{matrixTHHT29} follow from
the equalities \eqref{matrixTHHT} and \eqref{matrixTHHT2} and the representations of the Hankel-like and Toeplitz-like operators given in the paragraph preceding Lemma  \ref{WandOmega}.
For example to prove \eqref{matrixTHHT19} note that
\[
\begin{bmatrix} J_{\sX_+} & 0  \\ 0 & J_{\sY_+} \end{bmatrix}
\begin{bmatrix} \bT_{+,\a^\ast} \\  \bT_{+,\b^\ast} \end{bmatrix} =
\begin{bmatrix} I_p & 0 \\ 0 & T_{+,\a^\ast} \\ 0 & T_{+,\b^\ast} \end{bmatrix}  J_{\sX_+} ,
\]
and
\[
J_{\sX_+} \begin{bmatrix} \bH_{+,\b} &  \bH_{+,\a} \end{bmatrix} =
\begin{bmatrix} 0 & 0 & 0 \\ b & H_{+,\b} & H_{+,\a} \end{bmatrix}
\begin{bmatrix} J_{\sY_-} & 0  \\ 0 & J_{\sX_-} \end{bmatrix} .
\]
On the other hand
\[
\begin{bmatrix} J_{\sX_+} & 0  \\ 0 & J_{\sY_+} \end{bmatrix}
\begin{bmatrix} \bH_{+,\c^\ast} \\  \bH_{+,\d^\ast} \end{bmatrix} =
\begin{bmatrix} 0 & 0 \\ c^\ast & H_{+,\c^\ast} \\ d^\ast & H_{+,\d^\ast} \end{bmatrix} J_{\sY_-} ,
\]
and
\[
J_{\sY_-} \begin{bmatrix} \bT_{-,\d} &  \bT_{-, \c} \end{bmatrix} =
\begin{bmatrix} I_q & 0 & 0 \\ d & T_{-,\d} & T_{-, \c} \end{bmatrix}
\begin{bmatrix} J_{\sY_-} & 0  \\ 0 & J_{\sX_-} \end{bmatrix} .
\]
The equality \eqref{matrixTHHT19} now follows  from
\eqref{matrixTHHT}. The equality \eqref{matrixTHHT29} can be verified in the same manner.
\epr

\medskip
In what follows $M$ is the operator given by
\begin{equation}\label{invformM1}
M = \begin{bmatrix} M_{11} & M_{12} \\ M_{21} &  M_{22}  \end{bmatrix} :
\begin{bmatrix} L^1(\BR_+)^p \\ L^1(\BR_-)^q \end{bmatrix}  \to
\begin{bmatrix} L^1(\BR_+)^p \\ L^1(\BR_-)^q \end{bmatrix},
\end{equation}
where  $ M_{11} $, $ M_{12} $, $ M_{21} $, and $ M_{22} $ are the operators defined by  \eqref{M11a}--\eqref{M22a}.

\medskip
\begin{lem}\label{lem:invRjj}
Let $ \bR_{11} $, $ \bR_{12} $, $ \bR_{21} $, and $ \bR_{22} $ be defined by
\eqref{defR11}--\eqref{defR22},  and let $ M_{11} $, $ M_{12} $, $ M_{21} $, and $ M_{22} $ be defined by \eqref{M11a}--\eqref{M22a}.
Then
\begin{equation}\label{relRM}
\begin{bmatrix} J_{\sX_+} & 0 \\ 0 & J_{\sY_-} \end{bmatrix}
\begin{bmatrix} \bR_{11} & \bR_{12} \\ \bR_{21} & \bR_{22}  \end{bmatrix}
=
\begin{bmatrix} I_p & 0 & 0 & 0 \\ a & M_{11} & b  & M_{12} \\ 0 & 0 & I & 0 \\
c & M_{21} & d & M_{22} \end{bmatrix}
\begin{bmatrix} J_{\sX_+} & 0 \\ 0 & J_{\sY_-} \end{bmatrix} .
\end{equation}
In particular, $ \bR $ given by \eqref{invform1} is invertible if and only if $ M $ defined by
\eqref{invformM1} is invertible. Moreover,  $ M_{11} $ is invertible if and only if $ \bR_{11} $ is invertible and $ M_{22} $ is invertible if and only if $ \bR_{22} $ is invertible.
\end{lem}

\bpr
From the relations between Hankel-like operators and Toeplitz-like operators
on the one hand and Hankel integral operators and Wiener Hopf operators on the other hand
we have the identities:
\begin{align}
& J_{\sX_+} \bR_{11} = \begin{bmatrix} I_p & 0 \\ a & M_{11} \end{bmatrix} J_{\sX_+} ,
 \quad
& J_{\sX_+} \bR_{12} = \begin{bmatrix} 0 & 0 \\ b & M_{12} \end{bmatrix} J_{\sY_-},
\label{RM1112} \\
& J_{\sY_-} \bR_{21} = \begin{bmatrix} 0 & 0 \\ c & M_{21} \end{bmatrix} J_{\sX_+} ,
\quad
& J_{\sY_-} \bR_{22}  = \begin{bmatrix} I_q & 0 \\ d & M_{22} \end{bmatrix} J_{\sY_-} .
 \label{RM2122}
\end{align}
Putting together these equalities gives the equality \eqref{relRM}.
The equality \eqref{relRM} implies that $ \bR $ is invertible it and only if $ M $ is invertible.
The final statement follows from the first equality in \eqref{RM1112} and the second equality
in \eqref{RM2122}.
\epr

\medskip
We continue with specifying two other lemmas from Section \ref{sec:prelim}.

\begin{lem}\label{lem:altMij}
Since  conditions \textup{(C4)--(C6)} are satisfied, we have
\begin{align}
& M_{11} = I_p - H_{+,\a} H_{-,\a^\ast} + H_{+,\b} H_{-,\b^\ast} : L^1(\BR_+)^p \to L^1(\BR_+)^p,
 \label{formulaM11} \\
& M_{21} = T_{-,\d} H_{-,\b^\ast } - T_{-,\c} H_{-,\a^\ast} : L^1(\BR_+)^p \to L^1(\BR_-)^q ,
 \label{formulaM21}\\
& M_{12} =  T_{+,\a} H_{+,\c^\ast} - T_{+,\b}H_{+,\d^\ast}  : L^1(\BR_-)^q \to L^1(\BR_+)^p ,
 \label{formulaM12} \\
& M_{22} = I_q - H_{-,\d} H_{+,\d^\ast} + H_{-,\c} H_{+,\c^\ast} : L^1(\BR_-)^q \to L^1(\BR_-)^q.
\label{formulaM22}
\end{align}
\end{lem}

\bpr
The result is an immediate consequence of Lemma \ref{lemaltRij} and the relations
between the $ \bR_{ij} $ and $ M_{ij} $ in \eqref{RM1112} and \eqref{RM2122}.
\epr

\begin{lem}\label{lem:inverseM}
Let $ M_{ij} $, $ i,j=1,2$, be given by \eqref{M11a}--\eqref{M22a}, and let   $ M $ be given by \eqref{invformM1}. Since conditions \textup{(C1)--(C6)}  are satisfied, we have
\begin{equation}\label{equalityM1}
\begin{bmatrix} M_{11} & M_{12} \\ M_{21} &  M_{22}  \end{bmatrix}
\begin{bmatrix} I_p  & 0 \\ 0 &  -I_q  \end{bmatrix}
\begin{bmatrix} M_{11} & M_{12} \\ M_{21} &  M_{22}  \end{bmatrix}  =
\begin{bmatrix} M_{11} & 0 \\ 0 &  -M_{22}  \end{bmatrix}.
\end{equation}
In particular, $ M $ is invertible if and only if $ M_{11} $ and $ M_{22} $ are invertible.
Furthermore, in that case
\begin{equation}\label{invformM2}
M^{-1}  =
\begin{bmatrix} I_p & -M_{12} M_{22}^{-1} \\ -M_{21} M_{11}^{-1} & I_q \end{bmatrix} =
\begin{bmatrix} I_p & -M_{11}^{-1} M_{12}  \\ -M_{22}^{-1} M_{21}  & I_q \end{bmatrix}.
\end{equation}
\end{lem}

\bpr
The result is an immediate consequence of Lemma \ref{lem:inverseR} and the relations
 between  $ \bR_{ij} $ and $ M_{ij} $ given in \eqref{RM1112} and \eqref{RM2122}.\epr

\medskip
\noindent\textsc{Step 4. Results of Section \ref{sec:inversionthm}.}
We use Theorem \ref{thm:inversion1} to prove the ``if'' part of Theorem~\ref{thm:mainL1} and the identities \eqref{formula-Winv}, \eqref{formulaHg},  \eqref{formulaHgast}, and \eqref{formgg}.

First we check that  the various conditions appearing in Theorem \ref{thm:inversion1} are satisfied given our data. Since  $g\in L^1(\BR_+)^{p\ts q}=\sB_+$  is a solution to the  twofold EG inverse problem, Proposition \ref{cor:g-solinv} tells us that
\begin{equation*}
\boldsymbol\Omega\begin{bmatrix} \a \\ \c \end{bmatrix} =
\begin{bmatrix} e_\sA \\ 0 \end{bmatrix} \ands
\boldsymbol\Omega\begin{bmatrix} \b \\ \d \end{bmatrix} = \begin{bmatrix} 0 \\ e_\sD \end{bmatrix}.
\end{equation*}
Thus the identities in   \eqref{Omega-basic-eqs} are satisfied. Next,  note that the final conclusion of \textsc{Step 1} tells us that of items (a) and (b) in Theorem~\ref{thm:inversion1}  are also satisfied.

Thus Theorem~\ref{thm:inversion1}  tells us that the operator $\bOm$ is invertible.  But then
we can use Lemma \ref{WandOmega} to conclude  that the operator $W$ defined by \eqref{defZ1} is invertible too. This concludes the proof of the  ``if'' part of Theorem~\ref{thm:mainL1}.

Theorem \ref{thm:inversion1} also tells us that the   inverse $ \bR $ of  $ \bOm$ is given by \eqref{IGGIinverse}. From \eqref{OandW} and  \eqref{relRM} it then follows that the inverse of $W$ is  the operator $M$ defined by  \eqref{invformM1}. This proves identity
\eqref{formula-Winv}.

From Lemma \ref{lem:inverseM} we know that that $M_{11}$ and $M_{22}$ are invertible. The identities in \eqref{formulaHg} and \eqref{formulaHgast} are obtained by comparing the off diagonal entries of $W=M^{-1}$ in \eqref{formula-Winv}  and \eqref{invformM2}.

Finally to see  that the identities in \eqref{formgg} hold true, note that from  $ \bR \bOm = I $ it follows that
\[
\begin{bmatrix} I_p & 0 & 0 & 0 \\ a & M_{11} & b  & M_{12} \\ 0 & 0 & I & 0 \\
c & M_{21} & d & M_{22} \end{bmatrix}
\begin{bmatrix} I_p & 0 & 0 & 0 \\ 0 & I & g & H_{+,g} \\ 0 & 0 & I & 0 \\
g^\ast & H_{-,g^\ast} & 0 & I \end{bmatrix}  = I.
\]
In  particular, $ M_{11} g + b = 0$ and $M_{22} g^\ast + c = 0$. Using the invertibility of $M_{11}$ and $M_{22}$ we obtain the formulas for $g$ and $g^\ast$ in \eqref{formgg}. This completes the proof.
\epr


\subsection{Proof of Theorem \ref{thm:mainL2}}
Throughout  this subsection, as in Theorem \ref{thm:mainL2},   $\{a, b,c,d \}$ are the functions given by  in \eqref{eq:data1} and \eqref{eq:data2}, and  $\a$, $\b$, $\c$, $\d$ are  the functions given by \eqref{def:abcd-greek}. Furthermore, $e_p$ and $e_q$ are the functions on $\BR$ identically equal to the unit matrix $ I_p$  and $I_q$, respectively.  Finally $\sM=\sM_{\sA,\sB, \sC, \sD}$ is the admissible algebra constructed in the beginning three   paragraphs of the present section. In what follows we split the proof  of Theorem \ref{thm:mainL2} into two parts.

\smallskip
\noindent\textsc{Part 1.} In this part we assume that  the twofold  EG inverse problem associated with the data  set  $\{a, b,c,d \}$ has a solution, $g \in L^1(\BR_+)^{p \times q} =\sB_+$ say. Then we know  from Proposition \ref{prop:cond123} that conditions (C1) -- (C3)  are satisfied. But the latter, using  the final part of Remark \ref{rem:coincide}, implies that condition (L1) is satisfied. Furthermore, the second part of  Theorem   \ref{thm:mainL1} tells us that the operators $M_{11}$ and $M_{22}$ are invertible, and hence   condition (L2) is satisfied too. Finally, the two identities in \eqref{formgg} yield the two identities in \eqref{eq:1defgg*}. This concludes the
first part of the proof.

\smallskip
\noindent\textsc{Part 2.} In this part we assume that (L1) and (L2) are satisfied. Our aim is to show that the twofold  EG inverse problem associated with the data  set  $\{a, b,c,d \}$ has a solution.

We begin with some preliminaries. Recall that
\[
a_0=P_{\sA_d}\a=e_p \ands d_0=P_{\sD_d}\d=e_q.
\]
Furthermore, from the identities in \eqref{eq:L1abc} we know that (C1)--(C3) are satisfied. But then we can   repeat the arguments in Step 1 of the proof of Theorem ~\ref{thm:mainL1} to show that conditions (C4)--(C6)  are also satisfied. Thus all conditions (C1)--(C6)  are fulfilled. Finally, note   that in the paragraph  directly after Theorem~\ref{thm:mainL2} we showed  that (L1) and (L2) imply that $M_{11}$ and $M_{22}$ are invertible, and hence we can apply Lemma \ref{lem:inverseM} to see that  the inverse of
\[
M:=\begin{bmatrix}M_{11}& M_{12}\\ M_{21}& M_{22}\end{bmatrix}
\]
 is given by \eqref{invformM2}, i.e.,
 \begin{equation}\label{invformM2a}
M^{-1}  =
\begin{bmatrix} I_p & -M_{12} M_{22}^{-1} \\ -M_{21} M_{11}^{-1} & I_q \end{bmatrix} =
\begin{bmatrix} I_p & -M_{11}^{-1} M_{12}  \\ -M_{22}^{-1} M_{21}  & I_q \end{bmatrix}.
\end{equation}
It remains to show that there exists a $ g \in L^1(\BR_+)^{p \times q} $ such that
\begin{equation}\label{eq:Hankelops1}
 - M_{11}^{-1} M_{12} = H_{+,g} \ands - M_{22}^{-1} M_{21} = H_{-,g^\ast} .
\end{equation}
To do this we need (in the context of the present setting)  a more general version of Lemma \ref{isHankel}.  We cannot apply Lemma \ref{isHankel} directly because of the role of the constant functions in $ \BC^{ p \times p} e_p $ and $ \BC^{q \times q} e_q $.  The more general version of Lemma \ref{isHankel} will be given and proved  in the following intermezzo.

\smallskip
\noindent\textbf{Intermezzo.}
First we introduce the  required transition operators.  Let $ \tau \geq 0 $. Define $ V_{r,\tau} : L^2(\BR_+)^r \to L^2(\BR_+)^r $ by
\[
 \left( V_{r,\tau} f \right)(t) = \left\{ \begin{matrix} f(t-\tau), & t \geq \tau, \\
 \noalign{\smallskip} 0, & 0 \leq t \leq \tau .  \end{matrix} \right.
\]
Note that  its adjoint $ V_{r,\tau}^\ast $ is given  by $ \left( V_{r,\tau}^\ast f \right) (t) = f(t+\tau) $ for $ t \geq 0 $.
We also need the the flip over operator $J_r$  from $  L^2(\BR_+)^r $ to $  L^2(\BR_-)^r $ given by $ (J_rf)(t) = f(-t) $.
With some abuse of notation we also consider $ V_{r,\tau} $, $ J_r $ and their adjoints as  operators acting on $ L^1$-spaces.
For $ \varphi \in \BC^{k \times m} \dot{+} L^1(\BR_+)^{k \times m} $ we  then have
\begin{align*}
V_{k,\tau}^\ast H_{+,\varphi} = H_{+,\varphi} J_m V_{m,\tau} J_m,
& \ands
J_k V_{k,\tau}^\ast J_k H_{-,\varphi} = H_{-,\varphi} V_{m,\tau},  \\
V_{k,\tau}^\ast T_{+,\varphi} V_{m,\tau} = T_{+,\varphi}.
& \ands
J_k V_{k,\tau}^\ast J_k T_{-,\varphi} J_m V_{m,\tau} J_m = T_{-,\varphi}.
\end{align*}

\medskip
The following lemma is the more general version of Lemma \ref{isHankel} mentioned above.
The result  will  be used to show that $ M_{11}^{-1} M_{12} $ and $ M_{22}^{-1} M_{21} $ are classical Hankel integral operators.

\begin{lem}\label{lem:isHankel}
With $ M_{ij}$,  $ 1 \leq i,j \leq 2$, defined by \eqref{M11a}--\eqref{M22a}  we have  the following equalities:
\[
M_{11} V_{p,\tau}^\ast M_{12} = M_{12} J_q V_{q,\tau} J_q M_{22}
\ands  M_{22} j_q V_{q,\tau}^\ast J_q M_{21} = M_{21} V_{p,\tau} M_{11} .
\]
Moreover, if $ M_{11} $ and $ M_{22} $ are invertible, then
\[
V_{p,\tau}^\ast M_{12} M_{22}^{-1} = M_{11}^{-1} M_{12} J_q V_{q,\tau} J_q , \quad
J_q V_{q-,\tau}^\ast J_q M_{21} M_{11}^{-1}= M_{22}^{-1} M_{21} V_{p,\tau}.
\]
\end{lem}

\bpr
First we will prove that $ M_{11} V_{p,\tau}^\ast M_{12} = M_{12} J_q V_{q,\tau} J_q M_{22} $.
We start with deriving the equality
\begin{equation}\label{matrixTHSHT1}
\begin{bmatrix}  T_{+,\a^\ast} \\ T_{+,\b^\ast}\end{bmatrix} V_{p,\tau}^\ast
\begin{bmatrix} H_{+,\b} & H_{+,\a} \end{bmatrix} =
\begin{bmatrix}  H_{+,\c^\ast} \\ H_{+,\d\ast} \end{bmatrix} J_q V_{q,\tau} J_q
\begin{bmatrix} T_{-,\d} & T_{-,\c} \end{bmatrix} .
\end{equation}
Let $ f_- \in L^1(\BR_-)^q $ and $ h_- \in L^1(\BR_-)^p  $.
To obtain \eqref{matrixTHSHT1}, first notice that
\begin{align*}
 V_{p,\tau}^\ast \begin{bmatrix} H_{+,\b} & H_{+,\a} \end{bmatrix} \begin{bmatrix} f_- \\ h_- \end{bmatrix}
 & = \begin{bmatrix} H_{+,\b}  & H_{+,\a}  \end{bmatrix}
 \begin{bmatrix} J_q V_{q,\tau} J_q & 0  \\ 0 & J_p V_{p,\tau} J_p \end{bmatrix}
\begin{bmatrix} f_- \\ h_- \end{bmatrix} ; \\
J_q V_{q,\tau} J_q \begin{bmatrix} T_{-,\d} & T_{-,\c} \end{bmatrix} \begin{bmatrix} f_- \\ h_- \end{bmatrix}
& = \begin{bmatrix} T_{-,\d} & T_{-,\c} \end{bmatrix}
\begin{bmatrix} J_q V_{q,\tau} J_q & 0  \\ 0 & J_p V_{p,\tau} J_p \end{bmatrix}
\begin{bmatrix} f_- \\ h_- \end{bmatrix}.
\end{align*}
Then use \eqref{matrixTHHT19} to get that
\begin{align*}
& \begin{bmatrix}  T_{+,\a^\ast} \\ T_{+,\b^\ast} \end{bmatrix} V_{p,\tau}^\ast
\begin{bmatrix} H_{+,\b} & H_{+,\a} \end{bmatrix} \begin{bmatrix} f_- \\ h_- \end{bmatrix} \\
& =
{\begin{bmatrix}  T_{+,\a^\ast} \\ T_{+,\b^\ast}\end{bmatrix}
\begin{bmatrix} H_{+,\b} & H_{+,\a} \end{bmatrix}
\begin{bmatrix} J_q V_{q,\tau} J_q & 0  \\ 0 & J_p V_{p,\tau} J_p \end{bmatrix}}
\begin{bmatrix} f_- \\ h_- \end{bmatrix}\\
& =
\begin{bmatrix}  H_{+,\c^\ast} \\ H_{+,\d\ast}\end{bmatrix}
\begin{bmatrix} T_{-,\d} & T_{-,\c} \end{bmatrix}
\begin{bmatrix} J_q V_{q,\tau} J_q & 0  \\ 0 & J_p V_{p,\tau} J_p \end{bmatrix}
\begin{bmatrix} f_- \\ h_- \end{bmatrix}  \\
& =
\begin{bmatrix}  H_{+,\c^\ast} \\ H_{+,\d\ast}\end{bmatrix}  J_q V_{q,\tau} J_q
\begin{bmatrix} T_{-,\d} & T_{-,\c} \end{bmatrix} \begin{bmatrix} f_- \\ h_- \end{bmatrix}.
\end{align*}
We proved \eqref{matrixTHSHT1}.
By multiplying \eqref{matrixTHSHT1} on the left  and the right by
\[
\begin{bmatrix}  T_{+,\a} & T_{+,\b} \end{bmatrix}  \quad \mbox{and \ }
\begin{bmatrix} T_{-,\d^\ast } \\ T_{-,\c^\ast } \end{bmatrix},
\]
respectively, one gets $ M_{11} V_{p,\tau}^\ast M_{12} = M_{12} J_q V_{q,\tau} J_q M_{22} $.

The equality $ M_{22} J_q V_{q,\tau}^\ast J_q M_{21} = M_{21} V_{p,\tau} M_{11} $ can be proved in
a similar way.
The claim regarding the case that $M_{11}$ and $M_{22}$ are invertible follows trivially.
\epr

\medskip
\noindent
\textbf{We continue with the second part of the proof.} It remains to show   that there exists a $ g \in L^1(\BR_+)^{p \times q} $ such that the two identities in \eqref{eq:Hankelops1} are  satisfied. For this purpose we need Lemma \ref{lem:isHankel} and various  results presented in the Appendix.  In particular, in what follows we need the Sobolev space  $\sob(\BR_+)^n$  which consist  of  all functions $\va \in L^1(\BR_+)^n$ such that $\va$  is absolutely continuous on compact intervals of $\BR_+$  and    $\va^\prime\in L^1(\BR_+)^n$ (See  Subsection  \ref{ss:Hankel-L1}). Notice that  $ M_{12} J_q $ is a sum  of products of Wiener-Hopf operators and classical
Hankel integral operators. Therefore, by Lemma \ref{lem:Wiener1},
the operator  $ M_{12} J_q $ maps $ \sob(\BR_+)^q $ into $ \sob(\BR_+)^p $
and $ M_{12} J_q |_{\sob(\BR_+)^q}  $ is bounded as  an operator
from $ \sob(\BR_+)^q$ to $ \sob(\BR_+)^p $.
The operator $ M_{11} $ is of the form  \eqref{eq:defM12}.
Thus Lemma~\ref{lem:M-inv}  tells us that $ M_{11}^{-1} $ satisfies the condition (H1) in Theorem \ref{thm:main1}. We conclude that the operator $ - M_{11}^{-1} M_{12} J_q $ maps
$ \sob(\BR_+)^q $ into $ \sob(\BR_+)^p $ and $ - M_{11}^{-1} M_{12} J_q |_{\sob(\BR_+)^q} $  is bounded as an operator from $ \sob(\BR_+)^q $ to  $ \sob(\BR_+)^p $.
Also we know from Lemma~\ref{lem:isHankel} that
\[
 V_{p,\tau} (- M_{11}^{-1} M_{12} )J_q =  (- M_{11}^{-1} M_{12} )J_q V_{q,\tau}, \quad \forall  \tau \geq 0.
 \]
 According to Theorem \ref{thm:char-Hankel1} it follows that there exists a $ k \in L^\infty(\BR_+)^{p \times q} $ such that
$ - M_{11}^{-1} M_{12} J_q = H( k )$.
But then we can  apply Corollary \ref{cor:Hankel-int1}
to show that there exists a $ g \in L^1(\BR_+)^{p \times q} $ such that
 $ - M_{11}^{-1} M_{12} = H_{+,g} $.
In a similar way we prove that there exists a $ h \in L^1(\BR_-)^{q \times p} $
such that $ - M_{22}^{-1} M_{21} = H_{-,h} $.

Notice that the operators $ M_{ij} $ can be considered to be
operators acting between $ L^2 $-spaces.
This can be done because the Hankel and Wiener-Hopf operators that constitute the $ M_{ij} $
can be seen as operators between $ L^2 $ spaces.
Recall that $ H_{+,\rho}^\ast = H_{-,\rho^\ast} $ and $ T_{\pm,\rho}^\ast = T_{\pm,\rho^\ast} $.
Using Lemma~\ref{lem:altMij} and the definition of $ M_{12} $, we see that $ M_{12}^\ast = M_{21} $, and using the definitions of $ M_{11} $
and $ M_{22} $ we may conclude that $ M_{11}^\ast = M_{11} $ and $ M_{22}^\ast = M_{22} $.
From the equality \eqref{invformM2} one sees that $ H_{-,h} = H_{+,g}^\ast $.
Hence $ h = g^\ast $.

We need to show that for this $ g $ the inclusions \eqref{eq:probl1}  and \eqref{eq:probl2} are satisfied, or equivalently that \eqref{eq:probl-v2}  is satisfied.
Define for this $ g $ the operator $  W $ by \eqref{defZ1}.
We already know that $ W $ is invertible and that its  inverse is $ M $.
Let $ \tilde{a} $, $ \tilde{b} $, $  \tilde{c} $ and $  \tilde{d} $ be the solution of
\[
W \begin{bmatrix} \tilde{a} & \tilde{b} \\ \tilde{c} & \tilde{d} \end{bmatrix} = \begin{bmatrix} 0 & -g \\ -g^\ast & 0 \end{bmatrix}.
\]
Then put $ \tilde{\a} = e_p + \tilde{a} $, $ \tilde{\b} =\tilde{b} $,
$\tilde{\c} = \tilde{c} $ and $ \tilde{\d} = e_q + \tilde{d} $.
With the data $ \tilde{\a} $, $ \tilde{\b} $, $ \tilde{\c} $, and $ \tilde{\d} $ produce a new
$ \tilde{M}$ which, according to Theorem \ref{thm:mainL1}, is also the inverse of $W$. But then the the old $ M $ and new $ \tilde{M} $ are the same, and
hence $ g $ is the solution of the EG inverse problem associated with $\a$, $\b$, $\c$ and $\d$.
All the conditions of Theorem~\ref{thm:mainL1} are now satisfied and we conclude that
$ g = - M_{11}^{-1} b  $ and $ g^\ast = - M_{22}^{-1} c $.
\epr

\section{The EG inverse problem with additional invertibility conditions}\label{sec:can-case}

As before $ \a \in \sA_+ $, $ \b \in \sB_+ $, $ \c \in \sC_- $, and $ \d \in \sD_- $.
In this section we consider the case when $ \a $ is invertible in $ \sA_+ $ and
$ \d $ is invertible in $ \sD_- $.
Notice that in the example discussed in Subsection \ref{subsec:Mat} this condition
is satisfied whenever $a_0:=P_{\sA_d}\a$ and $d_0:=P_{\sD_d}\d$ are invertible.

\begin{thm}\label{thm:cancase}
Let $\a \in \sA_+ $, $\b \in \sB_+ $, $\c \in \sC_- $ and $\d  \in \sD_- $ and assume that
$ \a $ and $ \d $ are invertible in $ \sA_+ $ and $ \sD_-$, respectively.
If, in addition,  $\a$, $\b$, $\c$ and $\d$ satisfy the conditions \textup{(C1)}
and \textup{(C2)}, then $ g_1 =  - P_{\sB_+} \left( \a^{-\ast} \c^\ast \right) $
is the unique element of $\sB_+ $ that satisfies \eqref{inclu12} and
$ g_2 =  - P_{\sB_+} (\b \d^{-1} ) $ is the unique element of $ \sB_+ $ that satisfies
\eqref{inclu34}.
Moreover, in that case,  $g_1 = g_2 $ if and only if condition \textup{(C3)} is satisfied. In particular, if \textup{(C1)}-\textup{(C3)} hold, then $g=g_1=g_2$ is the unique solution to the twofold  EG inverse problem associated with the data set $\{\a,\b,\c,\d\}$.
\end{thm}
\bpr
The inclusion  $\a^{-1} \in \sA_+ $ implies that $ a_0 $ is invertible with inverse in
$ \sA_d $. Similarly, $ \d^{-1} \in \sD_- $ implies that $ d_0 $ is invertible with inverse in $ \sD_d $.

Let $g_1: =  - P_{\sB_+} \left( \a^{-\ast} \c^\ast \right)$. First we prove that $g_1 $ satisfies the second inclusion in \eqref{inclu12}. From the definition of  $g_1$ it follows that  $g_1+ \a^{-\ast} \c^\ast=\b_1$ for some $\b_1\in \sB_-$. Taking  adjoints we see that  $g_1^\ast + \c \a^{-1} =\b_1^\ast \in \sC_+$. Multiplying for the right by $\a$  and using the multiplication table in Section \ref{sec:setting} we see that
\[
g_1^\ast \a + \c  \in \sC_+ \sA_+ \subset \sC_+.
\]
So $g_1 $ is a solution of the second inclusion in \eqref{inclu12}. Notice that in this paragraph we did not yet use that $\a^{-1} \in \sA_+ $.

The next step is to show that $g_1$ is the unique element
of $\sB_+ $ that satisfies the second inclusion in \eqref{inclu12}.
Assume that $\varphi_1 \in \sB_+ $ and $ P_{\sC_-} ( \varphi_1^\ast \a +\c  ) =0 $.
We will prove that $ \varphi_1 = g_1 $.
Notice that
\[
\a^\ast (\varphi_1 - g_1 ) = ( \a^\ast \varphi_1 + \c^\ast ) - ( \a^\ast g_1 + \c^\ast ) \in \sB_-.
\]
Hence $ \varphi_1 - g_1 \in \sB_+ $ and $\a^\ast ( \varphi_1 - g_1 ) \in \sB_- $.
Since $ \a^{-\ast} \in \sA_- $ we have that
\[
 \varphi_1 - g_1 = \a^{-\ast} \a^\ast ( \varphi_1 -  g_1 ) \in \sA_- \sB_- \subset \sB_-
\]
and hence $ \varphi_1 - g_1 =0 $.

Next remark that $ P_{\sA_-}  ( \a^{-1} - a_0^{-1} ) =0 $. Indeed
\[
\a^{-1} - a_0^{-1} =  \a^{-1}( a_0 - \a ) a_0^{-1}  \in \sA_+ \sA_+^0 \sA_d \subset \sA_+^0.
\]
To show that $ g_1 $ satisfies the first inclusion in \eqref{inclu12}, note that
\[
\a + g_1 \c - e_\sA =
\a -  ( \a^{-\ast} \c^\ast ) \c + \left( P_{\sB_-} (\a^{-\ast} \c^\ast ) \right) \c -  e_\sA.
\]
Now use that $ \c^\ast \c = \a^\ast \a -  a_0 $ to see that
\begin{align*}
\a + g_1 \c - e_\sA & =
\a - \a^{-\ast} ( \a^\ast \a - a_0 ) + \left( P_{\sB_-} (\a^{-\ast} \c^\ast ) \right) \c -  e_\sA
\\
& =  \a^{-\ast} a_0 - e_\sA   + \left( P_{\sB_-} (\a^{-\ast} \c^\ast ) \right) \c \\
&= ( \a^{-\ast} -  a_0^{-1} ) a_0  + \left( P_{\sB_-} (\a^{-\ast} \c^\ast ) \right) \c.
\end{align*}
Since $ \c \in \sC_- $, we have that
$  P_{\sA_+} \left[ \left( P_{\sB_-} (\a^{-\ast} \c^\ast ) \right) \c \right] =0$,
and since $ P_{\sA_-}  ( \a^{-1} - a_0^{-1} ) =0 $, we also have that
$ P_{\sA_+} ( \a^{-\ast} - a_0^{-1} ) =0$.
We proved the first inclusion in \eqref{inclu12}.

Next, let $ g_2: =  - P_{\sB_+} (\b \d^{-1} ) $.  We will show that $g_2 $ is the unique element of $ \sB_+$ that satisfies the first inclusion in \eqref{inclu34}. To do this, note that $ g_2+\b \d^{-1}=\b_2$ for some $\b_2\in \sB_-$, which implies that
\[
g_2 \d + \b\in \sB_- \sD_- \subset \sB_-.
\]
We proved that $ P_{\sB_+} (g_2 \d + \b ) =0 $.
Assume that $\varphi_2 \in \sB_+ $ satisfies also the first inclusion in \eqref{inclu34}.
Then
\[
(\varphi_2 -g_2) \d = (\varphi_2 \d + \b ) - ( g_2 \d + \b ) \in \sB_-.
\]
Since $\d \in \sD_- $ we have that
$ (\varphi_2 -g_2) = (\varphi_2 -g_2) \d \, \d^{-1}  \in \sB_- \sD_- \subset \sB_- $.
Hence $ \varphi_2 - g_2 = 0 $ and $g_2 $ is the unique solution of the first inclusion in \eqref{inclu34}.


We proceed with showing  that $ g_2 $ also satisfies the second inclusion in \eqref{inclu34}. Indeed
\begin{align*}
\b^\ast g_2 + \d^\ast - e_\sD & =  \b^\ast ( - \b \d^{-1} + P_{\sB_-} ( \b \d^{-1}) )  + \d^\ast - e_\sD  \\
&= - \b^\ast \b \d^{-1} + \d^\ast - e_\sD  + \b^\ast P_{\sB_-}(\b \d^{-1}) \\
&= - ( \d^\ast \d - d_0 ) \d^{-1} + \d^\ast - e_\sD  + \b^\ast P_{\sB_-}( \b \d^{-1}) \\
&= d_0 \d^{-1} - e_\sD  + \b^\ast P_{\sB_-}( \b \d^{-1}) \in \sD_-^0.
\end{align*}
Here we used that
$ P_{\sD_+ } ( d_0 \d^{-1} - e_\sD ) = P_{\sD_+} \bigl( d_0^{-1} ( e_\sD -\d ) \d^{-1} \bigr) =0.$
We proved that $ g_2 $ satisfies the second inclusion in \eqref{inclu34}.

If $ \a^\ast \b - \c^\ast \d =0 $, then $ \b = \a^{-\ast} \c^\ast \d $.
It follows that
\begin{align*}
0 & = P_{\sB_+ } ( \b - \a^{-\ast} \c^\ast \d ) =
P_{\sB_+ } \Bigl( \b - P_{\sB_+ } \bigl( (\a^{-\ast} \c^\ast)  \d \bigr) \Bigr)
\\ &
= P_{\sB_+ } \bigl( \b - P_{\sB_+ } ( \a^{-\ast} \c^\ast ) \d  \bigr)
= P_{\sB_+ } ( \b + g_1  \d  ).
\end{align*}
So $ g_1 $ also solves the first inclusion in \eqref{inclu34} and the uniqueness of the solution
gives $ g_1 = g_2$.

Conversely, if $g_1 = g_2 $ then we have a solution of the inclusions
 \eqref{inclu12} and \eqref{inclu34}.
It follows from \cite[Theorem 1.2]{KvSch13} that the conditions (C1)--(C3) are satisfied and
in particular we get $ \a^\ast \b - \c^\ast \d =0 $.
\epr

\medskip
In the next proposition we combine the results of Theorem \ref{thm:cancase} with those of
Theorem~\ref{thm:inversion1}.

\begin{prop}\label{prop:inversion2}
Let $ \sM_{\sA,\sB,\sC,\sD} $ be an admissible algebra,  and
let $\a \in \sA_+ $, $\b \in \sB_+ $, $\c \in \sC_-$, $\d \in \sD_- $ be such that:
\begin{itemize}
\item[\textup{(a)}] $ \a $ is invertible in $ \sA_+ $ and $ \d $ is invertible in $ \sD_- $;
\item[\textup{(b)}] conditions \textup{(C1)--(C3)} are satisfied.
\end{itemize}
Let
\begin{equation}\label{def:phiOm}
\varphi = - P_{\sB_+}( a^{-\ast} \c^\ast) \ands
\boldsymbol{\Omega} =
\begin{bmatrix} I_{\sX_+} & \bH_{+,\varphi} \\ \bH_{-,\varphi^\ast} & I_{\sY_-} \end{bmatrix}.
\end{equation}
Then $ \varphi $ is the solution of the twofold EG inverse problem associated with
$ \a, \b , \c , \d $, the operator $ \boldsymbol{\Omega} $ is invertible, and
\begin{equation}\label{IHHIinverse}
\boldsymbol{\Omega}^{-1} =\begin{bmatrix} \bR_{11} & \bR_{12} \\ \bR_{21} & \bR_{22} \end{bmatrix},
\end{equation}
where the operators $\bR_{ij}$, $1 \leq i,j \leq 2$, are defined by \eqref{defR11}--\eqref{defR22}.
In particular, the operators  $ \bR_{11} $ and $  \bR_{22} $  are invertible and
\begin{equation} \label{formulaHH*}
\bH_{+,\varphi } = - \bR_{11}^{-1} \bR_{12} =
- \bR_{12} \bR_{22}^{-1}, \quad \bH_{-, \varphi^\ast } = - \bR_{21} \bR_{11}^{-1} =
- \bR_{22}^{-1} \bR_{21} ,
\end{equation}
and
\begin{equation}\label{eq:phibc}
\bR_{11} \varphi = - \b  \ands  \bR_{22} \varphi^\ast = -\c .
\end{equation}
\end{prop}

\bpr
Since the conditions of Theorem \ref{thm:cancase} are satisfied,   $ \varphi $ is
the unique solution of the twofold EG inverse problem associated with $  \a, \b, \c, \d $.
Also the fact that $ \a $ is invertible in $ \sA_+ $ and $ \d  $ is invertible in $ \sD_- $ gives
that $ a_0 $ and $ d_0 $ are invertible in $ \sA_d $ and $ \sD_d $, respectively, and
that according to Lemma~\ref{lem:C1-C6} also the conditions (C4)--(C6) are satisfied.
But then all the conditions of Theorem \ref{thm:inversion1} are satisfied.
The equalities \eqref{IHHIinverse}, \eqref{formulaHH*} and \eqref{eq:phibc} are now immediate from
Theorem~\ref{thm:inversion1}.
\epr

\medskip
Specifying Theorem \ref{thm:cancase} for  the example discussed in
Subsection \ref{subsec:Mat} yields the following corollary.

\begin{cor} Let $ \sA $, $\sB $, $\sC$, $\sD$, $ \sA_\pm $, $\sB_\pm $, $\sC_\pm $, $\sD_\pm $, and
$ \sA_d $, $\sB_d $, $\sC_d $, $\sD_d $ be as in Subsection \ref{subsec:Mat}, and let $\a \in \sA_+ $, $\b \in \sB_+ $, $\c \in \sC_- $, $\d  \in \sD_- $ be given. Assume that
\begin{itemize}
\item[\textup{(a)}]  $ \a_0 = P_{\sA_d} \a $ and $ \d_0 = P_{\sA_d} \d  $ are invertible.
\item[\textup{(b)}]   $\a$, $\b$, $\c$ and $\d$ satisfy  conditions \textup{(C1)}--\textup{(C3)}
\end{itemize}
Then $ g =  - P_{\sB_+} \left( \a^{-\ast} \c^\ast \right) $ is the unique element of $\sB_+ $ that satisfies \eqref{inclu12} and \eqref{inclu34}.
\end{cor}
\bpr
We only need to recall that the invertibility of the diagonal matrices $ a_0 $ and $ d_0$
implies invertibility of $ \a $ and $ \d$ in $ \sA_+ $ and $ \sD_- $, respectively.
\epr

\medskip
From  the above corollary it also follows that in the numerical Example \ref{ex:simple1}
the solution  $ g $,
\[
g = - P_{\sB_+} ( \a^{-\ast} \c^\ast ) =
\begin{bmatrix} 1 & 2 & 0 \\ 0 & 1 & 2 \\ 0  & 0 & 1 \end{bmatrix} ,
\]
is the unique solution of equations \eqref{inclu12b} and \eqref{inclu34b}.

\setcounter{equation}{0}
\section{Wiener algebra on the circle}\label{sec:wienerC}

In this section (as announced in Subsubsection \ref{ssec:WT}) we show how the solution of  the discrete twofold EG inverse problem, Theorem 4.1 in \cite{tHKvS1}, can be obtained as a corollary of  our abstract Theorem \ref{thm:mainthm}.

Let us first recall the discrete twofold EG inverse problem   as it was presented in \cite{tHKvS1}.  This requires some preliminaries. Throughout  $\sW^{n\ts m} $ denotes the space of $ n\ts m$ matrix functions with entries in the Wiener algebra on the unit circle which is denoted by $\sW$ and not by $\sW(\BT)$ as in Subsubsection \ref{ssec:WT}.  Thus a matrix function  $\va $ belongs to $\sW^{n\ts m} $  if and only if $\va$ is continuous on the unit circle and its Fourier coefficients
$\dots \va_{-1}, \va_0, \va_1, \ldots$ are absolutely summable.
We set
\begin{align*}
\sW_+^{n\ts m}&=\{\va\in \sW^{n\ts m}\mid \va_j= 0, \quad
\mbox{for } j=-1, -2, \dots\},   \\
\sW_-^{n\ts m}&=\{\va\in \sW^{n\ts m}\mid \va_j= 0, \quad \mbox{for } j=1, 2, \dots\}, \\
\sW_d^{n\ts m} &=\{\va\in \sW^{n\ts m}\mid \va_j= 0, \quad \mbox{for } j\not =0\}, \\
\sW_{+,0}^{n\ts m}&=\{\va\in \sW^{n\ts m}\mid \va_j= 0, \quad \mbox{for } j=0,-1, -2, \dots\},   \\
\sW_{-,0}^{n\ts m}&=\{\va\in \sW^{n\ts m}\mid \va_j= 0, \quad \mbox{for } j=0,1, 2, \dots\}.
\end{align*}
Given $\va \in \sW^{n\ts m}$ the function $\va^*$ is defined by $\va^*(\z)=\va(\z)^*$ for each $\z\in \BT$.
Thus the $j$-th Fourier coefficient of  $\va^*$ is given by $(\va^*)_j=(\va_{-j})^*$.
The map $\va \mapsto  \va^*$ defines an involution  which transforms
$\sW^{n\ts m}$ into $\sW^{m \times n}$, $ \sW_+^{n\ts m}$ into $ \sW_-^{m \ts n}$,
$ \sW_{-, 0}^{n\ts m}$ into $\sW_{+,0}^{m\ts n}$, etc.

The data of the discrete EG inverse problem  consist  of four functions, namely
\begin{equation}\label{Fabcd2}
\a\in \sW_+^{p\ts p}, \quad \b\in \sW_+^{p\ts q}, \quad \c\in \sW_-^{q\ts p}, \quad \d\in \sW_-^{q\ts q},
\end{equation}
and the problem is to find   $g\in \sW_+^{p\ts q}$ such that
\begin{align}
& \a + g \c - e_p \in  \sW_{-,0}^{ p \times p }  \ands
g^\ast \a  + \c \in   \sW _{+ ,0} ^{ q \times p }; \label{incluD12} \\
&  g \d + \b \in  \sW_{-,0}^{p \times q} \ands
\d+ g^\ast \b - e_q \in   \sW _{+,0}^{q \times q}.  \label{incluD34}
\end{align}
Here $ e_p $ and $e_q$ denote the functions identically equal to the identity matrices
$ I_p$  and $I_q$, respectively.  If $g$ has these properties, we refer to $g$ as a
\emph{solution to the discrete twofold EG inverse problem associated with the data set}
$\{  \a, \b, \c,\d \}$.  If a solution exists, then we know from Theorem 1.2  in \cite{KvSch13} that necessarily the following identities hold:
\begin{equation}\label{condD1}
\a^\ast \a - \c^\ast \c =  a_0, \quad \d^\ast \d - \b^\ast \b= d_0,\quad \a^\ast \b = \c^\ast \d.
\end{equation}
Here $a_0$ and $d_0$ are the zero-th Fourier coefficient of $\a$ and $\d$, respectively, and  we
identify the  matrices with $a_0$ and $d_0$ with the matrix functions on $\BT$ that are identically  equal to  $a_0$ and $d_0$, respectively.   In this section we shall assume that $ a_0  $ and $ d_0$ are invertible.
Then \eqref{condD1} is equivalent to
\begin{equation}\label{condD2}
\a a_0^{-1}  \a^\ast - \c a_0^{-1} \c^\ast =  e_p , \quad
\d d_0^{-1}  \d^\ast - \b d_0^{-1} \b^\ast =  e_q , \quad
\a a_0^{-1} \c^\ast = \b d_0^{-1} \d^\ast
\end{equation}
Finally, we associate with the data  $\a, \b, \c, \d$ the following  operators:
\begin{align}
& R_{11} = T_{+,\a} a_0^{-1} T_{+,\a^\ast} - T_{+,\b} d_0^{-1} T_{+,\b^\ast} : \sW_+^p \to \sW_+^p , \label{defR11D} \\
& R_{21} = H_{-,\c} a_0^{-1} T_{+,\a^\ast} - H_{-,\d} d_0^{-1} T_{+,\b^\ast} : \sW_+^p \to \sW_-^q , \label{defR21D}\\
& R_{12} = H_{+,\b} d_0^{-1} T_{-,\d^\ast} - H_{+,\a} a_0^{-1} T_{-,\c^\ast} : \sW_-^q \to \sW_+^p , \label{defR12D}  \\
& R_{22} = T_{-,\d} d_0^{-1} T_{-,\d^\ast} - T_{-,\c} a_0^{-1} T_{-,\c^\ast} : \sW_-^q \to \sW_-^q . \label{defR22D}
\end{align}
Here $T_{+,\a},\  T_{+,\a^\ast}, \ T_{+,\b}, \ T_{+,\b},\ T_{-,\c}, \ T_{-,\c^\ast}, \ T_{-,\d},\ T_{-,\d^\ast}$ are Toeplitz operators  and $H_{+,\a}, \ H_{+,\b}, \ H_{-,\c}, \ H_{-,\d}$  are Hankel operators. The definitions of these operators can be found in the final paragraph of this section.

The next theorem gives the solution of the discrete twofold EG inverse problem. By applying Fourier transforms it is straightforward to check that the theorem is just equivalent to \cite[Theorem 4.1]{tHKvS1}.

\begin{thm}\label{thm:mainthmDS1}
Let $\a$, $\b$, $\c$, $\d$ be the functions  given by \eqref{Fabcd2}, with both matrices $a_0$ and $d_0$ invertible. Then the discrete   twofold EG inverse problem associated with the data  set
$ \{ \a, \b, \c, \d \} $
has a solution if and only the following conditions are satisfied:
\begin{itemize}
\item[\textup{(D1)}] the identities in \eqref{condD1} hold true;
\item[\textup{(D2)}]
the operators $ R_{11} $ and $ R_{22} $  defined by  \eqref{defR11D}
and \eqref{defR22D} are one-to-one.
\end{itemize}
Furthermore, in that case  $ R_{11} $ and $ R_{22} $ are invertible,
the solution is unique and  the unique solution $g$ and its adjoint are  given by
\begin{equation} \label{unique1}
g = - R_{11}^{-1} \b \ands    g^\ast  = - R_{22}^{-1} \c.
\end{equation}
\end{thm}

The next step is to show how the above theorem can be derived as a corollary of
our abstract Theorem \ref{thm:mainthm}. This requires to put  the inverse problem in the context of  the general scheme of Sections  \ref{sec:preliminaries}--\ref{sec:solinverseprobl}.  To do this (cf., the first paragraph of  Subsubsection \ref{ssec:WT})  we use  the following choice of $ \sA$, $ \sB$, $\sC$, and $ \sD$:
\begin{equation} \label{def:DisABCD}
\sA = \sW^{p \times p }, \quad \sB = \sW^{p \times q} \quad
\sC= \sW^{q \times p }, \quad \sD = \sW^{q \times q} .
\end{equation}
The spaces $ \sA$, $ \sB$, $\sC$, $ \sD $ admit decompositions as in \eqref{decomp1}
and \eqref{decomp2} using
\begin{align*}
& \sA_+^0 = \sW_{+,0}^{ p \times p} ,  \quad \sA_-^0 = \sW_{-,0}^{ p \times p},  \quad
\sA_d = \{\eta e_p \mid \eta \in \BC^{p \times p} \}, \\
& \sB_+ = \sW_{+}^{ p \times q},  \quad \sB_- = \sW_{-,0}^{ p \times q},  \\
& \sC_- = \sW_{-}^{q \times p },  \quad \sC_+ = \sW_{+,0}^{q \times p },  \\
& \sD_+^0 = \sW_{+,0}^{q \times q},  \quad \sD_-^0 = \sW_{-,0}^{q \times q},  \quad
\sD_d = \{ \zeta e_q \mid \zeta \in \BC^{q \times q} \}.
\end{align*}
The algebraic structure is given by the algebraic structure of the Wiener algebra and by the
matrices with entries from the Wiener algebra. Note that
\begin{equation}\label{Fabcd2T}
\a\in \sA_+  , \quad \b\in \sB_+ , \quad \c\in \sC_- , \quad \d\in \sD_- ,
\end{equation}
and we are interested (cf., \eqref{inclu12} and \eqref{inclu34}) in finding  $g\in \sB_+$ such that
\begin{align}
 \a + g \c - e_p \in  \sA_-^0    &\ands g^\ast \a  + \c \in   \sC_+ , \label{incluD12T} \\
 g \d + \b \in \sB_-^0&\ands  \d+ g^\ast \b - e_q \in \sD_+^0. \label{incluD34T}
\end{align}
Furthermore, $ a_0 = P_{\sA_d} \a $ and $ d_0 = P_{\sA_d} \d $ are invertible in $ \sA_d $ and $ \sD_d $, respectively.

In the present context the spaces $\sX$ and $\sY$, $\sX_+$ and $\sY_+$, and $\sX_-$
and $\sY_-$  defined in the first paragraph of Section \ref{sec:preliminaries} are given  by
\begin{align*}
\sX= \sA \dot{+}\sB   = \sW^{p \times p} \dot{+}  \sW^{p \times q},&\quad \sY= \sC \dot{+} \sD   = \sW^{q \times p} \dot{+}  \sW^{q \times q}, \\[.1cm]
\sX_+ = \sA_+ \dot{+}\sB_+  = \sW_+^{p \times p} \dot{+}  \sW_+^{p \times q}, &\quad\sY_+ =\sC_+ \dot{+} \sD_+^0  = \sW_{+,0}^{q \times p} \dot{+} \sW_{+,0}^{q \times q}, \\[.1cm]
\sX_- = \sA_-^0 \dot{+} \sB_- = \sW_{-,0}^{p \times p} \dot{+} \sW_{-,0}^{p \times q},&\quad \sY_- = \sC_- \dot{+} \sD_- = \sW_-^{q \times p} \dot{+}  \sW_-^{q \times q}.
\end{align*}

\begin{rem}\label{rem:identspaces}
Note that the space $ \sX = \sA \dot{+} \sB =\sW^{p \times p} \dot{+}  \sW^{p \times q}$ can be identified in a canonical way with the space $ \sW^{p \times (p+q)} $,  and analogously the subspaces $ \sX_\pm $ can be identified in a canonical way with subspaces of $ \sW^{p \times (p+q)} $. For instance,  $ \sX_+  $ with $ \sW_+^{p \times (p+q)} $.  Similarly,  $ \sY = \sC \dot{+} \sD $ can be identified with $ \sW^{q \times (p+q)} $,  and the spaces $ \sY_\pm $ with subspaces of $ \sW^{q \times (p+q)} $. We will use these identifications in the proof of Theorem \ref{thm:mainthmDS1}
\end{rem}

\begin{rem}\label{rem:HRij}
Let $g\in \sW_+^{p\ts q}$, and let $R_{ij}$, $1\leq i,j\leq 2$, be the operators defined by \eqref{defR11D}--\eqref{defR22D}. Note that the operators $H_{+,g} $, $H_{-,g^\ast}$ and $R_{ij}$ act on vector spaces $ \sW_{\pm}^m $ with $ m = p $ or $ m =q $; see the final paragraphs of the present section. As usual we extend the action of these operators to spaces of matrices of the type $\sW_{\pm}^{m \times k }$. In this way (using the preceding remark) we see that  the operators $ H_{+,g} $, $ H_{-,g^\ast} $ and $ R_{ij} $ can be identified with the operators $ \bH_{+,g} $, $ \bH_{-,g^\ast} $ and $ \bR_{ij} $ as defined in Section~\ref{sec:prelim}, respectively.
\end{rem}

\noindent\textbf{Proof  of Theorem \ref{thm:mainthmDS1}.} We will apply Theorem \ref{thm:mainthm} and Lemma \ref{isHankel} using $\sA, \sB, \sC, \sD$  in \eqref{def:DisABCD}.  First we check that   the conditions in Theorem \ref{thm:mainthm} are satisfied. Condition (a) is satisfied by assumption.

Now assume that there exists a solution $g$ to the twofold EG inverse problem. Then conditions (C1)--(C6) are  satisfied too. Next put
\begin{align*}
\Omega&= \begin{bmatrix} I_{\sW_+^p} & H_{+,g} \\[.1cm] H_{-,g^*} & I_{\sW_-^q} \end{bmatrix}:
\begin{bmatrix}\sW_+^p \\[.1cm] \sW_-^q\end{bmatrix}\to \begin{bmatrix}\sW_+^p \\[.1cm] \sW_-^q\end{bmatrix},\\
\boldsymbol\Omega &= \begin{bmatrix} I_{\sX_+} & \bH_{+,g} \\[.1cm] \bH_{-,g^*} & \bI_{\sY_-} \end{bmatrix}:
\begin{bmatrix} \sX_+ \\[.1cm] \sY_- \end{bmatrix} \to \begin{bmatrix} \sX_+\\[.1cm] \sY_- \end{bmatrix}.
\end{align*}
Here $\sX_{\pm}$ and $\sY_{\pm}$ are the spaces defined in the paragraph preceding   Remark~\ref{rem:identspaces}.   Since  the conditions of  Theorem \ref{thm:inversion1} are satisfied, we know from  \eqref{IGGIinverse} that the operator $\boldsymbol\Omega$ is invertible and its inverse is given by $\boldsymbol\Omega^{-1}=\bR$, where
\[
\bR=\begin{bmatrix} \bR_{11} &  \bR_{12} \\ \bR_{21} &  \bR_{22} \end{bmatrix}.
\]
Thus $\boldsymbol\Omega\bR$ and $\bR\boldsymbol\Omega$ are identity operators. Using the similarity  mentioned in  Remark \ref{rem:HRij} above,  we see that $ \Omega R$ and $ R \Omega $ are also identity operators, and hence $\Omega$ is invertible.  Moreover,  the  fact that $ \bR_{11} $ and $ \bR_{22} $ are invertible implies that $ R_{11} $ and $ R_{22} $ are invertible. Finally, from
\eqref{formula-bHg}, \eqref{formula-bHgast} and \eqref{formggR}
(again using the above Remark \ref{rem:HRij}) we obtain  the identities  \begin{align}
H_{+,g}=-  R_{11}^{-1}  R_{12} =
- R_{12} R_{22}^{-1} , & \quad H_{-,g^*}=-  R_{21} R_{11}^{-1}   =
- R_{22}^{-1} R_{21};   \label{formulaDGG*} \\  \noalign{\smallskip}
  g= - R_{11}^{-1} \b, & \quad  g^*= - R_{22}^{-1}\c.  \label{formulaDgg*}
\end{align}
It follows that conditions (D1) and (D2)  are fulfilled.

    Conversely, assume that conditions (D1) and (D2) are satisfied.  Then the statements (b) and (i) in Theorem \ref{thm:mainthm} follow. To apply Lemma \ref{isHankel} we set $ \bV_{\sX, \pm } = S_{p,\pm} $ and $\bV_{\sY, \pm } = S_{q,\pm} $, where $ S_{q,\pm}$ and $S_{q,\pm}$ are the shift operators defined in the final paragraph of this section.  The intertwining of $\bV_{\sX, \pm }$ and $ \bV_{\sY, \pm } $ with the Hankel-like and Toeplitz-like operators are required for the application of Lemma \ref{isHankel} but these intertwining relations correspond  with the intertwining of $ S_{p,\pm} $ and $ S_{q,\pm} $ with the Hankel and Toeplitz operators $ H_{\pm,p}$, $ H_{\pm,q}$ and $ T_{\pm,p}$, $T_{\pm,q}$ appearing in the present section. From the final part of Lemma \ref{isHankel} we conclude that
\[
S_{*,p,+} R_{12} R_{22}^{-1} =  R_{11}^{-1} R_{12} S_{q,-}.
\]
Furthermore,   Lemma~\ref{lem:inverseR} tells us that $ R_{12} R_{22}^{-1} = R_{11}^{-1} R_{12} $. Hence
\[
S_{*,p,+} R_{11}^{-1} R_{12} =  R_{11}^{-1} R_{12} S_{q,-}.
\]
But then, using   Lemma \ref{lem:Hankel} below, it follows  that there exists a $ g \in \sW_+^{p \times q} $ such that $ H_{+,g} = - R_{11}^{-1} R_{12}  $.
Similarly we obtain that there exists a $ h \in \sW_-^{q \times p } $ such that
$ H_{-,h} = - R_{22}^{-1} R_{21} $. From Lemma \ref{lem:unitsR} it then follows that $ g = - R_{11}^{-1} \b $ and $ h = - R_{22}^{-1} \c $.

It remains to show that $ h = g^\ast $. To do this, we extend some of the operators from (subspaces of) Wiener spaces to subspaces of $L^2$-spaces over the unit circle $\BT$. More specifically, for $m=p,q$ write $L^2(\BT)^m$ for the space of vectors of size $m$ whose entries are $L^2$-functions over $\BT$, and write $L^2_+(\BT)^m$ and $L^2_-(\BT)^m$ for the subspaces of $L^2(\BT)^m$ consisting of functions in $L^2(\BT)^m$ such that the Fourier coefficients with strictly negative ($-1,-2,\ldots$) coefficients and positive ($0,1,2,\ldots$) coefficients, respectively, are zero. Then, with some abuse of notation, we extend the operators $H_{+,g}$, $H_{-,h}$ and $R_{ij}$, $i,j=1,2$, in the following way:
\begin{align*}
H_{+,g}: L^2_-(\BT)^q\to L^2_+(\BT)^p, & \quad H_{-,h}: L^2_+(\BT)^p\to L^2_-(\BT)^q,\\
R_{11}: L^2_+(\BT)^p\to L^2_+(\BT)^p, & \quad R_{12}: L^2_-(\BT)^q\to L^2_+(\BT)^p,\\
R_{21}: L^2_+(\BT)^p\to L^2_-(\BT)^q, & \quad R_{22}: L^2_-(\BT)^q\to L^2_-(\BT)^q.
\end{align*}
It then follows from the representations \eqref{defR11D}--\eqref{defR22D} and
\eqref{formulaR11}--\eqref{formulaR22}
that $ R_{11} = R_{11}^\ast $, $ R_{22} = R_{22}^\ast $, and $ R_{12}^\ast = R_{21} $.
We find that
$ -H_{+,g}= R_{11}^{-1} R_{12} = \left( R_{22}^{-1} R_{21} \right)^\ast = - H_{-,h} $,
and therefore $ h = g^\ast$.
We conclude that the solution of the twofold EG inverse problem is indeed given by
\eqref{unique1}.
\epr

\medskip
\paragraph{Toeplitz and Hankel operators.}
Throughout, for a function $\rho\in\sW^{n\times m}$, we write $M(\rho)$ for the {\em multiplication operator} of $\rho$ from $\sW^m$ into $\sW^n$, that is,
\[
M(\rho):\sW^m\to \sW^n,\quad  (M(\rho) f)(e^{it})= \rho(e^{it})f(e^{it})\quad (f\in\sW_m,\, t\in[0,2\pi]).
\]
We define Toeplitz operators $ T_{\pm,\rho} $ and Hankel operators $ H_{\pm,\rho} $ as compressions of multiplication operators, as follows. {\it Fix the dimensions $ p \geq 1 $ and $ q \geq 1$ for the remaining part  of this section.}

\smallskip\noindent
If $ \rho \in \sW^{p \times p }$, then
\begin{align*}
& T_{+,\rho} = P_{+,p} M( \rho ) : \sW_+^p \to \sW_+^p ,
\quad T_{-,\rho} = ( I-P_{+,p} ) M( \rho )  : \sW_{-,0}^p \to \sW_{-,0}^p , \\
& H_{+,\rho} = P_{+,p} M( \rho ) : \sW_{-,0}^p \to \sW_+^p,
\quad H_{-,\rho} = (I-P_{+,p})  M( \rho ) : \sW_{+}^p \to \sW_{-,0}^p .
\end{align*}
If $ \rho \in \sW^{p \times q }$, then
\begin{align*}
& T_{+,\rho} = P_{+,p} M( \rho ) : \sW_{+,0}^q \to \sW_+^p ,
\quad T_{-,\rho} = ( I-P_{+,p} ) M( \rho )  : \sW_{-}^q \to \sW_{-,0}^p , \\
& H_{+,\rho} = P_{+,p} M( \rho ) : \sW_{-}^q \to \sW_+^p,
\quad H_{-,\rho} = (I-P_{+,p}) M( \rho ) : \sW_{+,0}^q \to \sW_{-,0}^p .
\end{align*}
If $ \rho \in \sW^{q \times p }$, then
\begin{align*}
& T_{+,\rho} = (I - P_{-,q} ) M( \rho ) : \sW_{+}^p \to \sW_{+,0}^q,
\quad T_{-,\rho} = P_{-,q} M( \rho )  : \sW_{-,0}^p \to \sW_{-}^q , \\
& H_{+,\rho} = (I - P_{-,q} ) M( \rho ) : \sW_{-,0}^p \to \sW_{+,0}^q,
\quad H_{-,\rho} =P_{-,q}  M( \rho ) : \sW_{+}^p \to \sW_{-}^q .
\end{align*}
and for  $ \rho \in \sW^{q \times p } $ then
\begin{align*}
& T_{+,\rho} = (I - P_{-,q} ) M( \rho ) : \sW_{+,0}^q \to \sW_{+,0}^q,
\quad T_{-,\rho} = P_{-,q} M( \rho )  : \sW_{-}^q \to \sW_{-}^q , \\
& H_{+,\rho} = (I - P_{-,q} ) M( \rho ) : \sW_{-}^q \to \sW_{+,0}^q,
\quad H_{-,\rho} = P_{-,q}  M( \rho ) : \sW_{+,0}^q \to \sW_{-}^q .
\end{align*}

\smallskip
\paragraph{Shift operators} We also define the shift operators used in the present section.
Let $ \varphi \in \sW^{p \times p } $ and $ \psi \in \sW^{q \times q} $ be defined by
$ \varphi(z) = z e_p $ and $ \psi(z) = z e_q $,  with $ e_p $ and $ e_q $ the constant functions
equal to the unity matrix.
The shift operators that we need are now defined by
\begin{align*}
S_{p,+} = M(\varphi) : \sW_+^p \to \sW_+^p , &
\quad S_{q,+} = M(\psi ) : \sW_{+,0}^q \to  \sW_{+,0}^q , \\
S_{q,-} = M(\psi^{-1} ) : \sW_-^q \to \sW_-^q , &
\quad S_{p,-} = M(\varphi^{-1} ) : \sW_{-,0}^p \to  \sW_{-,0}^p , \\
S_{*,p,+} = T_{+,\varphi^{-1}} : \sW_+^p \to \sW_+^p , &
\quad S_{*,q,+} = T_{+,\psi^{-1}}  : \sW_{+,0}^q \to  \sW_{+,0}^q , \\
S_{*,q,-} = T_{-,\psi} : \sW_-^q \to \sW_-^q , &
\quad S_{*,p,-} = T_{-, \varphi} : \sW_{-,0}^p \to  \sW_{-,0}^p ,
\end{align*}
Then
\begin{align*}
& S_{*,p,+} S_{p,+} = I_{\sW_+^p} , \quad S_{*,p,-} S_{p,-} = I_{\sW_{-,0}^p} , \\
& S_{*,q,-} S_{q,-} = I_{\sW_-^q} , \quad S_{*,q,+} S_{q,+} = I_{\sW_{+,0}^q} .
\end{align*}
Also we have for $ m \in \{ p,q \}$ and $ n \in \{ p,q \} $, and $ \rho \in \sW^{n \times m} $ that
\begin{equation*}
H_{+,\rho} S_{m,-} = S_{*,n,+} H_{+,\rho} , \quad H_{-,\rho} S_{m,+} = S_{*,n,-} H_{-,\rho},
\end{equation*}
Finally for $ \rho \in \sW_+^{ n \times m} $ we have $ T_{+,\rho} S_{m,+} = S_{n,+} T_{+,\rho}$
and if $ \rho \in \sW_-^{ n \times m} $ then $ T_{-,\rho} S_{m,-} = S_{n,-} T_{-,\rho} $.

The following result is classical and is easy to prove using the inverse Fourier transform (see, e.g., \cite[Section 2.3]{BS06} or Sections XXII -- XXIV in \cite{GGK93}).
\begin{lem}\label{lem:Hankel}  Let $ G : \sW_-^q \to \sW_+^p $, and assume that $ G S_{q,-} = S_{*,p,+} G $. Then there exists a function $ g \in \sW_+^{p \times q} $ such that $ G = H_{+,g} $.
Similarly,  if $ H: \sW_+^p \to \sW_-^q $ and   $ H S_{p,+} = S_{*,q,-} H $, then
there exists a function $ h \in \sW_-^{q \times p} $ such that $ H = H_{-,h} $.
\end{lem}

\setcounter{equation}{0}
\appendix
\renewcommand{\theequation}{A.\arabic{equation}}

\section{Hankel and Wiener-Hopf integral operators}\label{appendix}

In this appendix, which consists of three subsections, we present a number of results  that play an essential role in the proof of  Theorem \ref{thm:mainL2}.  In Subsection~\ref{ss:prelim} we recall the definition of a Hankel operator  on  $L^2(\BR_+)$ and review some basic facts. In Subsection~\ref{ss:Hankel-L1} we present  a theorem (partially new) characterising   classical Hankel integral operators  mapping  $L^1(\BR_+)^p$ into $L^1(\BR_+)^q$. Two auxilarly   results are presented  in the final subsection.

\subsection{Preliminaries about Hankel operators}\label{ss:prelim}
We begin with some preliminaries about Hankel operators on  $L^2(\BR_+)$, mainly taken from  or  Section 1.8 in \cite{Peller03} or Section 9.1 in \cite{BS06}. Throughout  $J$ is the  flip over operator  on $L^2(\BR)$  defined by $(Jf)(t)=f(-t)$. Furthermore, $\sF$ denotes  the Fourier transform on $L^2(\BR)$  defined by
\[
(\sF f)(\l)=\frac{1}{\sqrt{2\pi}}\int_{\BR} e^{i \l t}f(t)\,dt.
\]
It is well-known (see, e.g.,  \cite[Section 9.1, page 482]{BS06}) that   $\sF$ is a unitary operator and
\[
\sF^*=\sF^{-1}=J\sF \ands J\sF=\sF J.
\]
 Given $\a\in L^\iy (\BR)$  we define the multiplier $m(\a)$ and the convolution operator $M(\a)$ defined by $\a$ to be the operators on $L^2(\BR)$ given by
\[
(m(\a) f)(t)= \a(t)f(t), \quad f\in L^2(\BR),\, t\in\BR, \ands
M(\a)=\sF^{-1}m(\a)\sF.
\]
Given $f\in L^2(\BR)$ we have
\begin{align}
(M(\a)f)(t)&=\left(\sF^{-1}m(\a)\sF f\right)(t) =\left(J \sF m(\a)\sF f\right)(t)\nn\\
&=\frac{1}{\sqrt{2\pi}}\int_\BR e^{-its}\left(m(\a)\sF f\right)(s)\,ds
=\frac{1}{\sqrt{2\pi}}\int_\BR e^{-its}\a(s)(\sF f)(s)\,ds\nn\\
&=\frac{1}{2\pi}\int_\BR e^{-its}\a(s)\left(\int_\BR e^{isr}f(r)\,dr\right)\,ds \nn\\
&=\frac{1}{2\pi}\int_\BR \int_\BR e^{is(r-t)}\a(s)f(r)\,dr\,ds, \quad t\in \BR.\label{eq:multiplier1}
\end{align}
By $P$ and $Q$ we denote the orthogonal projections on $L^2(\BR)$ of which the ranges $L^2(\BR_+) $ and $L^2(\BR_-)$, respectively.

\begin{DEF}\label{def:Hankel1} Let   $\a\in L^\iy (\BR)$.   Then the  \emph{Hankel operator} defined by $\a$ is the operator on $L^2(\BR_+)$  given by
\[
H(\a)=PM(\a)J|_{L^2(\BR_+) }: L^2(\BR_+)\to L^2(\BR_+).
\]
\end{DEF}
\medskip
The action of  the Hankel operator $H(\a)$ on $f\in L^2(\BR_+)$ is given by
\begin{align}
\left(H(\a)f\right)(t)&=\left(PM(\a)Jf\right)(t)=\left(M(\a)Jf\right)(t)\nn \\
&=\frac{1}{2\pi}\int_\BR \int_\BR e^{is(r-t)}\a(s)f(-r)\,dr\,ds \nn \\
&=\frac{1}{2\pi}\int_\BR \int_\BR e^{-is(t+r)}\a(s)f(r)\,dr\,ds \nn \\
&=\frac{1}{2\pi} \int_\BR \int_{0}^\infty e^{-is(t+r)} \a(s)f(r)\,dr\,ds, \quad t\geq 0. \label{eq:Hankel1}
\end{align}

The following result provides a characterization of which operators on $L^2(\BR_+)$ are Hankel operators; see \cite[Exercise (a) on page 199-200]{Nikolski17}.

\begin{thm}\label{thm:char-Hankel1}
A bounded linear operator $ K$ on  $L^2 (\BR_+) $ is a Hankel operator  if and only if
$ V_\tau^\ast K = K V_\tau $ for all $ \tau \geq 0 $, where for each $\t\geq 0$  the operator $V_\t$  is the transition operator on  $L^2 (\BR_+) $ defined by
\begin{equation}
( V_\tau f )(t) = \left\{ \begin{matrix} f(t-\tau), & \quad t \geq \tau, \\ 0, &\quad 0 \leq t \leq \tau, \end{matrix} \right. \label{defVtau1}
\end{equation}
\end{thm}

\begin{rem}\label{rem:story}
When we worked on this paper we assumed that the above theorem, which is a natural analogue of the intertwining shift relation theorem for discrete Hankel operators, to be true and that we only had to find a reference. The latter turned out to be a bit difficult. Various Hankel operator experts told us ``of course, the result is true.'' But no reference. We asked Vladimir Peller, and he mailed us how the result could be proved using the beautiful relations between $H^2$ on the disc and $H^2$ on the upper half plane, but again no reference. What to do? Should we include Peller's proof? November last year Albrecht B\"otcher solved the problem. He referred us to Nikolski's book \cite{Nikolski17} which appeared recently in Spring 2017 and contains the result as an exercise. Other references remain welcome.
\end{rem}

Next we consider the special case  when the defining function $\a$ is given by
\begin{equation} \label{eq:Wiener1}
\a(\l)=\int_\BR e^{i \l s} a(s)\,ds, \quad \mbox{where $a\in L^1(\BR_+)$}.\
\end{equation}
Then
\begin{equation}\label{eq:Hankel2}
(H(\a)f)(t)= \int_0^\iy a(t+s) f(s)\, ds, \quad  t\in \BR_+,  \  f\in L^2(\BR_+).
\end{equation}
In this case one calls $H(\a)$  the  \emph{classical Hankel integral operator} defined by $a$. To prove \eqref{eq:Hankel2} we may without loss of generality assume that  $f$ belongs to   $L^1(\BR_+)\cap L^2(\BR_+)$ and $\a$ is rational. In that case  using \eqref{eq:Wiener1} we have
\begin{align*}
\left(H(\a)f\right)(t)&=\frac{1}{2\pi}\int_\BR \left(\int_0^\iy e^{-is(t+r)} \a(s)f(r)\,dr\right)\,ds \\
&=\frac{1}{2\pi} \int_0^\iy\left(\int_\BR e^{-is(t+r)}\a(s)\, ds\right)f(r)\,dr\\
&=\int_0^\iy a(t+r)f(r)\,dr, \quad t\geq 0.
\end{align*}
If  $\a$ is given by \eqref{eq:Wiener1}, then $\a$ belongs to the Wiener algebra over $\BR$ and thus $H(\a)$ also defines a bounded linear operator on $L^1(\BR_+)$.

 We shall also deal with Hankel operators defined by  matrix-valued functions. Let $\a$ be a $q\ts p$ matrix whose entries $\a_{ij}$, $1\leq i\leq q$, $1\leq j\leq p$, are  $L^\iy(\BR_+)$ functions.  Then $H(\a)$  will denote the Hankel operator   from $L^2(\BR_+)^p$ to $L^2(\BR_+)^q$ defined by
\begin{equation}\label{def:block-Hankel}
H(\a)=\begin{bmatrix}H(\a_{11})&\cdots& H(\a_{1p})\\ \vdots&\cdots&\vdots\\
H(\a_{q1})&\cdots& H(\a_{qp})\end{bmatrix}.
\end{equation}
If the operators $H(\a_{ij})$, $1\leq i\leq q$, $1\leq j\leq p$, are all classical Hankel integral operators, then we call  $H(\a)$ a classical Hankel integral operator too.

\subsection{Classical Hankel integral operators on $L^1$ spaces}\label{ss:Hankel-L1}
The main  theorem of this section allows us to identify the classical Hankel integral operators among all operators  from  $L^1(\BR_+)^p$ to $L^1(\BR_+)^q$. We begin with some preliminaries about related   Sobolev spaces.

Let $n$ be a positive integer. By $\sob(\BR_+)^n$ we denote the Sobolev space consisting  of  all functions $\va \in L^1(\BR_+)^n$ such that $\va$  is absolutely continuous on compact intervals of $\BR_+$  and    $\va^\prime\in L^1(\BR_+)^n$.  Note that
\begin{equation}\label{eq:sob-prop}
\va \in \sob(\BR_+)^n \quad \Longrightarrow \quad \va(t)=-\int_t^\iy  \va^\prime (s) \,ds, \quad t\geq 0.
\end{equation}
The linear space $\sob(\BR_+)^n$ is a Banach space with norm
\begin{equation}\label{def:normSB}
\|\va\|_{\sob}=\|\va\|_{{L^1}}+\|\va^\prime\|_{{L^1}}.
\end{equation}
Furthermore,    $\sob(\BR_+)^n$    is continuously and densely embedded in $L^1(\BR_+)^n$. More precisely, the map $j:\sob(\BR_+)^n\to L^1(\BR_+)^n$ defined by $j\va=\va$ is a continuous linear map which is one-to-one and has dense range. From \eqref{def:normSB} we see that $j$ is a contraction.  We are now ready to state and proof the main theorem of this section.

\begin{thm}\label{thm:main1} An operator $K$ from $L^1(\BR_+)^p$ to  $L^1(\BR_+)^q$  is a classical Hankel integral operator  if and  only if  the following two conditions are satisfied:
\begin{itemize}
\item[\textup{(H1)}]
$K$ maps $\sob(\BR_+)^p$ boundedly into $\sob(\BR_+)^q$;

\item[\textup{(H2)}] there exists $k\in L^1(\BR_+)^{q\ts p}$ such that  $(K\va)^\prime+K\va^\prime =k(\cdot)\va(0)$ for each $\va \in \sob(\BR_+)^p$.
\end{itemize}
Moreover, in that case the operator $K$ is given by
\begin{equation}\label{def:hankel-cl}
(Kf)(t)=\int_0^\iy k(t+s)f(s)\, ds, \quad 0\leq t<\iy,
\end{equation}
where  $k\in L^1(\BR_+)^{q\ts p}$ is the matrix function from \textup{(H2)}.
\end{thm}
\bpr   We split the proof into three parts. In the first part we show that the conditions (H1) and (H2) are necessary. The proof  is taken from \cite{EG03}, and is given here for the sake of completeness. The second and third part concern  the reverse implication which seems to be new. In the second part we assume that $p=q=1$, and in the third part  $p$ and $q$ are arbitrary positive integers.

\smallskip
\noindent\textsc{Part 1.} Let $K$ on $L^1(\BR_+)$  be a classical Hankel integral operator, and assume  $K$ is given by \eqref{def:hankel-cl} with $k\in L^1(\BR_+)$.
Let $\va\in \sob(\BR_+)$. Then
\begin{equation}\label{eq:K-sob0}
(K\va)(t)=\int_0^\iy k(t+s)\va(s)\, ds =\int_t^\iy k(s) \va(s-t)\,ds.
\end{equation}
It follows that
\begin{align}
\left(\frac{d}{dt}K\va\right)(t)&=-\int_t^\iy k(s) \va^\prime(s-t)\,ds+k(t)\va(0) \nn \\
&=-\int_0^\iy k(t+s)\va^\prime(s)\, ds+k(t)\va(0)\label{eq:K-sob2}\\
&=-\left(K\frac{d}{dt}\va\right)(t)+k(t)\va(0). \label{eq:K-sob3}
\end{align}
This proofs (H2).
From the first identity in \eqref{eq:K-sob0} it follows that $K\va$ belongs to $L^1(\BR_+)$.
Since $\va^\prime\in L^1(\BR_+)$, we have $(K\va)^\prime= K\va' + k \va(0)\in L^1(\BR_+)$ and it follows that
$K\va\in \sob(\BR_+)$.
We conclude  that $K$ maps $\sob(\BR_+)$ into $\sob(\BR_+)$. Furthermore, from   \eqref{eq:K-sob0} we see that
\begin{align*}
\|K\va\|_{L^1} \leq \|k\|_{L^1} \|\va\|_{L^1}\leq \|k\|_{L^1}\|\va\|_{\sob}.
\end{align*}
From \eqref{eq:K-sob2} (using $\va(0)=-\int_0^\iy \va^\prime(s)\,ds$) it follows  that
\begin{align*}
\|(K\va)^\prime\|_{L^1}& \leq \|k\|_{L^1} \|\va^\prime\|_{L^1}+\|k\|_{L^1}\|\va^\prime\|_{L^1}\\[.2cm]
&\leq 2\|k\|_{L^1}\|\va^\prime\|_{L^1}\leq 2\|k\|_{L^1}\|\va\|_\sob.
\end{align*}
Hence $\|K\va\|_\sob\leq 3\|k\|_{L^1}\|\va\|_\sob$. Thus $K|_{\sob(\BR_+)}$ is a bounded  operator on $\sob(\BR_+)$, and   item (H1) is proved.

\smallskip
\noindent\textsc{Part 2.} In this part  $p=q=1$, and we assume that items (H1) and (H2) are satisfied. Given $k$ in item (H2), let $H$ be the operator on $L^1(\BR_+)$ defined by
\[
 (Hf)(t)=\int_0^\iy k(t+s)f(s)\, ds, \quad 0\leq t<\iy.
 \]
Then  $H$ is a classical Hankel integral operator, and the first part of the proof  tells us that $(H\va)^\prime+H\va^\prime =k(\cdot)\va(0)$ for each $\va \in \sob(\BR_+)$. Now put $M=K-H$. Then $M$ is an operator  on $L^1(\BR_+)$, and $M$ maps $ \sob(\BR_+)$  into $ \sob(\BR_+)$. Furthermore, we have
\begin{equation}\label{eq:propM1}
(M\va)^\prime=-M\va^\prime,\quad \va \in \sob(\BR_+).
\end{equation}
It suffices to prove that $M$ is zero.

For $n=0,1,2, \ldots$ let $\va_n$ be the function on $\BR_+$ defined by  $\va_n(t)=t^n e^{-t}$, $0\leq t <\iy$. Obviously, $\va_n\in \sob(\BR_+)$. By induction we shall prove  that $M\va_n$  is zero for each  $n=0,1,2, \ldots$.  First we show that $M\va_0=0$. To do this note that $\va_0^\prime (t)=-e^{-t}=-\va_0 (t)$. Using \eqref{eq:propM1} it follows  that $\psi_0:=M\va_0$ satisfies
\[
\psi_0^\prime=(M\va_0)^\prime=-M\va_0^\prime=M\va_0=\psi_0.
\]
Thus $\psi_0$ satisfies the differential equation $\psi_0^\prime=\psi_0$, and hence $\psi_0(t)=c e^t$  on $[0, \iy)$ for some $c\in \BC$. On the other hand,  $\psi_0=M\va_0\in \sob(\BR_+)\subset L^1(\BR_+)$. But then $c$ must be zero, and we conclude that $M\va_0=0$.

Next, fix   a positive   integer  $n\geq 1$, and assume  that $M\va_j=0$ for  $j=0, \ldots, n-1$. Again we use \eqref{eq:propM1}. Since
\[
\va_n^\prime(t)=n t^{n-1}e_{-t}- t^n e^{-t}=n\va_{n-1}-\va_n,
\]
we obtain
\[
(M\va_n)^\prime=-M\va_n^\prime=n M\va_{n-1}+M\va_n.
\]
But, by assumption, $M\va_{n-1}=0$. Thus $(M\va_n)^\prime=M\va_n$, and hence $\psi_n:=M\va_n$ satisfies the differential equation $\psi_n^\prime=\psi_n$. It follows that $\psi_n(t)=ce^t$ on $[0, \iy)$ for some $c\in \BC$. On the other hand,  $\psi_n=M\va_n\in \sob(\BR_+)\subset L^1(\BR_+)$. But then $c=0$, and we conclude that $M\va_n=0$.

By induction we obtain $M\va_j=0$ for each $j=0,1,2, \ldots$. But then $Mf=0$ for any $f$ of the form $f(t)=p(t)e^{-t}$, where $p$ is a polynomial.  The set of all these functions is dense in $L^1(\BR_+)$. Since $M$ is an operator on $L^1(\BR_+)$, we conclude that $M=0$.

\smallskip
\noindent\textsc{Part 3.}  In this part  we use the result of the previous part to prove the  sufficiency  of the conditions  (H1) and (H2). Assume $K$ from $L^1(\BR_+)^p$ to  $L^1(\BR_+)^q$, and write $K$ as a $q\times p$ operator matrix
\begin{equation}\label{eq:matKij}
K=\begin{bmatrix}K_{11}&\cdots& K_{1p}\\ \vdots&\cdots&\vdots\\
K_{q1}&\cdots& K_{qp}\end{bmatrix},
\end{equation}
where $K_{ij}$ is an operator  on $L^1(\BR_+)$ for    $1\leq j\leq p$ and $1\leq i\leq q$. let
\begin{align*}
&\t_j :L^1(\BR_+)\to L^1(\BR_+)^p, \quad \t_j f=[\delta_{j,k} f]_{k=1}^p\quad (f\in L^1(\BR_+);\\
&\pi_i  : L^1(\BR_+)^q\to L^1(\BR_+), \quad \pi_j f= f_j, \quad (f=\begin{bmatrix}f_1\\ \vdots\\ f_q  \end{bmatrix}\in  L^1(\BR_+)^q).
\end{align*}
Note that $K_{ij}= \pi_i K\t_j$ for   each $i, j$. Furthermore, we have
\[
\t_j\sob(\BR_+)\subset \sob(\BR_+)^p\ands \pi_j \sob(\BR_+)^q\subset \sob(\BR_+)
\]
Now fix a pair $i, j$,   $1\leq j\leq p$ and $1\leq i\leq q$. Then conditions (H1) and (H2) tell us  that
\begin{itemize}
\item[\textup{(i)}] $K_{ij}$ maps $\sob(\BR_+)$ into $\sob(\BR_+)$ and
$K_{ij}|_{\sob(\BR_+)}$ is a bounded operator  from $\sob(\BR_+)$ to $\sob(\BR_+)$;
\item[\textup{(ii)}] there exists $k_{ij}\in L^1(\BR_+) $ such that  $(K_{ij}\va)^\prime+K_{ij}\va^\prime =k_{ij}(\cdot)\va(0)$ for each $\va \in \sob(\BR_+)$.
\end{itemize}
But then we can use the result of the second part of the proof which covers the case when $p=q=1$.  It follows that $K_{ij}$ is a classical Hankel integral operator. Moreover, $K_{ij}$ is given by
\[
(K_{ij}f)(t)=\int_0^\iy  k_{ij}(t+s)f(s)\,ds, \quad 0\leq t <\iy \ands f\in  L^1(\BR_+).
\]
Here $k_{ij}\in L^1(\BR_+)$ is the function appearing in item (ii) above.  Recall that $K$ is given by \eqref{eq:matKij}. Since the pair $i,j$ is arbitrary, we see that   $K$ is a classical Hankel integral operator, and
\[
(Kf)(t)=\int_0^\iy k(t+s)f(s)\, ds, \quad 0\leq t<\iy \ands  f\in  L^1(\BR_+)^p,
\]
where
\[
k:=\begin{bmatrix}k_{11}&\cdots& k_{1p}\\ \vdots&\cdots&\vdots\\
k_{q1}&\cdots& k_{qp}\end{bmatrix}\in L^1(\BR_+)^{q\ts p}.
\]
This completes the proof.  \epr

\medskip
The following corollary shows that if the operator $ K $ in Theorem \ref{thm:main1} is assumed to be a Hankel operator, i.e. $ K = H(\a) $ for some  $\a \in L^\infty(\BR)^{q \times p} $, then it suffices to verify (H1) to conclude that $ H(\a) $ is a classical Hankel operator.

\begin{cor}\label{cor:Hankel-int1}Let $\a\in L^\iy (\BR)^{q\ts p}$, and  assume that  $H(\a)$ maps $L^1(\BR_+)^p$ into $L^1(\BR_+)^q$. Furthermore, assume  that  $K=H(\a)$ satisfies condition \textup{(H1)} in Theorem \ref{thm:main1}, i.e.,  $H(\a)$ maps $\sob(\BR_+)^p$ into $\sob(\BR_+)^q$ and  the operator $H(\a)|_{\sob(\BR_+)^p}$ is a bounded operator from $\sob(\BR_+)^p$ into $\sob(\BR_+)^q$. Then $K=H(\a)$ also satisfies condition \textup{(H2)} in Theorem \ref{thm:main1}, and thus there exists $k\in L^1(\BR_+)^{q\ts p}$ such that
\begin{equation}\label{eq:kdefinesal}
\a(\l)= \int_\BR e^{i\l s} k(s)\,ds, \quad \l \in \BR.
\end{equation}
In particular,  $H(\a)$ is a classical Hankel integral operator.
\end{cor}
\bpr  We split the proof into two parts. In the first part we assume  that $p=q=1$. In the second part $p$ and $q$ are arbitrary positive integers, and we reduce the problem  to the case considered in the first part.

\smallskip
\noindent\textsc{Part 1.} In this part we prove the theorem for the case when $p=q=1$.
Note that we assume that condition (H1)  in Theorem \ref{thm:main1} is satisfied for $H(\a)$ in place of $K$.  Take $\va\in \sob(\BR_+)$.   Since $H(\a)$  maps  $\sob(\BR_+)$ into $\sob(\BR_+)$, the function  $H(\a)\va$ also belongs to $\sob(\BR_+)$, and hence $(H(\a)\va)^\prime$ belongs to $L^1(\BR_+)$.  On the other hand, since   $H(\a)$ maps  $L^1(\BR_+)$ into $L^1(\BR_+)$ and $\va ^\prime\in L^1(\BR_+)$, we also have $H(\a)\va ^\prime\in L^1(\BR_+)$.  Hence
\[
\frac{d}{dt} \big(H(\a)\va  \big)-  H(\a)\frac{d}{dt}\va \in L^1(\BR_+).
\]
Using \eqref{eq:Hankel1} we obtain
\begin{align}
\frac{d}{dt} \Big(H(\a)\va (t)\Big)&=
 \frac{1}{2\pi}\int_\BR (-is)e^{-its}\a(s)\left(\int_0^\iy  e^{-isr}\va (r)\,dr\right)\,ds \nn \\
&=  \frac{1}{2\pi}\int_\BR e^{-its}\a(s)\left(\int_0^\iy  \left(\frac{d}{dr}e^{-isr}\right)\va (r)\,dr\right)\,ds, \label{eq:integrals}
\end{align}
and
\[
\Big(H(\a)\frac{d}{dt}\va  \Big)(t)=\frac{1}{2\pi}\int_\BR e^{-its}\a(s)\left(\int_0^\iy  e^{-isr}\va ^\prime(r)\,dr\right)\,ds.
\]
Now integration by parts yields
\begin{align*}
\int_0^\iy  \Big(\frac{d}{dr}e^{-isr}\Big) \va (r)\,dr&=-\int_0^\iy  e^{-isr}\va ^\prime(r)\,dr+ \left(e^{-isr}\va (r)\Big|_0^\iy\right)\\[.1cm]
&=-\int_0^\iy  e^{-isr}\va ^\prime(r)\,dr-\va (0).
\end{align*}
Hence
\[
\int_0^\iy  \Big(\frac{d}{dr}e^{-isr}\Big) \va (r)\,dr +\int_0^\iy  e^{-isr}\va ^\prime(r)=-\va (0).
\]
This implies that
\begin{equation}\label{indepPhi}
-\frac{\va (0)}{2\pi}\int_{\BR} e^{-its}\a(s)\, ds
= \left(\frac{d}{dt} \big(H(\a)\va  \big) +  H(\a)\frac{d}{dt}\va \right)(t),\quad t\geq 0.
\end{equation}
In case $\va (0)\not= 0$, it follows that (H2) holds for $K=H(\a)$ with $k\in L^1(\BR_+)$ given by
\[
k(t):=\frac{1}{\va(0)}\left(\frac{d}{dt} \big(H(\a)\va \big)+ H(\a)\frac{d}{dt}\va \right)(t)
=-\frac{1}{2\pi}\int_{\BR} e^{-its}\a(s)\, ds,
\]
which is independent of the choice of $\va$. On the other hand, in case $\va(0)=0$, then \eqref{indepPhi} shows that (H2) still holds for this choice of $k$. We can thus use Theorem \ref{thm:main1} with $K=H(\a)$ to conclude that $H(\a)$ is a classical Hankel integral operator and $\a$ is defined by \eqref{eq:kdefinesal}.

\smallskip
\noindent\textsc{Part 2.} In this part $p$ and $q$ are arbitrary positive integers. Since $\a\in L^\iy (\BR)^{q\ts p}$, the function $\a$ is a $q\ts p$ matrix function of which the (i,j)-th entry
$\a_{i,j}$  belongs to $L^\iy (\BR_+)$. It follows that
\[
K=H(\a)=\begin{bmatrix}H(\a_{11})&\cdots& H(\a_{1p})\\ \vdots&\cdots&\vdots\\
H(\a_{q1})&\cdots& H(\a_{qp})\end{bmatrix}.
\]
Put $K_{ij}=H(\a_{ij})$, where  $1\leq j\leq p$ and $1\leq i\leq q$.  Now  fix (i,j). Since $K$ maps    $L^1(\BR_+)^p$ into $L^1(\BR_+)^q$, the operator $K_{ij}$ maps     $L^1(\BR_+)$ into $L^1(\BR_+)$. Furthermore, since  $K=H(\a)$ satisfies condition \textup{(H1)} in Theorem \ref{thm:main1}, the operator  $K_{ij}=H(\a_{ij})$ satisfies condition \textup{(H1)} in Theorem \ref{thm:main1} with $p=q=1$. But then the result of the first part of the proof tells us that $K_{ij}=H(\a_{ij})$ satisfies condition \textup{(H2)} in Theorem \ref{thm:main1} with $p=q=1$. Thus, using \eqref{eq:kdefinesal}, there exists $k_{ij}\in L^1 (\BR_+)$ such that
\[
\a_{ij}=\int_\BR e^{i \l s} k_{ij}(s)\, ds, \quad \l \in \BR.
\]
The latter holds for each  $1\leq j\leq p$ and $1\leq i\leq q$. It follows that
\[
\a(\l)=\int_\BR e^{i \l s} k(s)\, ds,\hspace{.2cm}  \mbox{where}\hspace{.2cm}  k=\begin{bmatrix}k_{11} &\cdots& k_{1p}\\ \vdots&\cdots&\vdots\\
k_{q1}&\cdots& k_{qp}\end{bmatrix}\in L^1 (\BR_+)^{q\ts p}.
\]
This completes the proof. \epr

\subsection{Two auxiliary results}\label{ss:auxil2}
We present two lemmas   concerning condition  (H1) in Theorem \ref{thm:main1}. We begin with some preliminaries. Let
\begin{equation}\label{eq:defM12}
M=I+H_{11}H_{12}+ H_{21}H_{22},
\end{equation}
where
\begin{align*}
&H_{11}:L^1(\BR_+)^q\to L^1(\BR_+)^p, \quad  H_{12}:L^1(\BR_+)^p\to L^1(\BR_+)^q,\\[.2cm]
&H_{21}:L^1(\BR_+)^r\to L^1(\BR_+)^p, \quad  H_{22}:L^1(\BR_+)^p\to L^1(\BR_+)^r,\
\end{align*}
and  we assume that $H_{ij}$ is a classical Hankel integral operator  for each $1\leq i,j\leq 2$. We are interested in  computing the inverse of $M$, assuming the inverse  exists.  Put $\wt{M}=I+\wt{H}_1\wt{H}_2$, where
\begin{align*}
\wt{H}_1&=\begin{bmatrix}H_{11}&H_{21} \end{bmatrix}:
\begin{bmatrix}L^1(\BR_+)^q\\[.1cm] L^1(\BR_+)^r \end{bmatrix}\to  L^1(\BR_+)^p, \\
\wt{H}_2&=\begin{bmatrix}H_{12}\\H_{22}  \end{bmatrix}: L^1(\BR_+)^p\to
\begin{bmatrix}L^1(\BR_+)^q\\ L^1(\BR_+)^r \end{bmatrix}.
\end{align*}
Note that the entries of $\wt{H}_1$ and $\wt{H}_2$ are classical Hankel integral operators, and
\[
\wt{M}=I+\wt{H}_1\wt{H}_2=I+H_{11}H_{12}+ H_{21}H_{22}=M.
\]
It follows that $M$ is invertible if and only if $ \wt{M}$ is invertible, and in that case
\begin{equation}\label{eq:M-inverse1}
 \wt{M}^{-1}=M^{-1}.
\end{equation}
Theorem 0.1 in \cite{GrK05} tells us how to compute $ \wt{M}^{-1}$. This yields the following result.

\begin{lem}\label{lem:M-inv} Assume $M$ given by \eqref{eq:defM12} is invertible. Then
\begin{equation}\label{eq:Minverse2}
M^{-1}=I+K_1+K_2+K_3+K_4,
\end{equation}
where for each $j=1,2,3, 4$ the operator $K_j$ is a product of two classical Hankel integral operators.  In particular, $M^{-1}=I+K$, where $K$ is an operator on $L^1(\BR_+)^p$ satisfying condition  (H1)  in Theorem \ref{thm:main1}.
\end{lem}
\bpr From Theorem 0.1 in \cite{GrK05} we know that
\begin{equation}\label{eq:wtM-inverse}
 \wt{M}^{-1}=I+AB+CD,
\end{equation}
where the operators $A,B,C, D$ have the following operator matrix representation:
\[
A=\begin{bmatrix} A_{11}&A_{12} \end{bmatrix}, \
B=\begin{bmatrix} B_{11}\\ B_{21} \end{bmatrix}, \  C=\begin{bmatrix} C_{11}&C_{12} \end{bmatrix}, \   D=\begin{bmatrix} D_{11}\\ D_{21}\end{bmatrix}.
\]
and for each $i,j=1,2$ the entries $A_{ij}$, $B_{ij}$, $C_{ij}$, $D_{ij}$ are classical Hankel integral operators. Using \eqref{eq:M-inverse1} and  \eqref{eq:wtM-inverse} it follows that
\[
M^{-1}=\wt{M}^{-1}= I+\begin{bmatrix} A_{11}&A_{12}\end{bmatrix}
\begin{bmatrix} B_{11}\\B_{21}\end{bmatrix}+
\begin{bmatrix} C_{11}&C_{12}\end{bmatrix}
\begin{bmatrix} D_{11}\\D_{21}\end{bmatrix}.
\]
Thus \eqref{eq:Minverse2} holds true with
\[
K_1=A_{11}B_{11}, \quad K_{2}=A_{12}B_{21}, \quad K_{3}=C_{11}D_{11}, \quad K_{4}=C_{12}D_{22}.
\]
Clearly for each $j=1,2,3, 4$ the operator $K_j$ is a product of two classical Hankel integral operators. Recall that  for each classical Hankel integral operator $H$ from $L^1(\BR_+)^n$ to $L^1(\BR_+)^n$ for some $n$ and $m$ we have  $H$ maps $\sob(\BR_+)^n$ into $\sob(\BR_+)^m$ and $H|_{\sob(\BR_+)^n}$ is bounded as an operator  from $\sob(\BR_+)^n$ into $\sob(\BR_+)^m$.  It follows that the same is true if $H$ is a sum or a product   of  classical Hankel integral operators. But then condition  (H1)  in Theorem \ref{thm:main1} is satisfied
for $K=K_1+K_2+K_3+K_4$.\epr

\begin{lem}\label{lem:Wiener1} Let $\t\in  L^1(\BR_-)^{q\ts p}$, and let $T$ be the Wiener-Hopf integral operator  mapping  $L^1(\BR_+)^p$ into $L^1(\BR_+)^q$ defined by
\begin{equation}\label{eq:defT}
(Tf)(t)= \int_t^\iy \t(t-s)f(s)\, ds,  \quad  0\leq t<\iy\quad (f\in L^1(\BR_+)^p).
\end{equation}
Then $T$ maps $\sob(\BR_+)^p$ into $ \sob(\BR_+)^q$, and $T|_{\sob(\BR_+)^p}$ is a bounded linear operator from  $\sob(\BR_+)^p$ into $\sob(\BR_+)^q$.
\end{lem}
\bpr
We split the proof into two parts.  In the first part we prove the lemma for the case when $p=q=1$. In the second $p$ and $q$ are arbitrary  positive integers, and we reduce the problem  to the case considered in the first part.

\smallskip
\noindent\textsc{Part 1.} In this part we prove the lemma for the case when $p=q=1$. To do this  take $\va \in \sob(\BR_+)$. Then
\begin{align*}
(T\va)(t)&=  -\int_t^\iy \t(t-s) \left(\int_s^\iy  \va^\prime (r) \,dr\right)\, ds\\
&=-\int_t^\iy \left(\int_t^r \t (t-s)\, ds\right)\va^\prime (r)\, dr, \quad 0\leq t<\iy.
\end{align*}
Put
\[
\rho(r-t):=\int_0^{r-t}\t(-u)\,du=\int_t^r \t (t-s)\, ds, \quad 0\leq t\leq r<\iy.
\]
Note that $\rho(0)=0$. Furthermore,
\begin{equation}\label{eq:Tva}
(T\va)(t)=-\int_t^\iy \rho(r-t)\va^\prime(r)\, dr, \ands \rho^\prime(t)=\t(-t)\quad (0\leq t<\iy).
\end{equation}
Using \eqref{eq:defT} with $f=\va$ we see that $\psi:= T\va$ belongs to $L^1(\BR_+)$.  Furthermore, from the first identity in \eqref{eq:Tva} it follows that $\psi$  is absolutely continuous on compact intervals of $\BR_+$.  Using the Leibnitz rule and  the second  identity in \eqref{eq:Tva}, we obtain
\begin{align}
\psi^\prime(t)&=-\frac{d}{dt}\int_t^\iy \rho(r-t)\va^\prime(r)\, dr \nn \\
& =-\int_t^\iy \frac{\partial}{\partial t}\rho(r-t)\va^\prime(r)\, dr +\rho(t-t)\va^\prime(t)\nn \\
&=  - \int_t^\iy \t(t-r)\va^\prime(r)\, dr.  \label{eq:Wiener2}
\end{align}
Since $\tau\in L^1(\BR_-)$ and $\va^\prime\in L^1(\BR_+)$, it follows that $\psi^\prime\in L^1(\BR_+)$. We conclude that $\psi$ belongs to $\sob(\BR_+)$.

It remains to show that $T|_{\sob(\BR_+)}$ is bounded on $\sob(\BR_+)$.  As before let  $\va\in L^1(\BR_+)$, and let $\psi=T\va$. By $\|T\|$ we denote the norm of $T$ as an operator on $L^1(\BR_+)$. From  the definition of $T$ in \eqref{eq:defT} and using \eqref{def:normSB}, we see that
\[
\|\psi\|_{L^1}\leq  \|T\| \|\va\|_{L^1}\leq  \|T\| \|\va\|_{\sob}
\]
On the other hand from \eqref{eq:Wiener2} and using \eqref{def:normSB} we obtain
\[
\|\psi^\prime\|_{L^1}\leq  \|T\|  \|\va^\prime\|_{L^1}\leq  \|T\| \|\va\|_{\sob}.
\]
Together these inequalities show (using \eqref{def:normSB}) that $\|\psi\|_{\sob}\leq  \|T\| \|\va\|_{\sob}$. Thus $\|T\|_{\sob(\BR_+)}\leq \|T\|$. This proofs the lemma for the case when $p=q=1$.

\smallskip
\noindent\textsc{Part 2.} In this part $p$ and $q$ are arbitrary positive integers.  Since $\t\in  L^1(\BR_-)^{q\ts p}$, the function $\t$ is a $q\ts p$ matrix function of which the (i,j)-th entry  $\t_{ij}$  belongs to $L^1 (\BR_-)$. It follows that
\begin{equation}\label{eq:matTij}
T=\begin{bmatrix}T_{11}&\cdots& T_{1p}\\ \vdots&\cdots&\vdots\\
T_{q1}&\cdots& T_{qp}\end{bmatrix},
\end{equation}
where for    $1\leq j\leq p$ and $1\leq i\leq q$ the operator $T_{ij}$ is the Wiener-Hopf integral operator  on  $L^1(\BR_+)$ given gy
 \[
 (T_{ij})f(t)= \int_t^\iy \t_{ij}(t-s)f(s)\, ds,  \quad  0\leq t<\iy \quad (f\in L^1(\BR_+)).
 \]
From the first part of the proof  we know that for each $i,j$ the operator $T_{ij}$  maps $\sob(\BR_+)$ into $ \sob(\BR_+)$, and $T|_{\sob(\BR_+)}$ is a bounded linear operator from  $\sob(\BR_+)$ into $\sob(\BR_+)$. Now recall that $T$ is given by \eqref{eq:matTij}. It follows that $T$  maps $\sob(\BR_+)^p$ into $ \sob(\BR_+)^q$, and $T|_{\sob(\BR_+)^p}$ is a bounded linear operator from  $\sob(\BR_+)^p$ into $\sob(\BR_+)^q$, which completes the proof. \epr

\bigskip
\paragraph{\bf Acknowledgement}
This work is based on the research supported in part by the National
Research Foundation of South Africa. Any opinion, finding and conclusion or
recommendation expressed in this material is that of the authors and
the NRF does not accept any liability in this regard.

\end{document}